\documentclass[]{amsart}
 \usepackage{graphicx}
\usepackage{amsmath}
\usepackage{amssymb}
\usepackage{amsfonts}
\usepackage{amsthm}
\usepackage{epstopdf}

\newtheorem{theorem}{Theorem}

\newtheorem{conjecture}[theorem]{Conjecture}
\theoremstyle{remark}
\newtheorem{remark}{Remark}[section]

\begin{document}
\title[Zakharov-Kuznetsov equation in two dimensions]
{Numerical study of Zakharov-Kuznetsov\\ equations in two dimensions}

\author[C. Klein]{Christian Klein}
\address{Institut de Math\'ematiques de Bourgogne, UMR 5584; \\
                Universit\'e de Bourgogne-Franche-Comt\'e, 9 avenue Alain Savary, 21078 Dijon
                Cedex, France} 
\email{Christian.Klein@u-bourgogne.fr}

\author[S. Roudenko]{Svetlana Roudenko}
\address{Department of Mathematics \& Statistics\\Florida International University,  Miami, FL, 33199, USA}
\curraddr{}
\email{sroudenko@fiu.edu}
    
\author[N. Stoilov]{Nikola Stoilov}
\address{Institut de Math\'ematiques de Bourgogne, UMR 5584; \\
                Universit\'e de Bourgogne-Franche-Comt\'e, 9 avenue Alain Savary, 21078 Dijon
                Cedex, France} 
\email{Nikola.Stoilov@u-bourgogne.fr}


\begin{abstract}
We present a detailed numerical study of solutions to the (generalized) Zakharov-Kuznetsov  equation in two spatial dimensions with various power nonlinearities. In the $L^{2}$-subcritical case, numerical evidence is presented for the stability of solitons and the soliton resolution for generic initial data. In the $L^2$-critical and supercritical cases, solitons appear to be unstable against both dispersion and blow-up. It is conjectured that blow-up happens in finite time and that blow-up solutions have some resemblance of being self-similar, i.e., the blow-up core forms a rightward moving self-similar type rescaled profile 
with the blow-up happening at infinity in the critical case and at a finite location in the supercritical case. 
In the $L^{2}$-critical case, the blow-up appears to be similar to the one in the $L^{2}$-critical generalized Korteweg-de Vries  equation with the profile being a dynamically rescaled soliton. 
\end{abstract}


 

\maketitle

\section{Introduction}

We are interested in the 2D generalized Zakharov-Kuznetsov 
(ZK) equation
\begin{equation}
    u_t + (u_{xx}+u_{yy} + u^p)_x  = 0, \quad p=2,3,4.
    \label{ZK}
\end{equation}
This equation is a two-dimensional generalization of the well-known 
Korteweg-de Vries (KdV) equation, which is spatially limited as the 1D model of weakly nonlinear waves in shallow water. 
The 2D quadratic ($p=2$) ZK equation governs, for example, weakly 
nonlinear ion-acoustic waves in a plasma comprising cold ions and hot 
isothermal electrons in the presence of a uniform magnetic field 
\cite{MP1999}. In \cite{MM1989} this equation appears as the 
amplitude equation for two-dimensional long waves on the free surface of a thin film flowing down a vertical plane with moderate values of the fluid surface tension and large viscosity. While originally the equation was proposed by Zakharov and Kuznetsov in the 3D setting, see \cite{ZK}, the first rigorous derivation was done by Lannes, Linares and Saut in \cite{LLS} from the Euler-Poisson system. In this paper we initiate numerical investigations of the two-dimensional ZK equation with pure power nonlinearities $p=2,3,4$.

The wellposedness theory for the Cauchy problem for the ZK equation 
with $H^1$ initial data was initiated in \cite{F95}, followed by 
lower regularity improvements in \cite{LP2009}, \cite{FLP} 
\cite{RV2012}, \cite{K2018}. From the local theory it follows that 
solutions to the ZK equation have a maximal forward lifespan $[0,T)$ 
with either $T=+\infty$ or $T<+\infty$. In the later case in the 2D 
setting one has $\|\nabla u(t)\|_{L^2(\mathbb R^2)} \nearrow \infty$ 
as $t \to T$, though the unbounded  growth of the gradient might also happen in infinite time.     

During their existence, solutions to ZK have several conserved quantities, relevant to this work is the $L^2$ norm (or mass), and the energy (or Hamiltonian):
$$
M[u(t)] = \int_{\mathbb R^2} u^2(t) = M[u(0)],
$$
\begin{equation}\label{E}
E[u(t)]=\dfrac{1}{2}\int_{\mathbb 
R^2}\left(u_{x}^{2}(t)+u_{y}^2(t)\right) - \dfrac{1}{p+1}\int_{\mathbb R^2} u^{p+1}(t) = E[u(0)]. 
\end{equation}
Unlike the 1D KdV or modified KdV, the  ZK equation is not integrable for any power $p$. 

One of the useful symmetries in the evolution equations is {\it the scaling invariance}, which states that an appropriately rescaled version of the original solution is also a solution of the equation. For the equation \eqref{ZK} it is
\begin{equation}\label{E:scaling}
u_\lambda(x,y,t)=\lambda^{\frac{2}{p-1}} u(\lambda x, \lambda y,\lambda^3 t).
\end{equation}
This symmetry makes invariant the Sobolev norm $\dot{H}^s$ with $s = 1- \frac2{p-1}$, since 
$\|u_\lambda\|_{\dot{H}^s} = \lambda^{\frac{2}{p-1}+s-1} \|u\|_{\dot{H}^s}$. Moreover, 
the index $s$ gives rise to the critical-type classification of \eqref{ZK}: when $s<0$, or $p<3$, the equation \eqref{ZK} is called the $L^2$-subcritical equation (in this paper a representative of this case is $p=2$); if $p>3$, or $s>0$, the equation is $L^2$-supercritical (we use $p=4$), and with $p=3$, or $s=0$, it is $L^2$-critical. This classification is important when one studies long time behavior of solutions for various nonlinearities. For that we need the notion of solitons. 

The 2D ZK equation has a family of localized traveling waves (or solitary waves, often referred to as solitons), which travel only in $x$ direction
\begin{equation}\label{E:Q}
u(x,y,t) = Q(x-ct,y)
\end{equation}
satisfying
\begin{equation}\label{Q}
-cQ+Q_{xx}+Q_{yy}+Q^{p}=0;
\end{equation}
and defining the ground state solution (i.e., the unique radial positive $H^1$ solution vanishing at infinity, for which the existence, uniqueness and various other properties are well-known, see for example, \cite{SS1999}). 
We note that $Q \in C^{\infty}(\mathbb R^2)$, $\partial_r Q(r) <0$ 
for any $r = |(x,y)|>0$, and that $Q$ has exponential decay $|\partial^\alpha Q(x,y)| \leq c_\alpha \, e^{-r}$ for any multi-index $\alpha$ and any $(x,y) \in \mathbb R^2$.
The solitons \( Q_{c}(x,y) \) are related to the soliton \(Q_{1}(x,y)=:Q(x,y) \) 
for \( c > 0 \) via
\begin{equation} \label{Qscal}
Q_{c}(x,y)=c^{\frac1{p-1}} \,Q(\sqrt{c}\,x,\sqrt{c}\,y),
\end{equation}
thus, it suffices to consider 
\( c=1 \). 

In the $L^2$-subcritical case the local theory together with the Gagliardo-Nirenberg inequality implies that the $H^1$ norm of solutions remains bounded, and thus, all solutions in the subcritical case exist globally in time. In the $L^2$-critical case, using the energy and mass conservation together with the Gagliardo-Nirenberg inequality and its sharp constant expressed in terms of the soliton mass, one has $\|\nabla u\|^2_{L^2} \leq (1-\frac{\|u\|^2_{L^2}}{\|Q\|^2_{L^2}})^{-1} \, E[u]$. Thus, if $\|u_0\|_{L^2} < \|Q\|_{L^2}$, then solutions with the initial condition $u_0$, exist also globally in time, while the blow-up might be possible if the initial mass $\|u_0\|^2_{L^2}$ is greater or equal to that of the soliton $Q$. 

The main aim of this work is to investigate behavior of solutions in various cases of the 2D ZK equation numerically. In particular, we are interested in stability of solitons and in their interaction in the subcritical case, in the scattering and blow-up behavior in the critical and 
supercritical cases. 
For that we mention that the orbital stability of solitons in the 
context of the generalized ZK equation \eqref{ZK} was obtained by de 
Bouard \cite{dB} showing that the traveling waves are orbitally 
stable in the 2D case for $p<3$ and unstable for $p>3$. The 
instability of solitons in the critical case $p=3$ was shown by the 
second author and her collaborators in \cite{FHR3}, see also 
\cite{FHR2} for an alternative proof of instability in the 
supercritical ZK case. The more refined asymptotic stability was 
obtained for $p=2$ by Cote, Mu\~noz, Pilod and Simpson in \cite{CMPS} 
(in fact, for $2 \leq p < p^* \approx 2.15$), in that work the 
authors also studied the interaction of N well-separated solitons. 
The instability of solitons in the critical ($p=3$) case led to 
showing the existence of blow-up in the 2D critical ZK equation, the 
first such work in a higher dimensional generalization of generalized 
KdV (gKdV) equation, see \cite{FHRY}. In that work the blow-up is shown for initial data with negative energy and the mass slightly above the ground state mass. We note that unlike other dispersive equations such as the nonlinear Schr\"odinger equation (NLS), the KdV-type equations (including ZK equation) do not have a convenient virial identity, which gives a straightforward proof of existence of blow-up solutions. Therefore, the proof of existence of blow-up solutions via analytical tools has only been done via construction of such solutions, for example, for the blow-up in 1D critical gKdV see \cite{M}, \cite{MM}. 
\smallskip

In this paper we investigate the following conjectures about the stability of solitons in the $L^2$-subcritical case, about scattering and the stable blow-up dynamics in the $L^2$-critical and supercritical cases\footnote{In our conjectures and simulations we consider exponentially decaying initial data, it will be interesting to investigate slower decay conditions.}.
\begin{conjecture}[$L^2$-subcritical case] \label{C:1} 
Consider the subcritical 2D ZK equation, in particular, when $p=2$ in \eqref{ZK}.
\begin{enumerate}
\item
The soliton solutions \eqref{E:Q}-\eqref{Q}-\eqref{Qscal} are orbitally and asymptotically stable. 
\item 
Solutions of \eqref{ZK} with general sufficiently localized initial data and of sufficient smoothness decompose as $t\to\infty$ into solitons and radiation.
\end{enumerate}
\end{conjecture}

\begin{conjecture}[$L^2$-critical case] \label{C:2} 
Consider the critical 2D ZK equation \eqref{ZK} with $p=3$. 
\begin{enumerate}
\item
If $u_0 \in \mathcal{S}(\mathbb{R}^{2})$ is such that 
$\|u_0\|_{2} < \|Q\|_{2}$, then the solution $u(t)$ to \eqref{ZK} is dispersed. 
\item    
If $u_0 \in \mathcal{S}(\mathbb{R}^{2})$ is sufficiently localized and such that 
$\|u_0\|_{2} > \|Q\|_{2}$, then the solution blows up in finite time $t=t^*$ and
such that as $t \to t^*$
\begin{equation}\label{selfs}
u(x,y,t)- \frac{1}{L(t)}\,Q\left(\frac{x-x_{m}(t)}{L(t)}, \frac{y-y_{m}(t)}{L(t)} \right) \to \tilde{u}\in L^{2},
\end{equation}
with 
\begin{equation}\label{ux2}
\|u_{x}(t)\|_{2} \sim \frac{1}{L(t)}, ~~ L(t)\sim \sqrt{t^*-t}, \quad \mbox{and} \quad 
x_{m}(t) \sim \frac{1}{t^*-t}, ~~ y_m(t) \to y^* <\infty.
\end{equation}
\end{enumerate}
\end{conjecture}

\begin{conjecture}[$L^2$-supercritical case] \label{C:3}
Consider the supercritical 2D ZK equation, in particular, when $p=4$ in \eqref{ZK}. 
Let $u_0 \in \mathcal{S}(\mathbb R^2)$ be of sufficiently large mass and energy 
\footnote{We have not investigated numerically the precise value. For some thresholds, for example, see \cite{FLP}.} and of some localization. 
Then ZK evolution $u(t)$ blows up in finite time $t^{*}$ and finite location $(x^{*},y^{*})$, i.e., the blow-up core resembles a self-similar structure with  
\begin{equation}\label{selfss}
u(x,y,t)- \frac{1}{L^{\frac2{p-1}}(t)}\, P\left(\frac{x-x_{m}(t)}{L(t)}, \frac{y-y_{m}(t)}{L(t)}\right)\to \tilde{u}\in L^{2},
\end{equation}
where $P(x,y)$ is a localized solution to \eqref{ZKresinfty} (which 
is conjectured to exist),  
$$
x_m(t) \to x^*, \quad y_m(t) \to y^{*},
$$
and 
\begin{equation}\label{ux2s}
\|u_{x}(t)\|_{2}\sim \frac{1}{L^{\frac2{p-1}}(t)} \quad \mbox{with} \quad  L(t)\sim (t^{*}-t)^{1/3} \quad \mbox{as} \quad t \to t^*.
\end{equation}
\end{conjecture}

\begin{remark}\label{buremark}
We note that numerical blow-up computations are extremely challenging, 
since they push the limits of the best currently available methods, approaches and computational power. This is especially true for dispersive equations, where also the radiation should be correctly approximated. In a sense, this can be seen as an invitation to analytical studies of the phenomena shown in this paper (see some work in this direction \cite{FHRY}, \cite{FHR3}, \cite{CMPS}). 
Nonetheless, the techniques applied here have been successfully tested on an example of the gKdV equation, for which the analytical description is much better understood, though far from being complete. 
\end{remark}

The paper is organized as follows: 
In Section \ref{S:Num} we present the numerical tools used to solve the ZK equation. 
Examples for the $L^{2}$-subcritical case are discussed in Section \ref{S:sub}. 
The $L^{2}$-critical case is studied in Section \ref{S:cubic}. 
In Section \ref{S:last} we discuss examples for the $L^{2}$-supercritical case. 

\subsection{Acknowledgements}
CK and NS were partially supported by 
the ANR-FWF project ANuI - ANR-17-CE40-0035, the isite BFC project 
NAANoD, the ANR-17-EURE-0002 EIPHI and by the 
European Union Horizon 2020 research and innovation program under the 
Marie Sklodowska-Curie RISE 2017 grant agreement no. 778010 IPaDEGAN.
SR was partially supported by the NSF grant DMS-1815873/1927258, she would also like to thank the AROOO (`A Room of Ones's Own') initiative for focused research time for this project. 


\section{Numerical methods}\label{S:Num}
In this section we review the numerical methods to be applied in 
the rest of the paper. First, we construct the ZK solitons via an 
iterative approach. Then we introduce the integration of the ZK 
equation with a Fourier spectral method for the spatial coordinates 
and a fourth order scheme in time. 
Finally, we review the dynamic rescaling method, which is used to track blow-up solutions. 

\subsection{Solitons}
We first obtain the soliton solutions for the equations \eqref{ZK} by 
solving equation \eqref{Q}. Since this is also the defining equation 
for the solitons of the NLS equation in 2D, it is known that its 
solutions have radial symmetry. Here, we do not use this fact, since we 
intend to apply Fourier methods throughout the paper, and thus, directly 
construct the solitons on the grids for the time evolution. 

To this end, we use discrete Fourier transforms in both $x$ and $y$, 
which is, loosely speaking, equivalent to approximating a function via 
a truncated Fourier series. Since it is known that the NLS solitons are 
rapidly decreasing functions, they can be treated as periodic smooth functions on sufficiently 
large periods within the finite numerical precision. 
We work with $x\in L_{x}[-\pi,\pi]$ and $y\in L_{y}[-\pi,\pi]$, 
where $L_{x}$ and $L_{y}$ are positive real numbers, chosen so that the Fourier coefficients decrease both in $x$ and $y$ to machine precision (which is of the order of 
$10^{-16}$ in double precision). We denote the dual Fourier variables 
to $x$ and $y$ by $k_{x}$ and $k_{y}$, respectively, and write
\begin{equation}
    Q(x,y)\approx \sum_{k_{x}=-N_{x}/2+1}^{N_{x}/2}
    \sum_{k_{y}=-N_{y}/2+1}^{N_{y}/2}\hat{Q}(k_{x},k_{y})\, e^{i(k_{x}x+k_{y}y)}
    \label{Qhat};
\end{equation}
the \emph{discrete Fourier transform} $\hat{Q}=\mathcal{F}u$ can be conveniently 
computed with a \emph{fast Fourier transform} (FFT). An advantage of 
Fourier methods is that the numerical resolution can be controlled 
via the decay of the Fourier coefficients, the highest coefficients 
indicate the numerical error introduced by the truncation of the 
series. 

\begin{figure}[!htb]
\includegraphics[width=0.32\hsize]{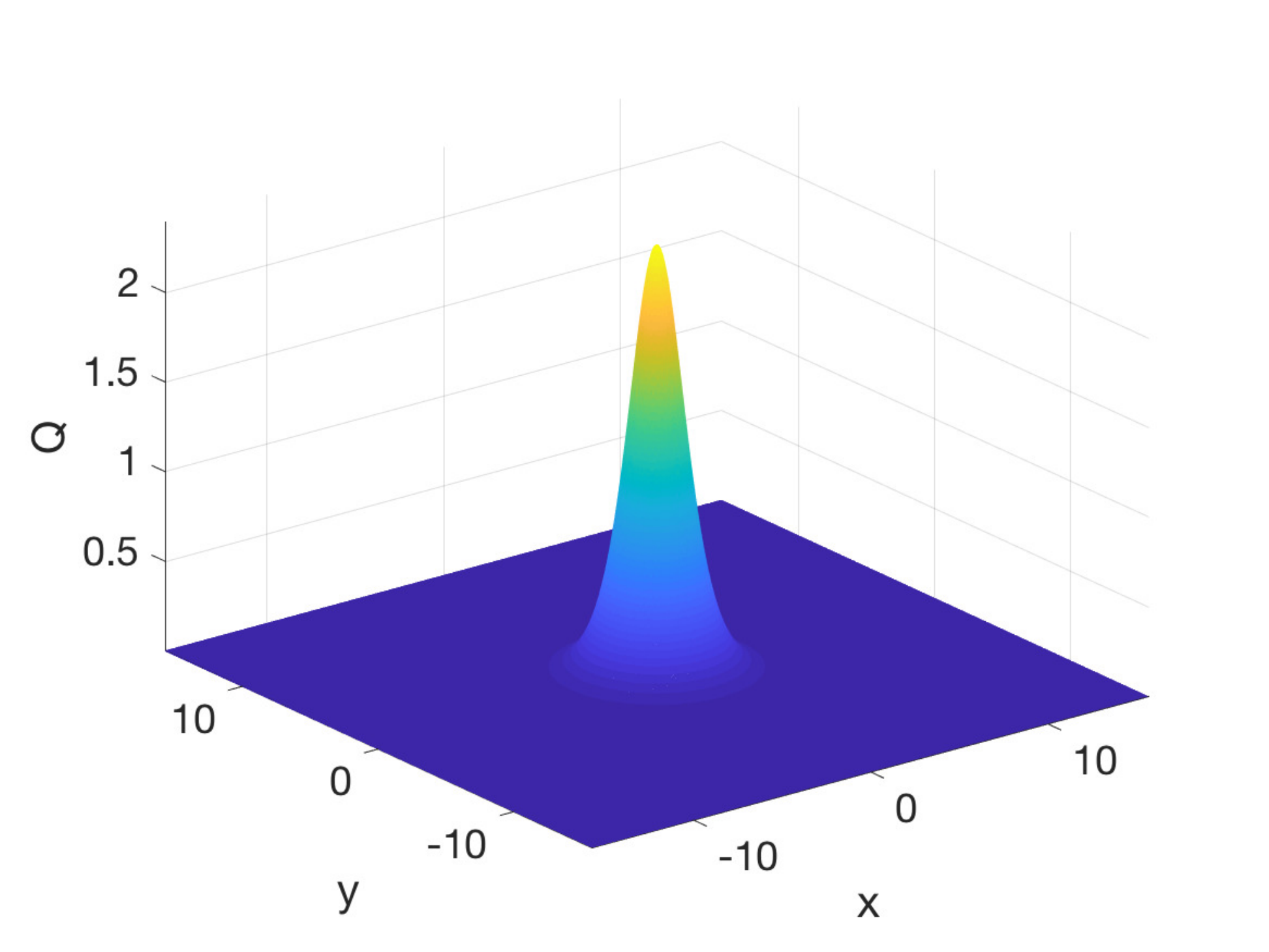}
\includegraphics[width=0.32\hsize]{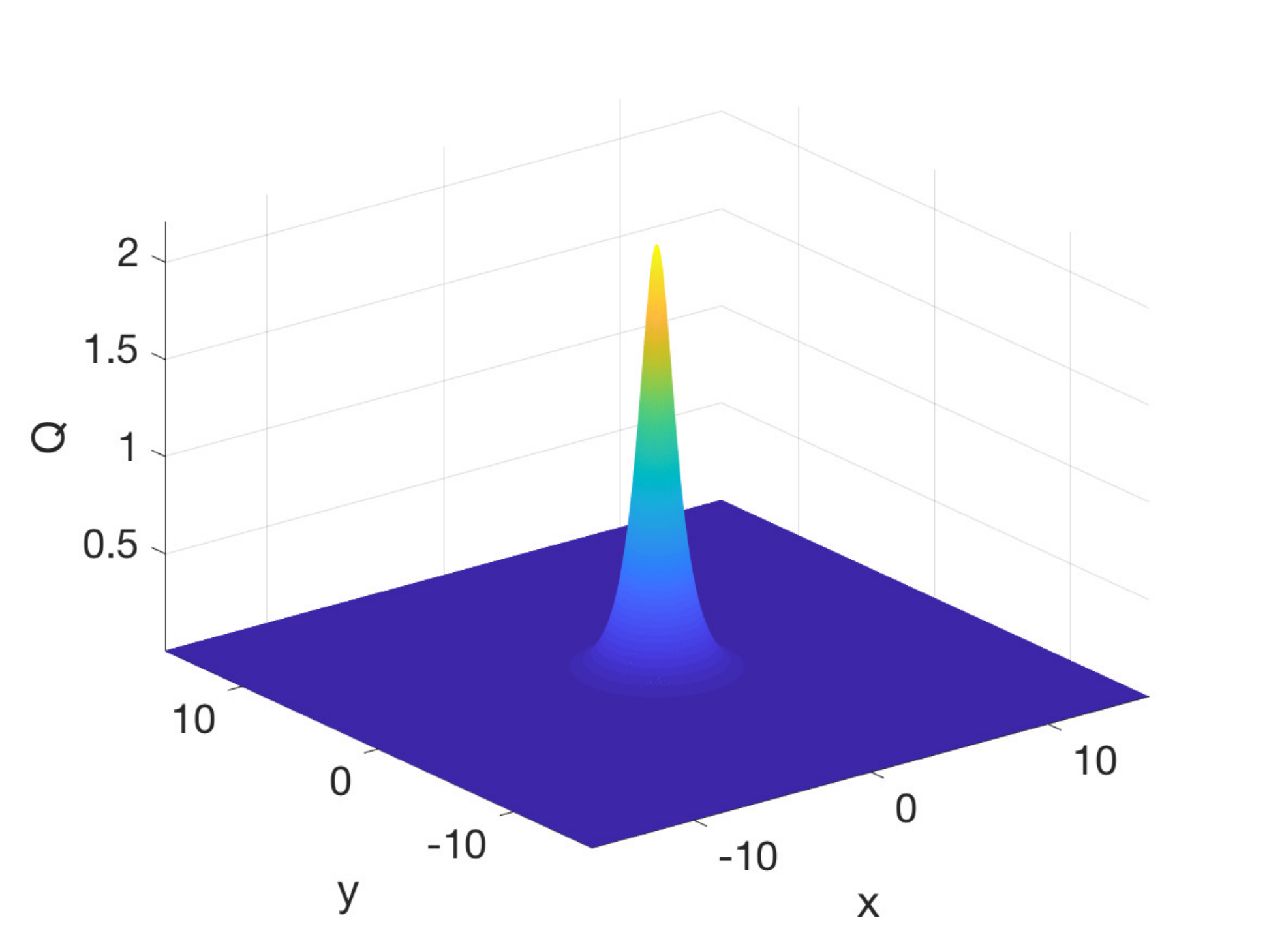}
\includegraphics[width=0.32\hsize]{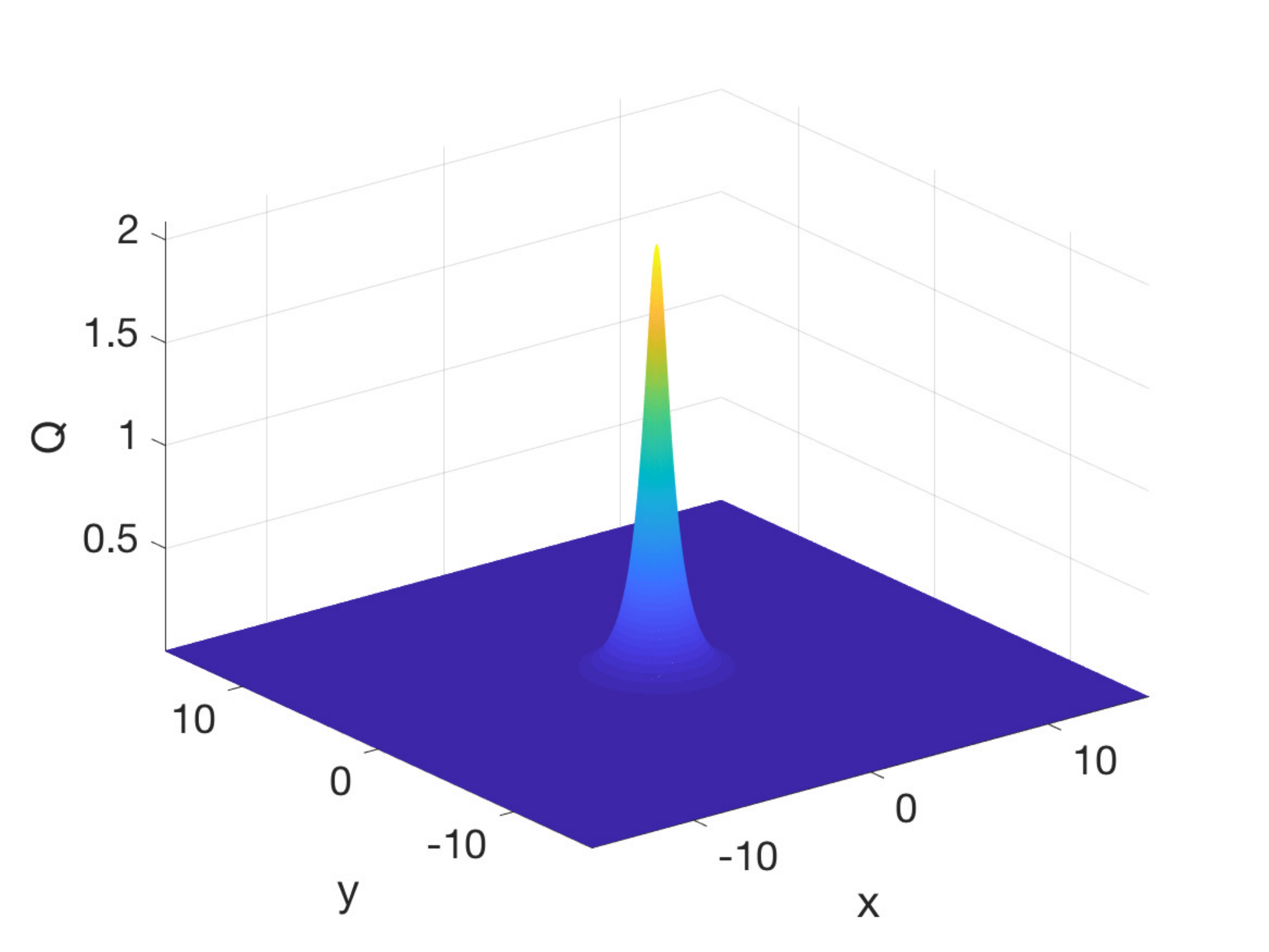}
\caption{Solitons to the ZK equation (\ref{Q}) for $c=1$ and 
$p=2,3,4$ from left to right. }
\label{solfig}
\end{figure}

With this Fourier discretization,  equation \eqref{Q} is approximated by an $N_{x}N_{y}$ dimensional system of nonlinear equations for the $\hat{Q}$. The latter will be iteratively solved by a Newton-Krylov iteration. This means that we invert the 
Jacobian via Krylov subspace methods as in \cite{AKS}, here GMRES \cite{gmres}.  
We use $N_{x}=N_{y}=2^{10}$, $L_{x}=L_{y}=10$ and $Q=2\,e^{-x^{2}-y^{2}}$ 
as initial iterates in all cases. The iteration is stopped when the 
residual is smaller than $10^{-10}$. The solitons for $c=1$ and 
$p=2,3,4$ are shown in Fig.~\ref{solfig}. It can be seen that they 
become more localized and slightly smaller with increasing nonlinearity (the maximum value, which is also the value at zero, decreases: $Q(0) = 2.3920$ if $p=2$, 
$Q(0) = 2.2062$ if $p=3$, and $Q(0) = 2.0853$ if $p = 4$.  
The Fourier coefficients decrease in all cases to machine precision, 
see Fig.~\ref{test} on the left for $p=3$, which implies that the 
solution is spatially well resolved. 

\subsection{Time evolution}
The same Fourier discretization as for the soliton  above is used for 
the full ZK equation \eqref{ZK}, which is thus approximated by an 
$N_{x}N_{y}$ dimensional system of ordinary differential equations in $t$ of the form
\begin{equation}
    \hat{u}_{t}= \mathcal{L}\hat{u}+\mathcal{N}[\hat{u}]
    \label{sys},
\end{equation}
where $\mathcal{L}=ik_{x}(k_{x}^{2}+k_{y}^{2})$ and 
$\mathcal{N}[\hat{u}]=-ik_{x}\mathcal{F}(u^{p})$. Because of the 
appearance of third derivatives in $x$ and $y$, 
this system is \emph{stiff}, implying that explicit methods will be 
inefficient due to stability conditions as they necessitate prohibitively small times 
steps in order to stabilize the code. 
Implicit schemes are less restrictive in this sense, but are 
computationally expensive, since the resulting nonlinear equation has 
to be solved in each time step. Therefore, we have compared in 
\cite{etna,KR} various adapted integrators for stiff systems with a 
diagonal $\mathcal{L}$ as we have here, which 
are explicit and of fourth order. It turned out that 
\emph{exponential time differencing} (ETD) schemes, see \cite{HO} for a comprehensive 
review with many references, are most efficient  in the context 
of KdV-type equations. There are various fourth order ETD methods, which 
all showed a similar performance in our tests. Here, we apply the method 
by Cox and Matthews \cite{CM} in the implementation described in 
\cite{etna, KR}. The accuracy of the time integration scheme can be 
controlled via the conserved energy of the equation. Due to 
limitations in the accuracy of numerical methods, the 
computed energy (again Fourier 
techniques are applied to \eqref{E}) will not be 
exactly conserved. The quantity $\Delta E = |E(t)/E(0)-1|$ can be 
used as discussed in \cite{etna,KR} as an estimate of the numerical 
error. Typically it overestimates the accuracy of the numerical 
solution by 1-2 orders of magnitude. 

As far as the blow-up is concerned, it is numerically very challenging to study blow-up solutions. For the generalized KdV and KP equations, this was done, for instance, in \cite{KP1,KP2}. There it was shown that the integration of the dynamically rescaled 
equation \eqref{ZKresc} is problematic if Fourier methods are used. 
Instead in \cite{KP1,KP2,KS} the equations were integrated without 
rescaling, and then a postprocessing of the results was done according to \eqref{resc} to identify the type of the blow-up. The same strategy will be applied here to ZK. However, the generalized KdV 
equations have the additional complication that the blow-up occurs at 
infinite values of $x$, and that the blow-up profile is leaving the 
initial location with infinite velocity. To treat such cases, in \cite{KP1, AKS} we introduced a reference frame, in which the maximum of the solution is stationary at some point 
$x_{m}$ during the whole computation, i.e., an accelerated reference frame. 
This means we apply \eqref{resc} with $L=1$ and $y_{m}=0$ and solve
\begin{equation}\label{ZKnum}
u_t + (u_{xx}+u_{yy} + u^p)_x -v_{x}u_{x}  = 0, \quad p=2,3,4,
\end{equation}
where $u(x_{m})$ is taken to be a maximum of the solution for all times. By 
differentiating equation (\ref{ZKnum}) with respect to $x$ and 
evaluating it for $x=x_{m}$, we get 
\begin{equation}
    v_{x}= \left.\frac{(u_{xx}+u_{yy} + u^p)_{xx}}{u_{xx}}
    \right|_{x=x_{m}}
    \label{vx}.
\end{equation}
Since it is computationally expensive to compute $v_{x}$ in each time 
step, we only apply this approach for blow-up computations in the 
$L^{2}$ critical case. 

\subsection{Test}

\begin{figure}[!htb]
 \includegraphics[width=0.49\hsize]{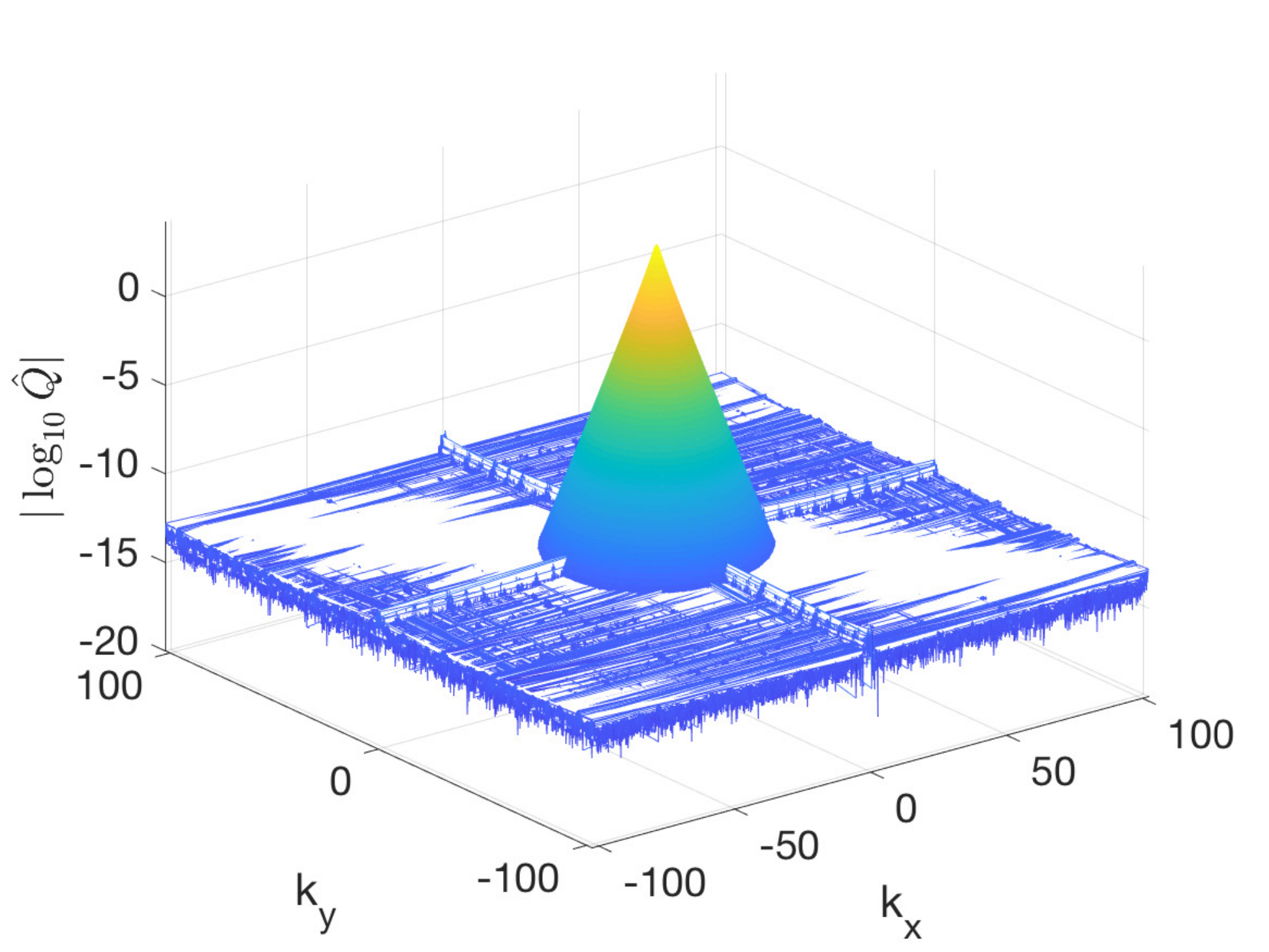}
 \includegraphics[width=0.32\hsize]{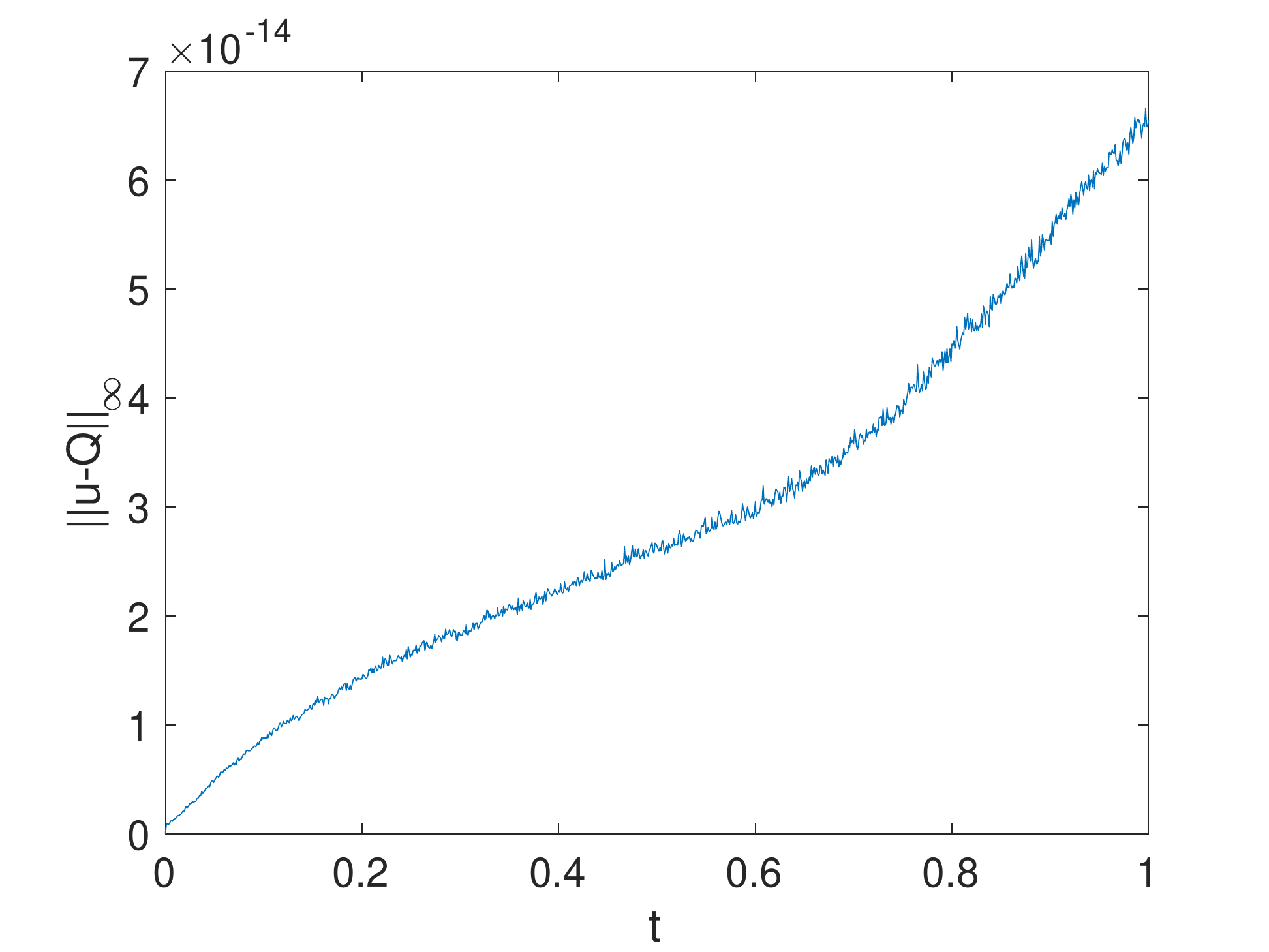}
\caption{Left: Modulus of the Fourier coefficients of the ZK soliton 
with $c=1$ and $p=3$. Right: the difference between the numerically 
computed solution to \eqref{ZK} with $p=3$ with the soliton initial 
condition and the soliton from \eqref{Q}, depending on time. }
\label{test}
\end{figure}

To test the time evolution code and the soliton at the same time, 
we consider the soliton in the $L^2$-critical case ($p=3$) as initial data 
and a co-moving frame with $v_{x}=c=1$. For $t\in[0,1]$ we apply 
$N_{t}=1000$ time steps. The numerically computed energy is conserved 
to the order of $10^{-14}$. The difference between the numerically 
computed solution and the soliton can be seen on the right of 
Fig.~\ref{test}. It increases with time, but is of the order of 
$10^{-14}$ as the energy conservation. Though we show later that the 
soliton is unstable against both dispersion and blow-up, the code is 
able to propagate it on the considered time intervals with 
essentially machine precision. 

\subsection{Dynamic Rescaling}
Recalling the scaling invariance \eqref{E:scaling} for the equation \eqref{ZK},
one can use this symmetry 
in the context of blow-up in the form of a \emph{dynamical rescaling}
\begin{equation}\label{resc}
\begin{array}{c}
X = \frac{x-x_{m}(t)}{L(t)}, \quad Y = \frac{y-y_{m}(t)}{L(t)}, \quad 
    T=\int_{0}^{t}\frac{dt'}{L^{3}(t')},\\
    ~\\
U(X,Y,T) = L^{\frac2{p-1}}(t) \, u(x,y,t).
\end{array}
\end{equation}
The dynamically rescaled ZK equation reads
\begin{equation}\label{ZKresc}
U_T - a\bigg(\frac2{p-1} U+X U_{X}+Y U_{Y}\bigg) - v_{X}U_{X}-v_{Y}U_{Y}+ \bigg(U_{XX}+U_{YY} + U^p\bigg)_X  = 0,
\end{equation}
where 
\begin{equation}\label{a}
a \equiv a(T) = \frac{d\ln L}{dT},\quad v_{X}=\frac{x_{m,T}}{L},\quad v_{Y}=\frac{y_{m,T}}{L}.
\end{equation}
It is assumed that blow-up happens as $T\to\infty$, and that \( 
U_{T} \) vanishes in this limit. Thus, the equation 
\eqref{ZKresc} in the limit becomes 
\begin{equation}\label{ZKresinfty}
-\overset{\infty}{a}\bigg(\frac2{p-1}\overset{\infty}{U}+ X\overset{\infty}{U}_X+Y\overset{\infty}{U}_Y \bigg) - 
 v_{\underset{\infty}{X}}\overset{\infty}{U}_{X}-v_{\underset{\infty}{Y}}\overset{\infty}{U}_{Y}+ \bigg(\overset{\infty}{U}_{XX}+\overset{\infty}{U}_{YY} + \overset{\infty}{U}~^p \bigg)_X = 0,
\end{equation}
where the sub/superscript \( \infty \) denotes that the quantity is taken in the limit as $T\to\infty$ and $ \overset{\infty}{U} $ stands for a blow-up profile. 

Two possible stable blow-up mechanisms are expected in KdV-type equations: either an algebraic dependence of \( L \) on \( T \), or an exponential 
one. In the former case the quantity \( \overset{\infty}{a} \) in 
\eqref{a} will vanish, and equation \eqref{ZKresinfty} will be 
identical to the equation for the soliton if 
$v_{\underset{\infty}{Y}}=0$;
this mechanism is expected in the $L^{2}$-critical case. If 
$L\propto 1/T$ as in the $L^{2}$-critical gKdV case, recalling \eqref{resc}, we get 
\begin{equation}\label{Lcritical}
L\propto \sqrt{t^{*}-t}.
\end{equation}
In the supercritical case, one expects an exponential decay of $L$ 
with $T$, that is, $L\propto \exp(-\gamma T)$ with $\gamma>0$, and from \eqref{resc} we have
\begin{equation}\label{Lsuper}
L\propto (t^{*}-t)^{1/3}.
\end{equation}
In our simulations we trace the $L^{\infty}$ norm of $u$, the 
$L^{2}$ norm of $u_{x}$, and in the $L^2$-critical case the velocity  $v_x$.  The first two norms are proportional to $L^{\frac2{p-1}}(t)$ via rescaling in \eqref{resc}, see details in Sections \ref{S:cubic} and Section \ref{S:last}. 
 We also note that this is similar to the blow-up situation in NLS-type equations, for example, see \cite{YRZ2019}, \cite{SS1999}.

\section{The $L^{2}$-subcritical case}\label{S:sub}
In this section we study the ZK equation in the subcritical case \( 
p=2 \). We consider the stability of the solitons, the interaction of 
solitons and the appearance of solitons in the long time evolution of 
general localized initial data. 

We work with \( L_{x}=L_{y}=10 \) and \( 
N_{x}=N_{y}=2^{10} \) 
Fourier modes and \( N_{t}=2000 \) time steps on the considered time 
intervals. In all studied cases the Fourier coefficients decrease at 
least to the order of $10^{-5}$, and the relative energy is conserved 
at least to the same order (except for the examples in the last 
subsection, these numbers are in general of the order of $10^{-10}$). 
This means that the numerical error is in all cases much smaller than 
plotting accuracy. 

The results of this section give positive confirmation to the Conjecture \ref{C:1}.

\subsection{Soliton stability}
We first address stability of solitons by considering initial 
data of the form $u_0 \equiv u(x,y,0)=\lambda Q(x,y)$, $\lambda \in \mathbb R$. 
We use a co-moving frame, i.e., we solve
\begin{equation}
u_t + (u_{xx}+u_{yy} + u^p - vu)_x  = 0, \quad p=2,3,4,
    \label{ZKcom}
\end{equation}
with $v=1$.

\begin{figure}[!htb]
 \includegraphics[width=0.49\hsize]{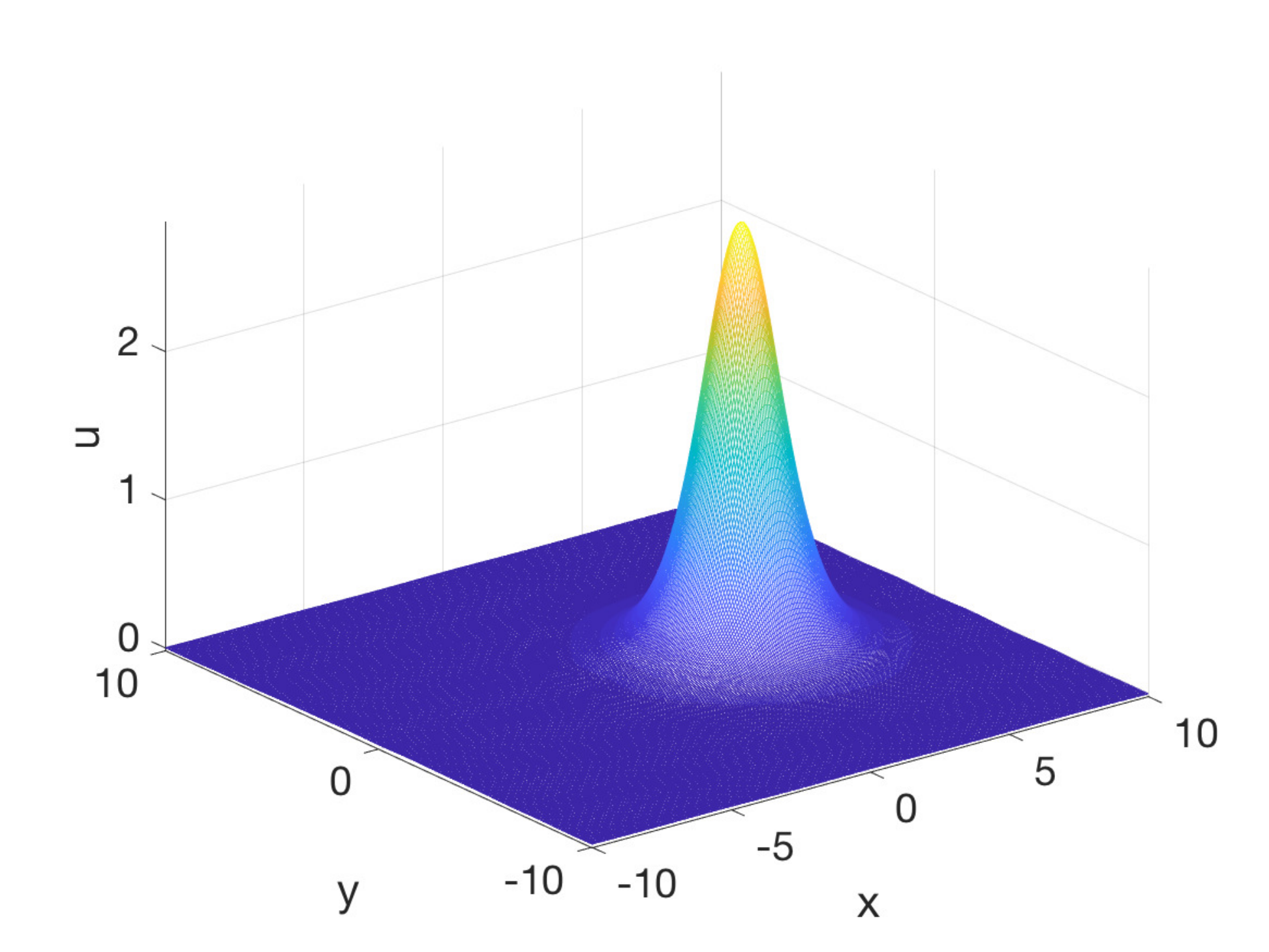}
 \includegraphics[width=0.49\hsize]{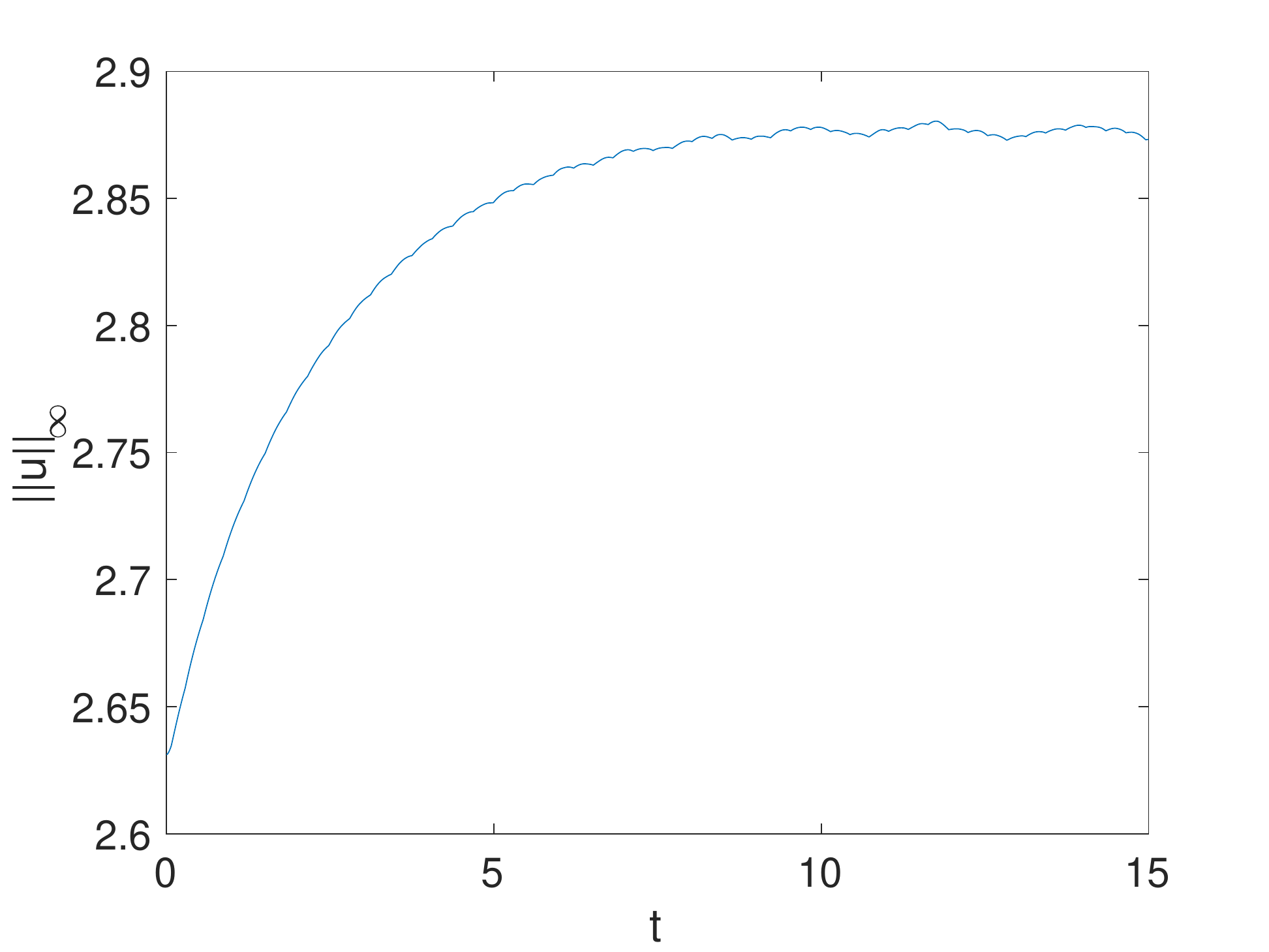}
\caption{Solution to \eqref{ZK} with $u(x,y,0)=1.1Q(x,y)$: on the left the solution at $t=15$; on the right the $L^{\infty}$ norm of the solution depending on time. }
\label{ZKp2sol11t15}
\end{figure}

The solution to \eqref{ZK} with $\lambda=1.1$ can be seen at $t=15$ on the 
left of Fig.~\ref{ZKp2sol11t15}. The perturbed soliton (with $u_0 = 1.1Q$) visibly moves 
faster than the original soliton ($u_0 = Q$), since it is not stationary in the 
co-moving frame. There is also some radiation propagating into the negative $x$-direction (which is more visible in Fig.~\ref{ZKp2soldiff} below). The $L^{\infty}$ norm of the solution 
on the right of Fig.~\ref{ZKp2sol11t15} also appears to saturate at a 
higher value than the initial value. Thus, it seems that the 
perturbation with higher mass than the original soliton leads to a  
larger and faster moving to the right soliton and some radiation moving to the left. 

\begin{figure}[!htb]
 \includegraphics[width=0.49\hsize]{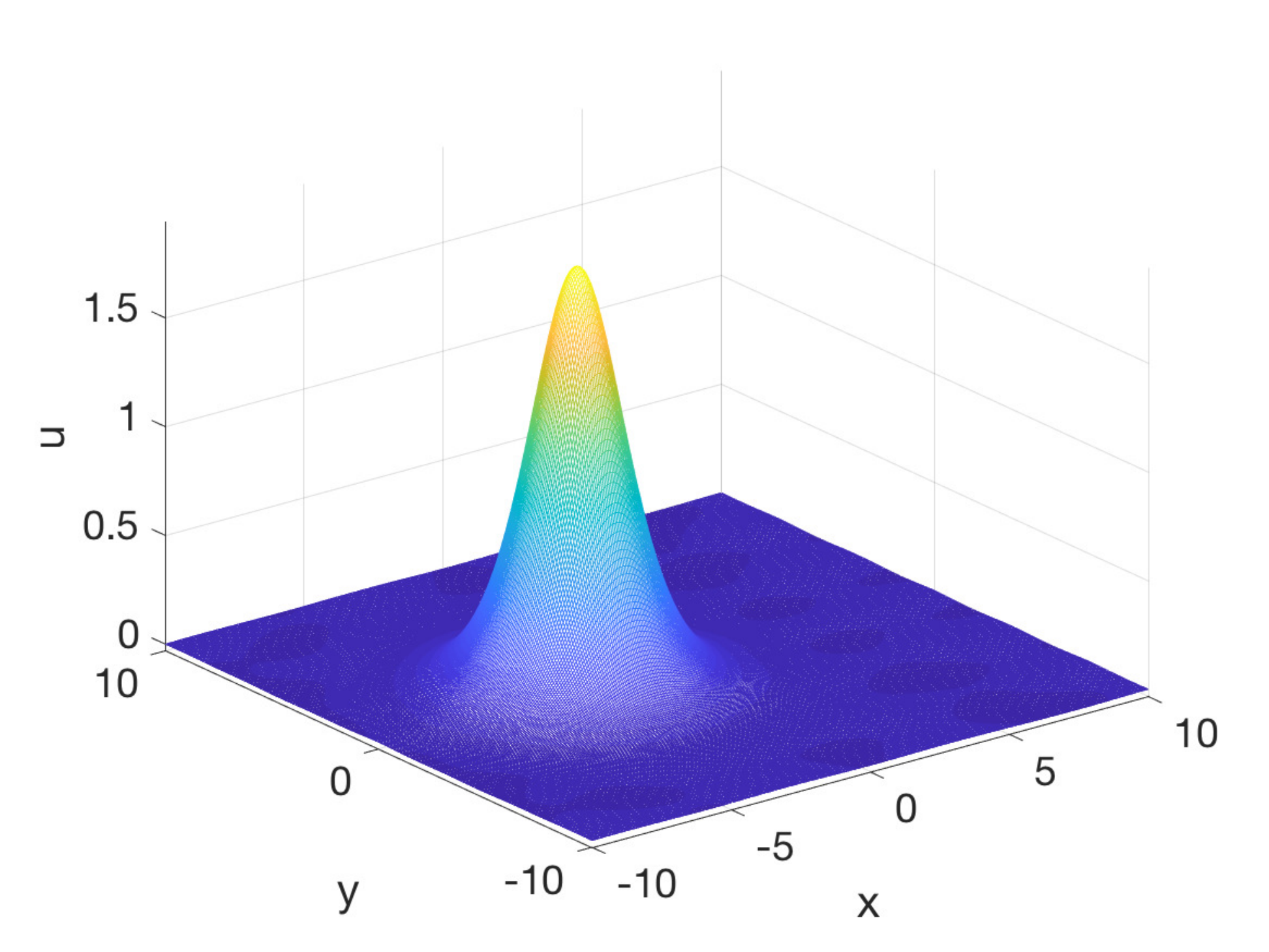}
 \includegraphics[width=0.49\hsize]{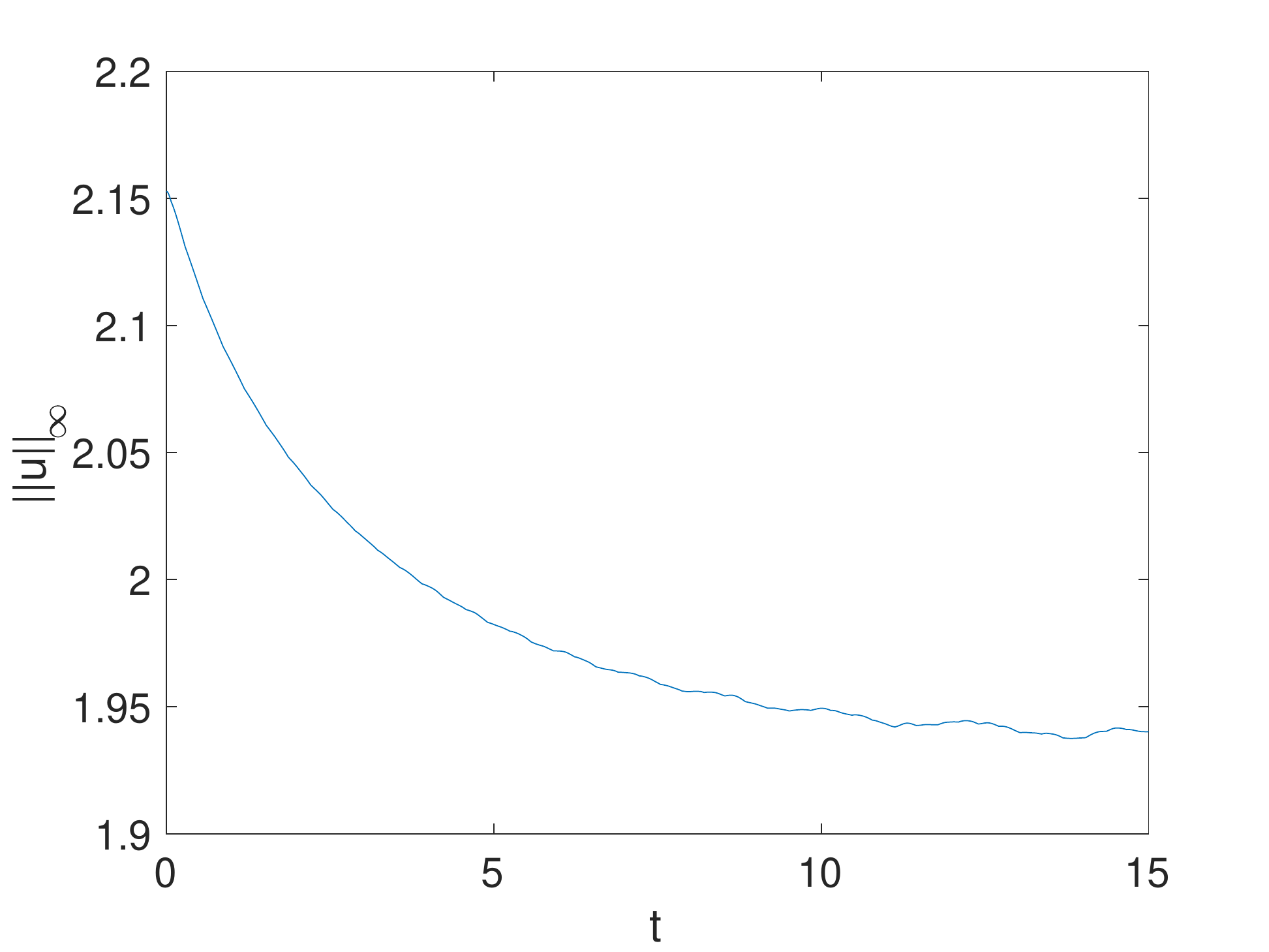}
\caption{Solution to \eqref{ZK} with $u(x,y,0)=0.9Q(x,y)$: 
on the left the solution for $t=15$; on the right the $L^{\infty}$ norm of the 
solution depending on time. }
\label{ZKp2sol09t15}
\end{figure}

On the other hand, a perturbation with smaller mass $(u_0 = 0.9Q)$ is shown on the 
left of Fig.~\ref{ZKp2sol09t15}. The fact that the hump moves to the 
left of the origin in the co-moving frame indicates that the resulting 
soliton has even smaller mass than the perturbed one. This is in 
accordance with the $L^{\infty}$ norm of this solution shown on 
the right of Fig.~\ref{ZKp2sol09t15}: it appears to decrease to a lower height than the initial value. Thus, the perturbed initial 
data of a mass smaller than the perturbed soliton seem to lead to a 
soliton of smaller mass plus radiation. Figs.~\ref{ZKp2sol11t15} and 
\ref{ZKp2sol09t15} indicate that the ZK soliton is stable, as expected, when $p=2$. 

\begin{remark}
Note that in this paper we systematically approximate   situations on 
$\mathbb{R}^{2}$ by a setting on $\mathbb{T}^{2}$. Within machine 
precision, this does not make a difference for stationary localized 
solutions as the solitons of the ZK equation, if the periods are 
chosen sufficiently large. However, if radiation appears, as 
it happens in this and in the following sections, one would have to 
choose prohibitively large computational domains to avoid the 
reappearance of emitted radiation (always emitted in the negative 
$x$-direction) for positive values of $x$. This is acceptable as long 
as this radiation has much smaller amplitudes than the studied bulk of 
the solution. Effects of the radiation can be seen in 
Fig.~\ref{ZKp2sol11t15} and \ref{ZKp2sol09t15} in the variations of 
the $L^{\infty}$ norms for large times (this is also connected with the 
determination of the $L^{\infty}$ norm on a discrete grid, which means 
that the real location of the maximum of the solution might not be on 
a grid point). 
\end{remark}

One may ask to which extend the final states of the solutions shown 
in Fig.~\ref{ZKp2sol11t15} and \ref{ZKp2sol09t15} are solitons if the 
radiation cannot escape the computational domain. To address this 
question we show in Fig.~\ref{ZKp2soldiff} the difference between the 
solutions for $t=15$ and a fitted soliton rescaled according to 
(\ref{Qscal}) ($c$ is determined via 
$c=||u||_{\infty}/||Q||_{\infty}$). In Fig.~\ref{ZKp2soldiff}, one 
can clearly see how the radiation forms a background in the 
computational domain, and that the difference between the bump and a 
soliton is smaller than the radiation background. Thus, 
the conclusion that the final state is a soliton plus 
radiation, escaping on $\mathbb{R}^{2}$ to infinity is justified.
Even more can be seen in the previous figures: the radiation escapes to the `left' of the moving rightward solution at an angle of $30^0$ with the negative $x$-axis (for a total opening of $60^0$). This is in confirmation of the asymptotic stability result in \cite{CMPS}.    
\begin{figure}[!htb]
 \includegraphics[width=0.49\hsize]{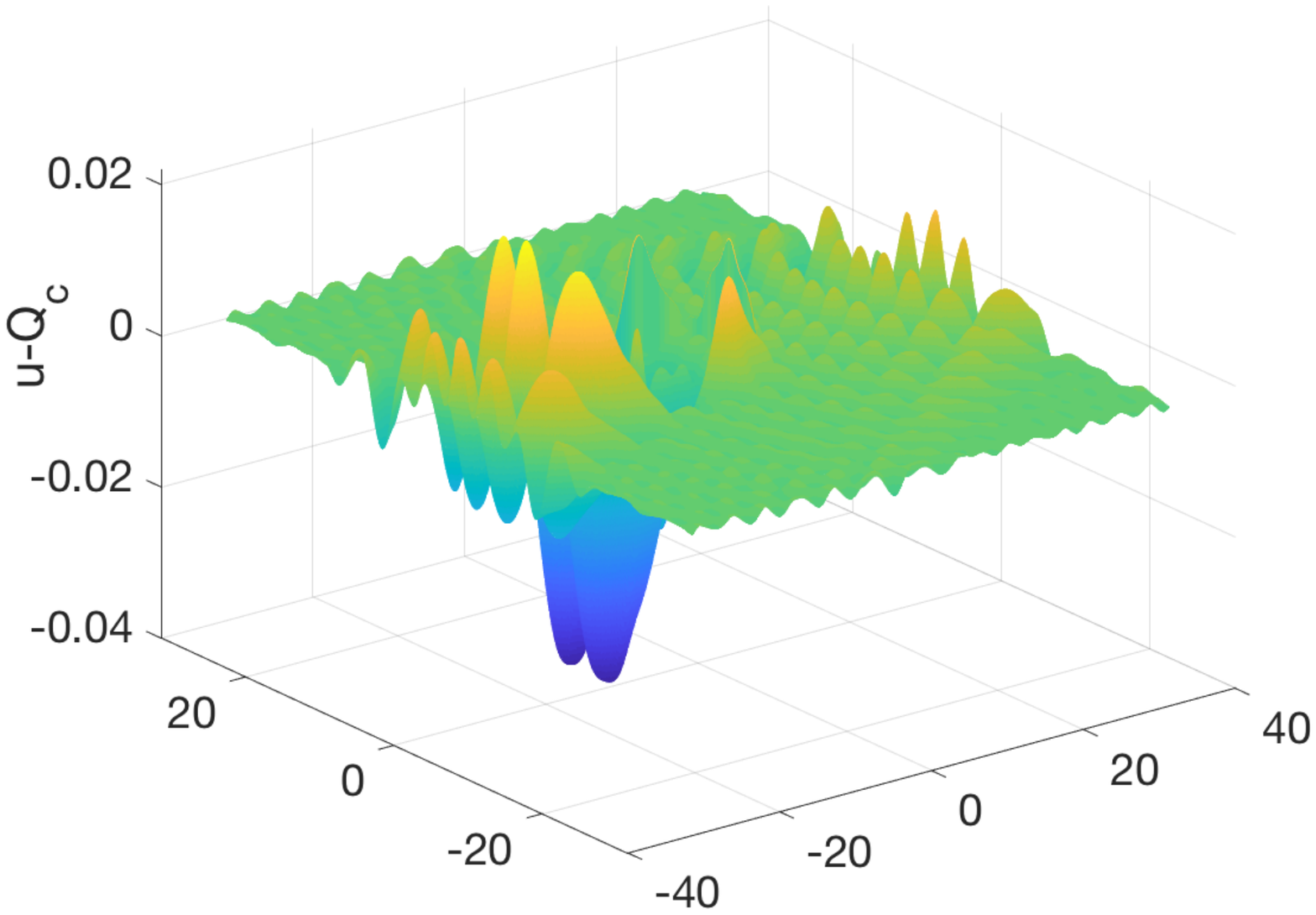}
 \includegraphics[width=0.49\hsize]{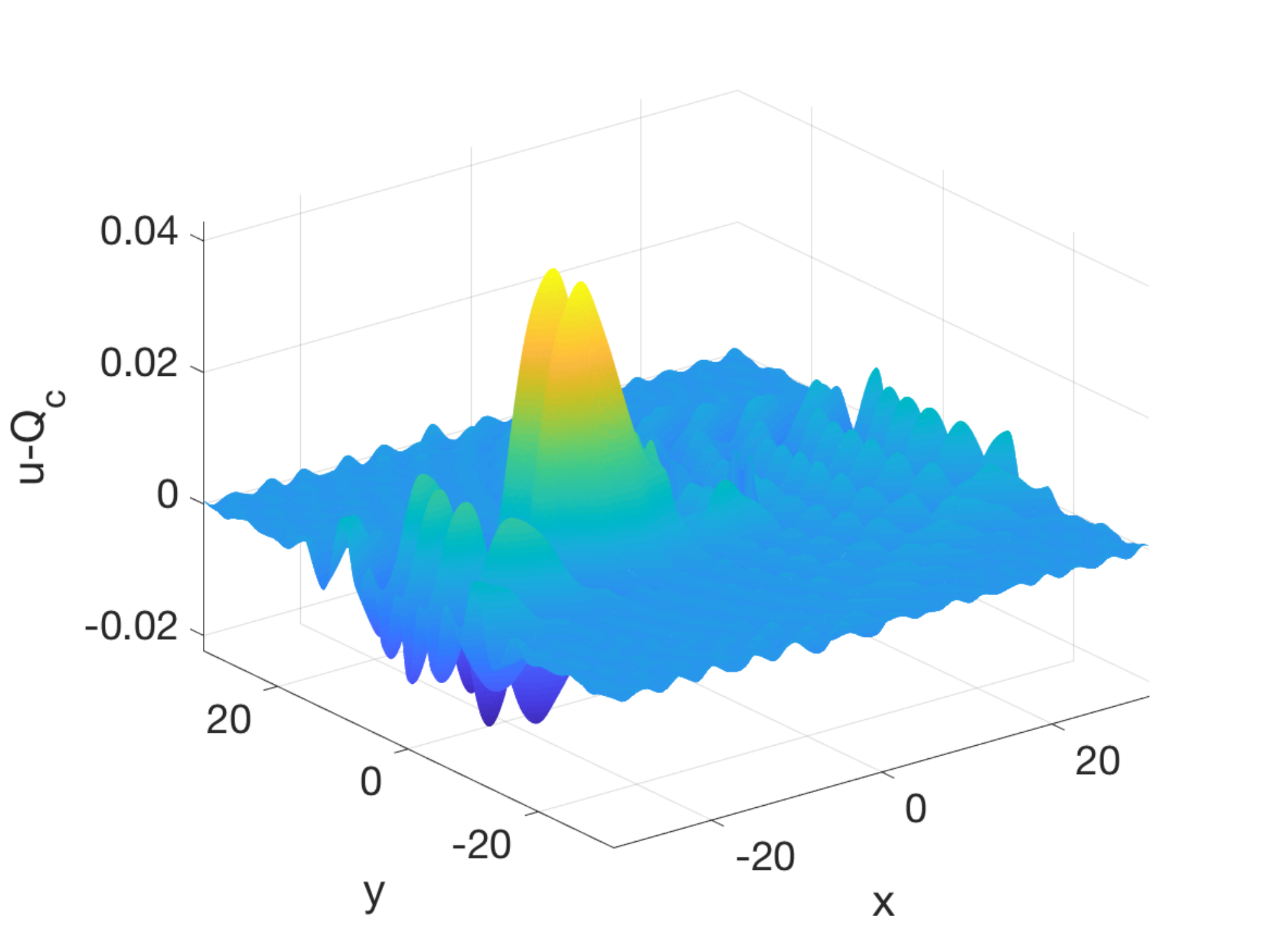}
\caption{Difference of the solution to \eqref{ZK} with $u(x,y,0)=\lambda Q(x,y)$ and a fitted rescaled soliton from \eqref{Qscal}: on the left $\lambda=0.9$; on the right $\lambda=1.1$. }
\label{ZKp2soldiff}
\end{figure}

\subsection{Soliton interaction}
Since the soliton solutions are rapidly decreasing for \( 
x^{2}+y^{2}\to\infty \), one can 
study their interactions by  considering initial  data which are the 
sum of displaced solitons. This allows to study multi-soliton 
solutions to what is clearly a non-integrable equation. Initial data 
can be constructed by superimposing two one-soliton solutions that are sufficiently far. 
Indeed, since solitons have an exponential decay, their contribution far away from their joint center of mass is zero within the numerical precision.

\begin{figure}[!htb]
 \includegraphics[width=0.49\hsize]{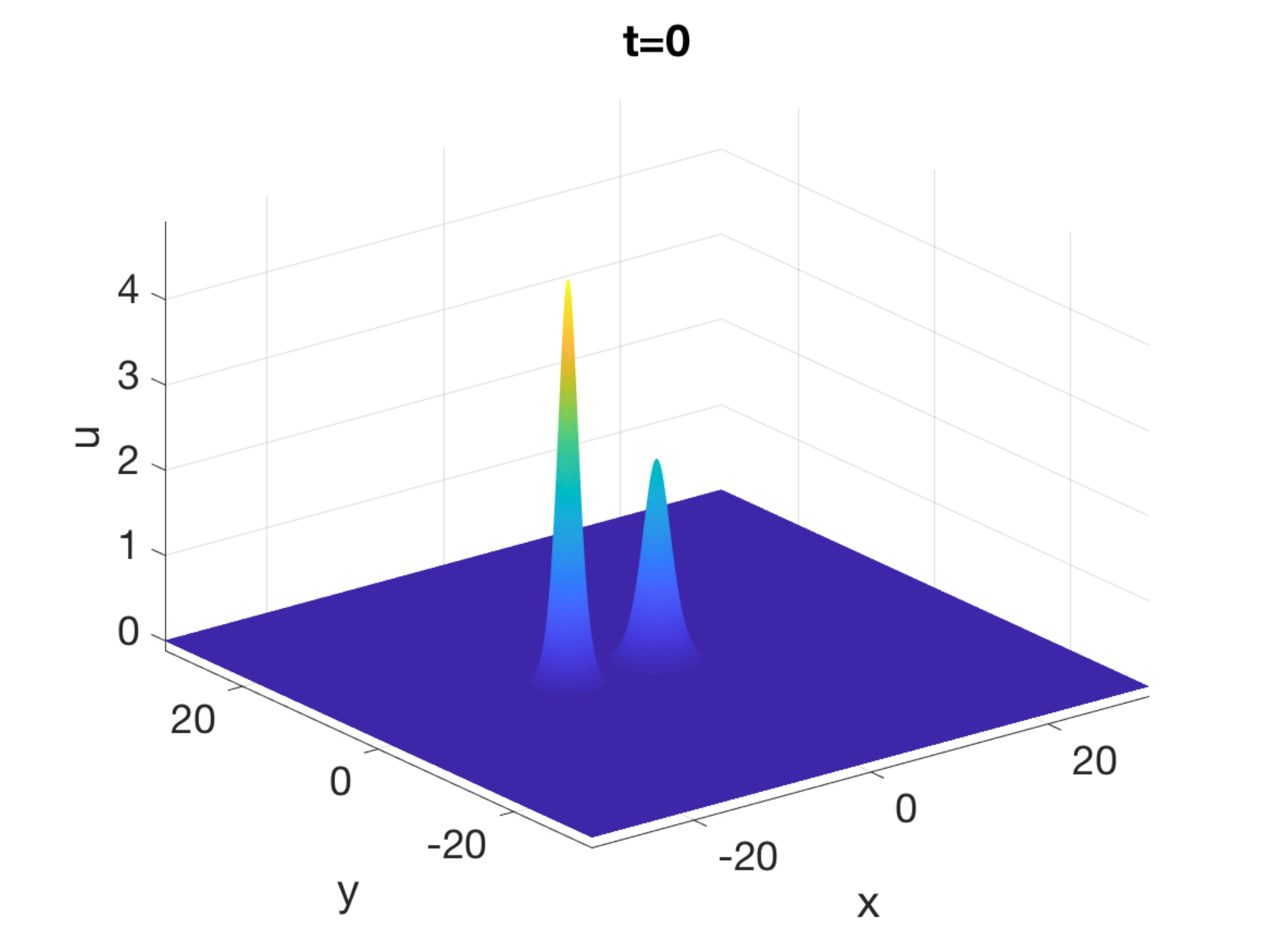}
 \includegraphics[width=0.49\hsize]{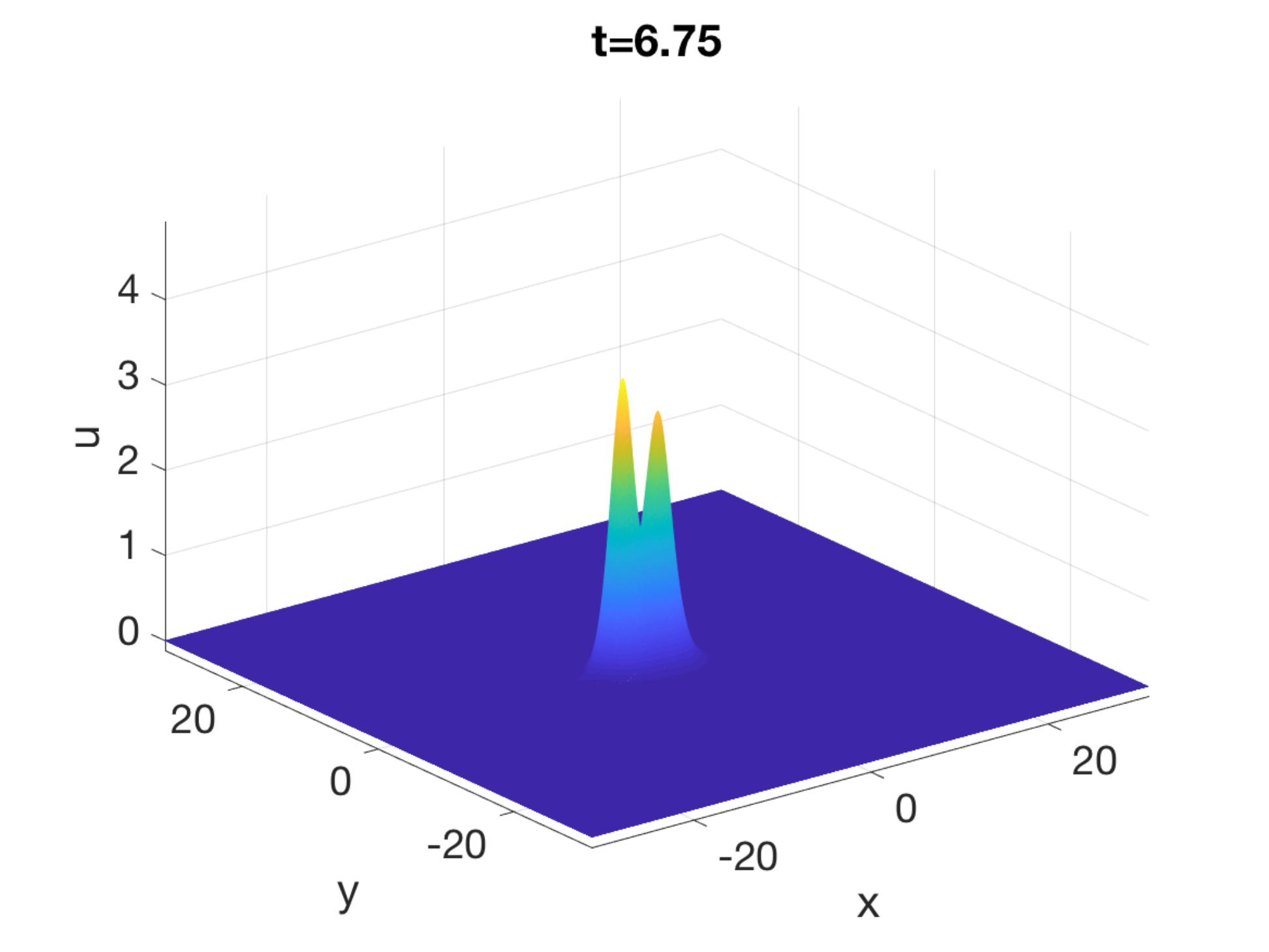}\\
 \includegraphics[width=0.49\hsize]{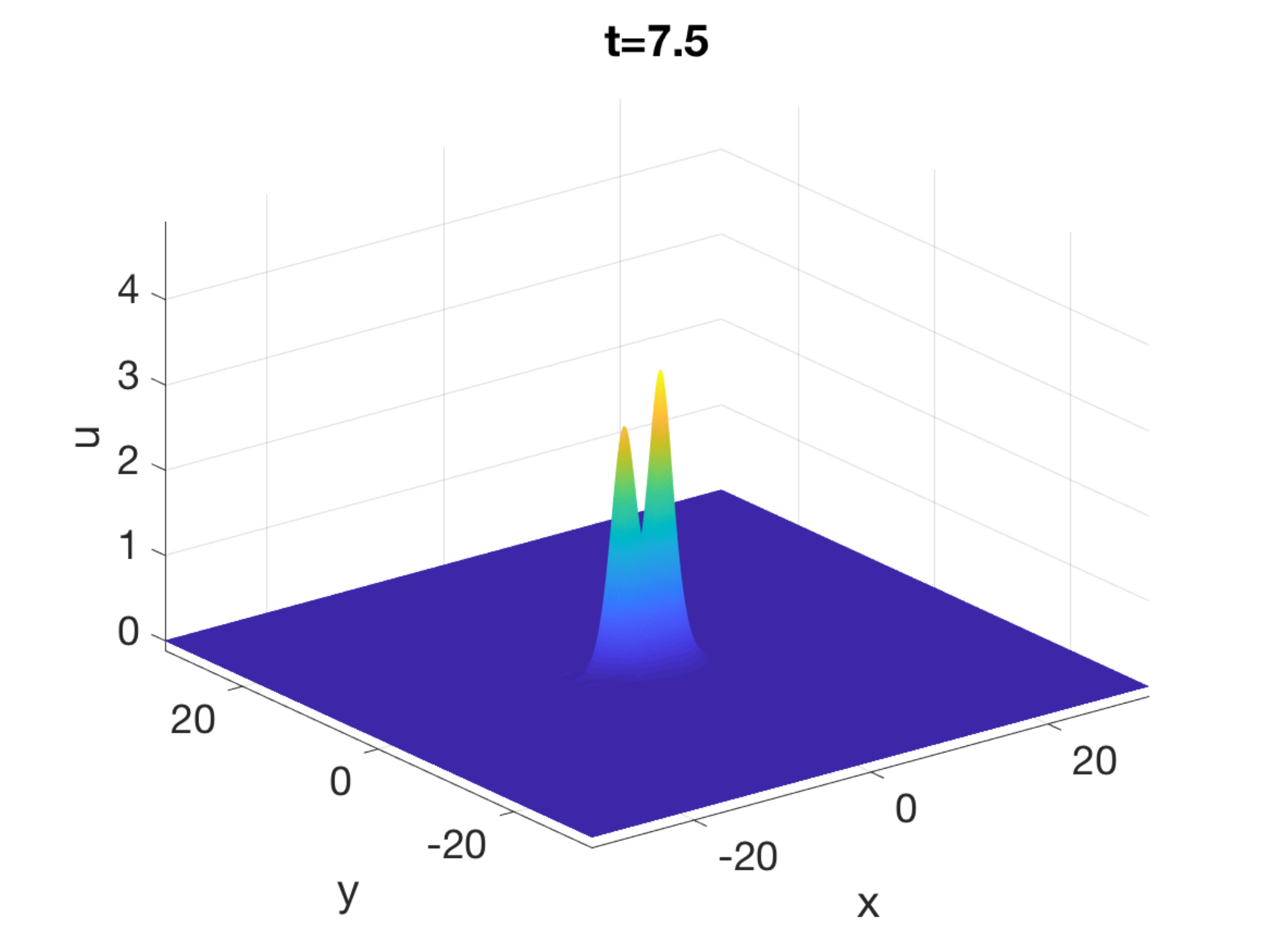}
 \includegraphics[width=0.49\hsize]{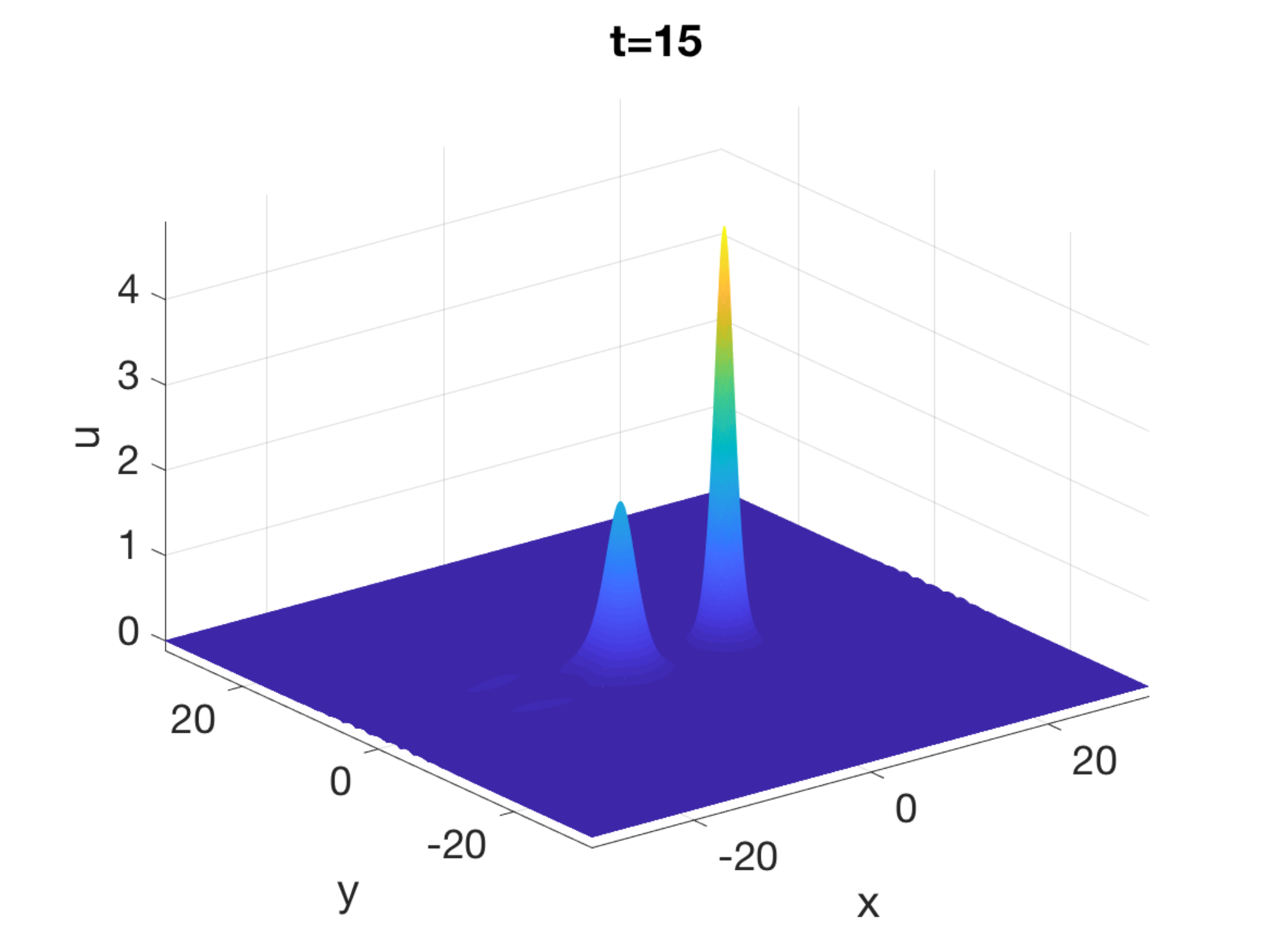}
\caption{Snapshots of the solution to \eqref{ZK} with the initial condition being a 
superposition of the soliton with $c=2$ centered at $x=-10$ 
and a soliton with $c=1$ centered at the origin. }
\label{ZKp22sol}
\end{figure}

As shown in Fig.~\ref{ZKp22sol} on the left of the first row, we 
consider the initial condition with two localized but somewhat separated solitons: one of them is the soliton \eqref{Qscal} with $c=2$ centered at $x=-10$ and another one is with $c=1$ centered at the origin. In our simulations we actually solve \eqref{Q} with $c=2$ to obtain the appropriate soliton, however, as an alternative, we could have used the scaling property \eqref{Qscal} to get the soliton with $c=2$. 

It can be seen in Fig.~\ref{ZKp22sol} that the faster soliton (note 
that we are still in a co-moving frame with \( c=1 \)) will hit the 
slower soliton around \( t=7 \). The collision is essentially elastic, 
the solitons appear to keep their shape after the collision.

\begin{figure}[!htb]
 \includegraphics[width=0.49\hsize]{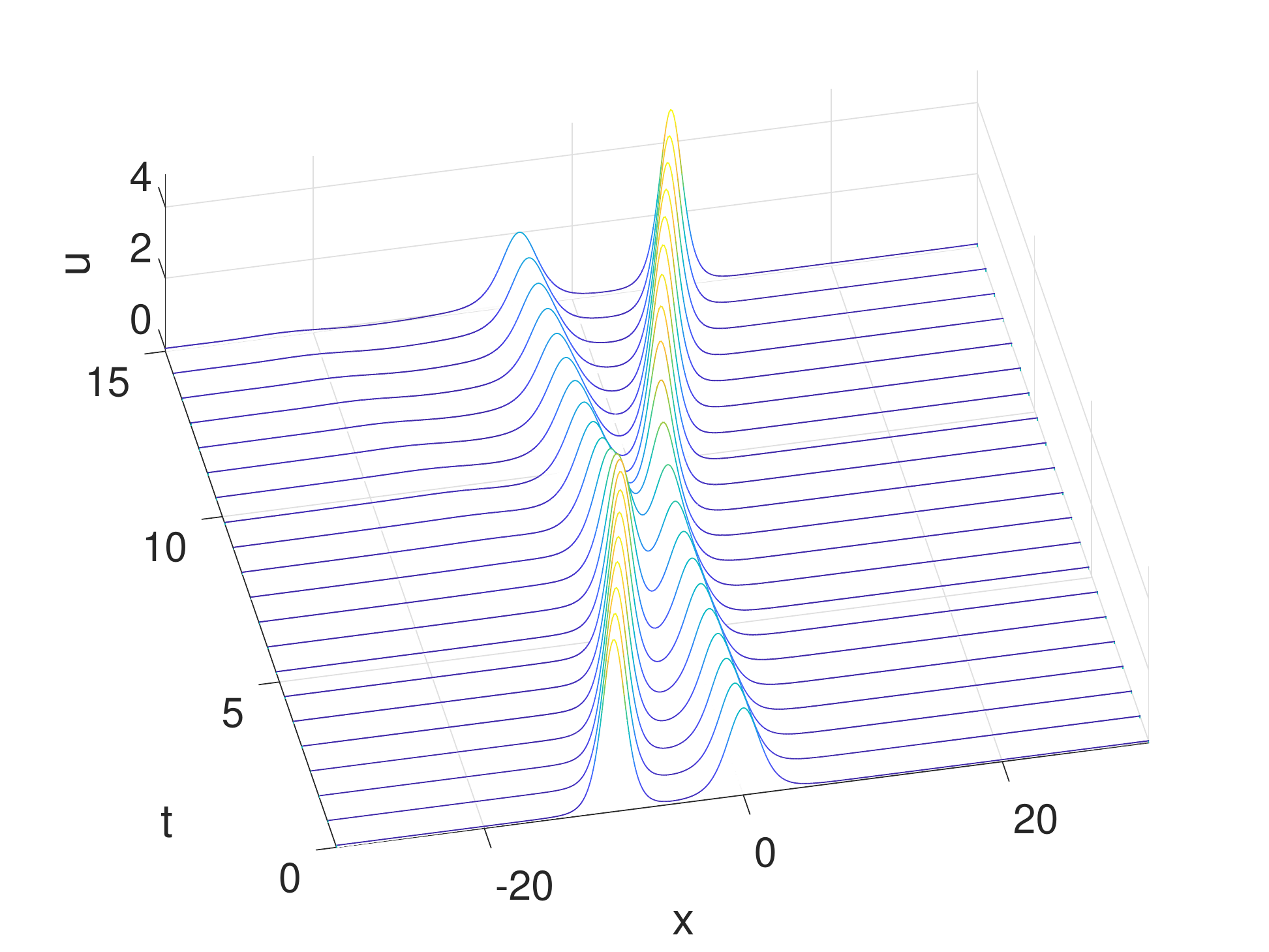}
 \includegraphics[width=0.49\hsize]{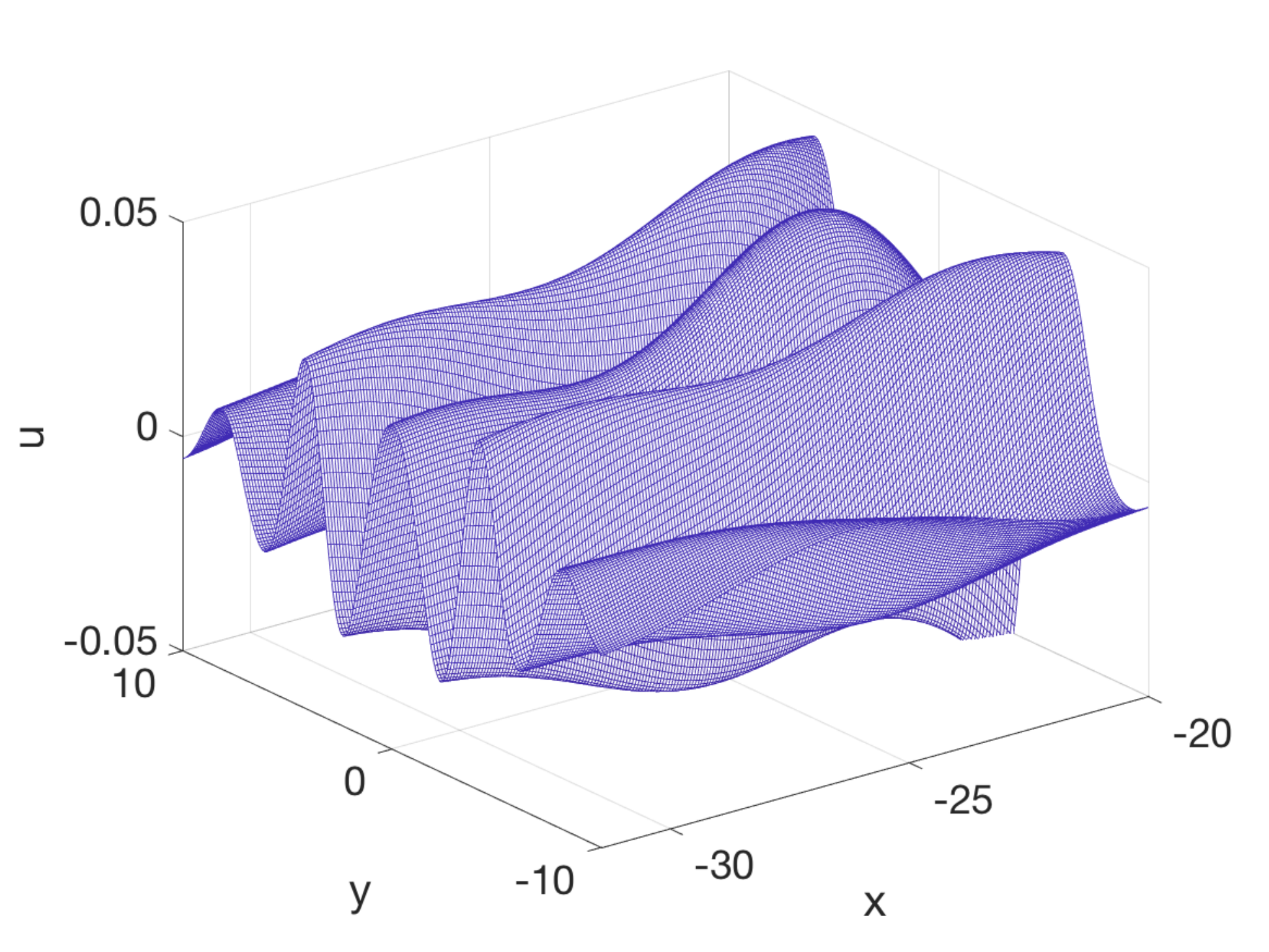}
\caption{Solution to \eqref{ZK} with $u(x,y,0)$ being the superposition of 
a soliton with $c=2$ centered at $x=-10$ and a soliton with $c=1$ centered 
at the origin on the $x$-axis for various times on the left, and a close-up of the bottom right subplot of Fig.~\ref{ZKp22sol} (at $t=15$) on the right. }
\label{ZKp22solx}
\end{figure}

The collision is best seen in a video (available online),  
 or on the $x$-axis as 
shown in Fig.~\ref{ZKp22solx}, since the motion is exclusively in the $x$-direction. In Fig.~\ref{ZKp22solx} we show the profile of this 2-soliton solution on the $x$-axis for various times. We note that the figure (except for the small radiation towards infinity) resembles 
closely the KdV 2-soliton.  However, the appearance of some radiation as seen from the close-up (right subplot) of the bottom right figure of Fig.~\ref{ZKp22sol},  shows  that the ZK 
equation is indeed not integrable.

\begin{figure}[!htb]
 \includegraphics[width=0.49\hsize]{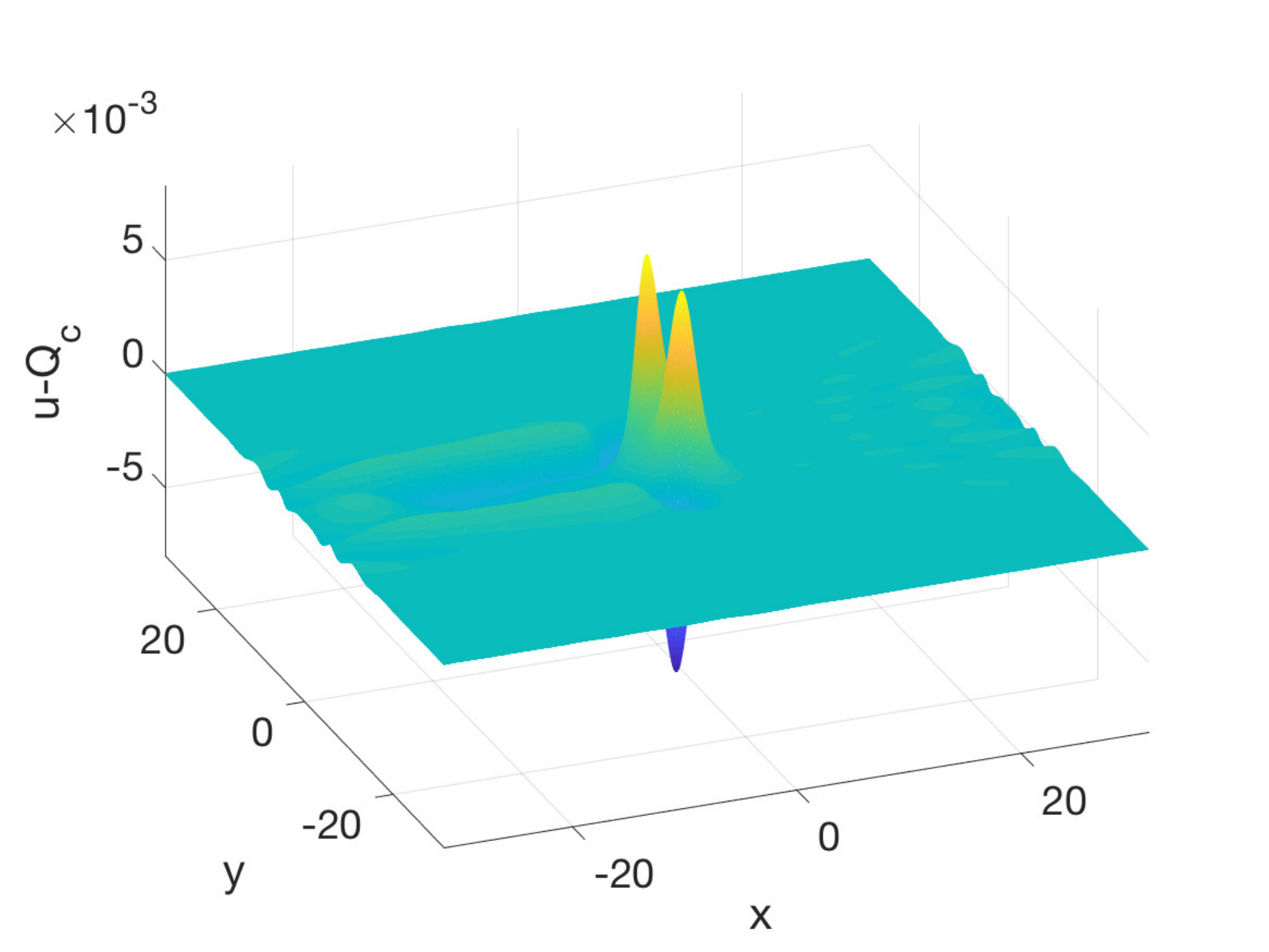}
 \includegraphics[width=0.49\hsize]{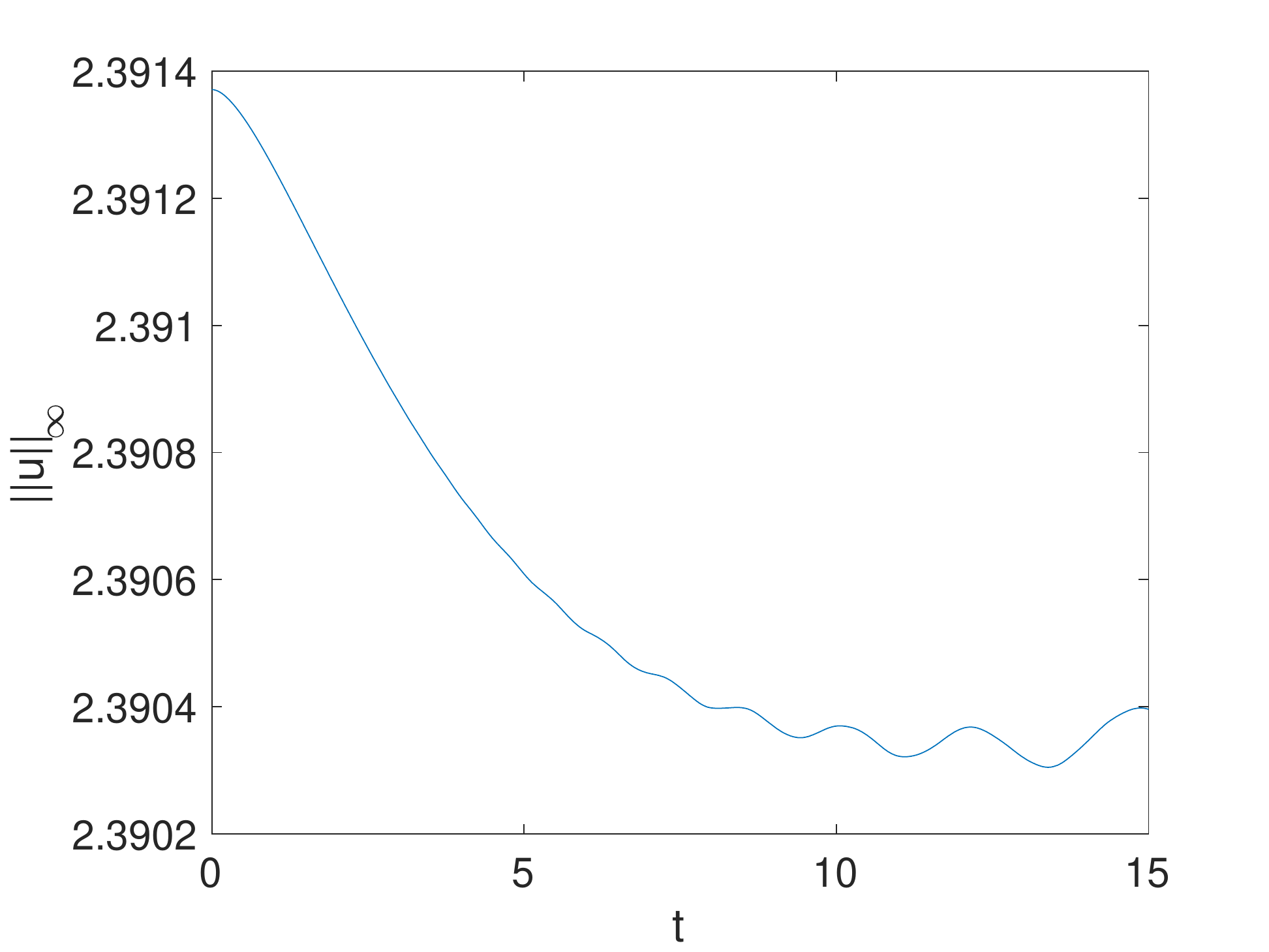}
\caption{Solution to \eqref{ZK} with $u(x,y,0) = 
Q(x,y-a)+Q(x,y+a)$ and $a=5$; on left the difference of the 
numerical solution at $t=15$ and the initial data, on the right the 
$L^{\infty}$ norm of the solution. }
\label{ZKp2x2d5}
\end{figure}

Next, we consider initial data of the form $u(x,y,0) = 
Q(x,y-a)+Q(x,y+a)$ with $a>0$ a constant in order  to study how 
the solitons interact if they have equal speeds, but are separated in the 
$y$-direction. In Fig.~\ref{ZKp2x2d5} we show the difference of the solution for $a=5$ 
at $t=15$ and the initial data in a co-moving frame with $c=1$. In 
this case the value of each soliton at the maximum of the other 
is on the order of $10^{-4}$. The interaction is therefore minimal, 
and the difference shown in Fig.~\ref{ZKp2x2d5} on the left is of the 
order of $10^{-3}$. 

\begin{figure}[!htb]
 \includegraphics[width=0.49\hsize]{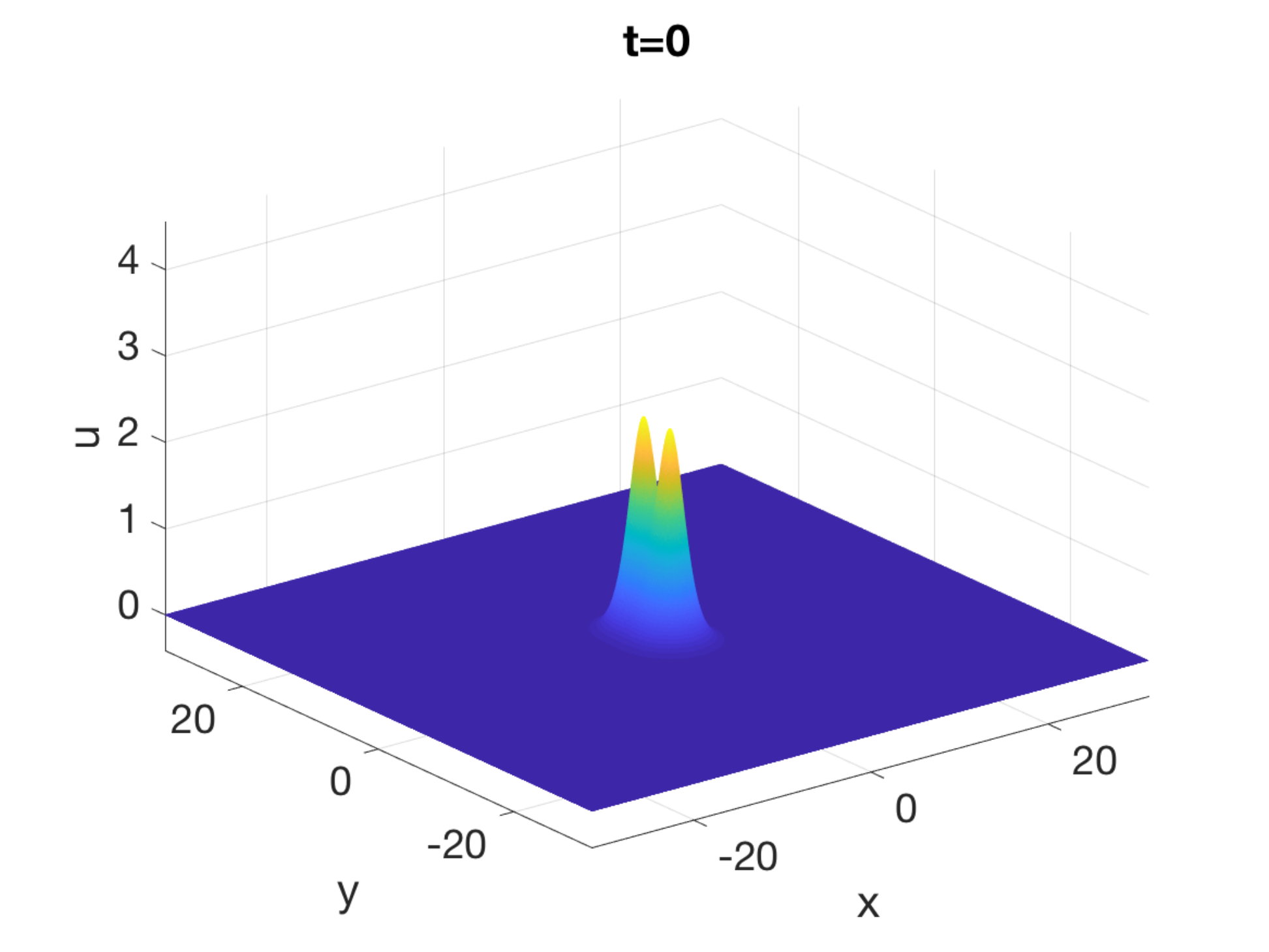}
 \includegraphics[width=0.49\hsize]{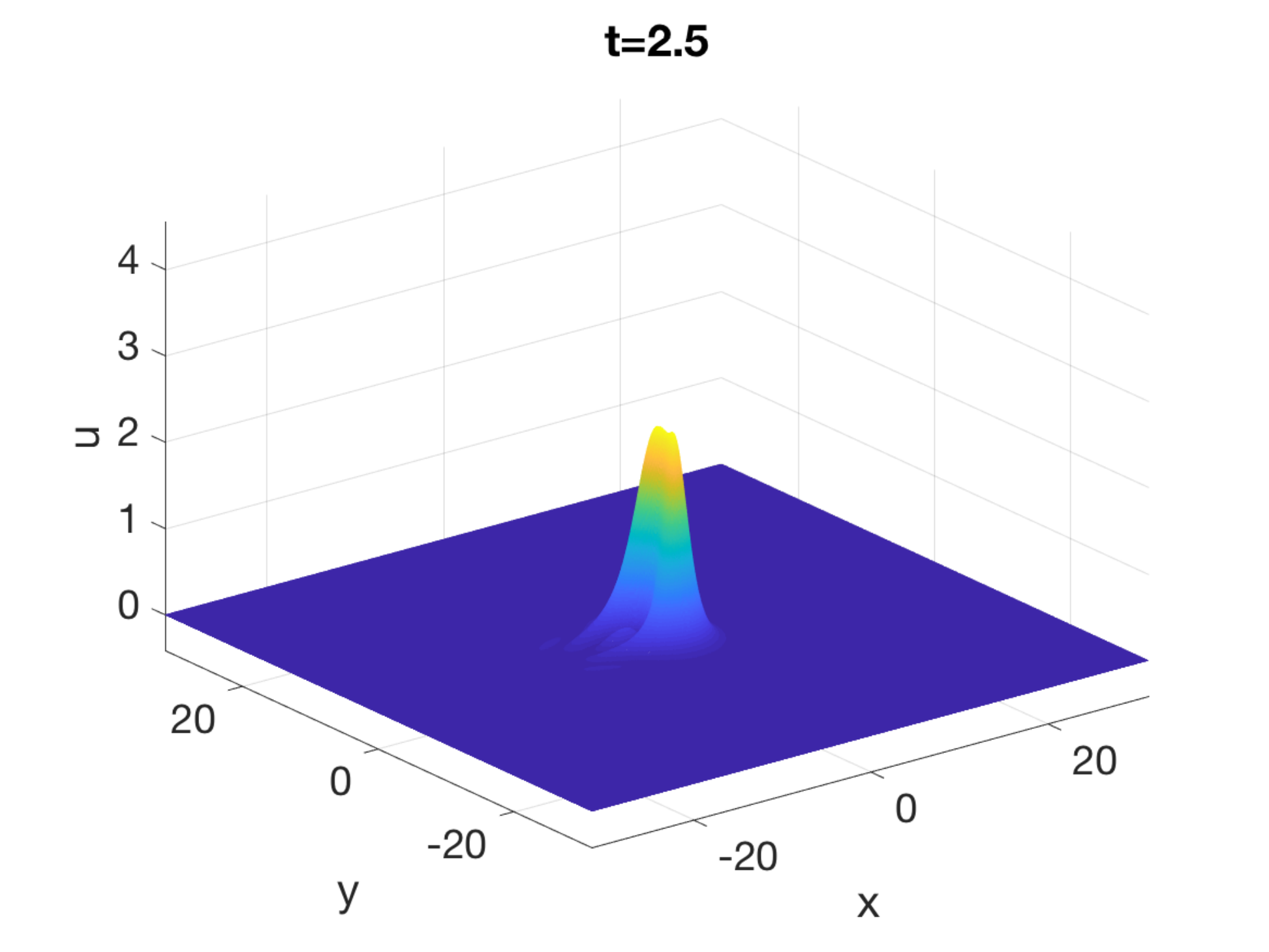}\\
 \includegraphics[width=0.49\hsize]{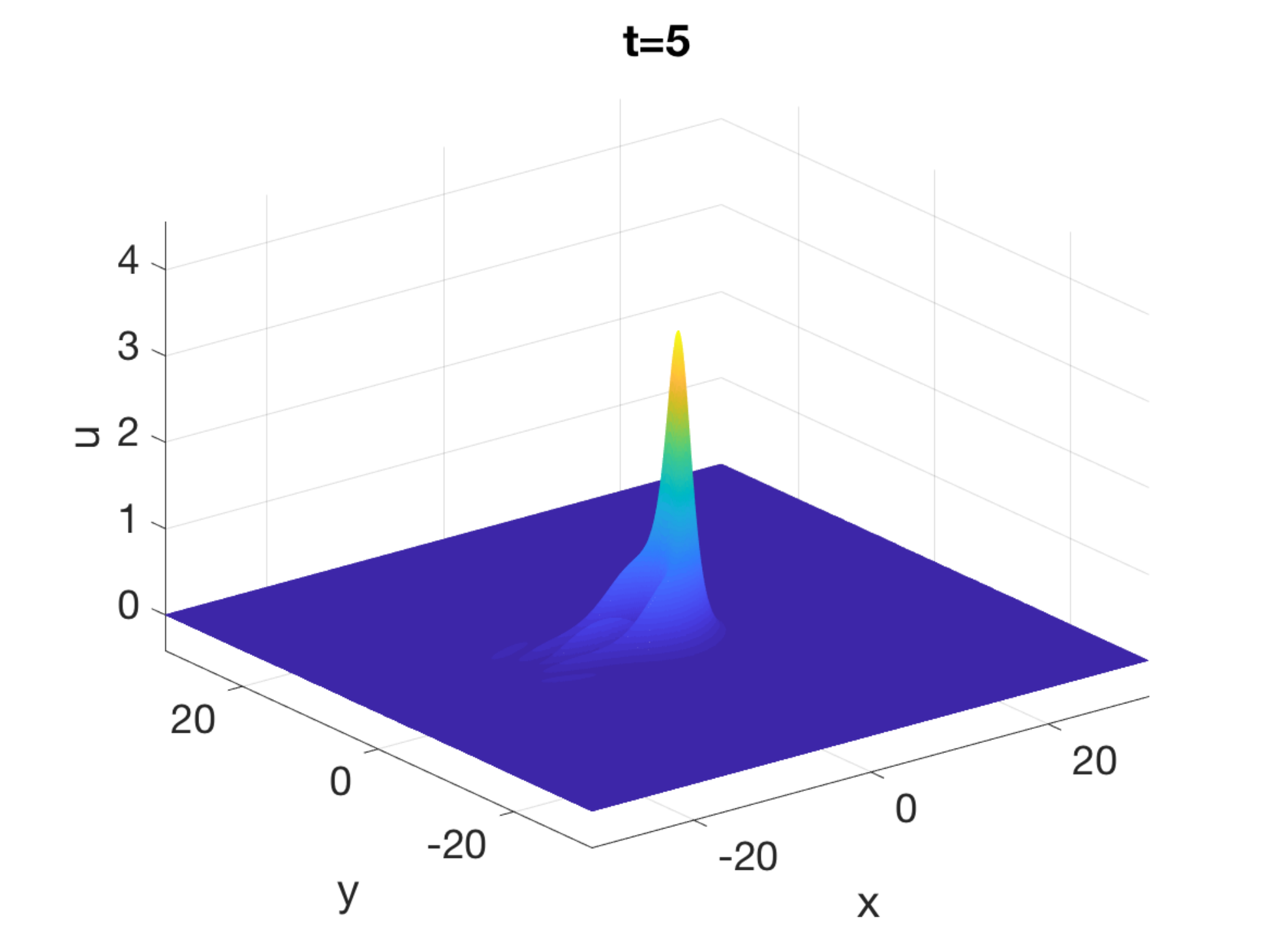}
 \includegraphics[width=0.49\hsize]{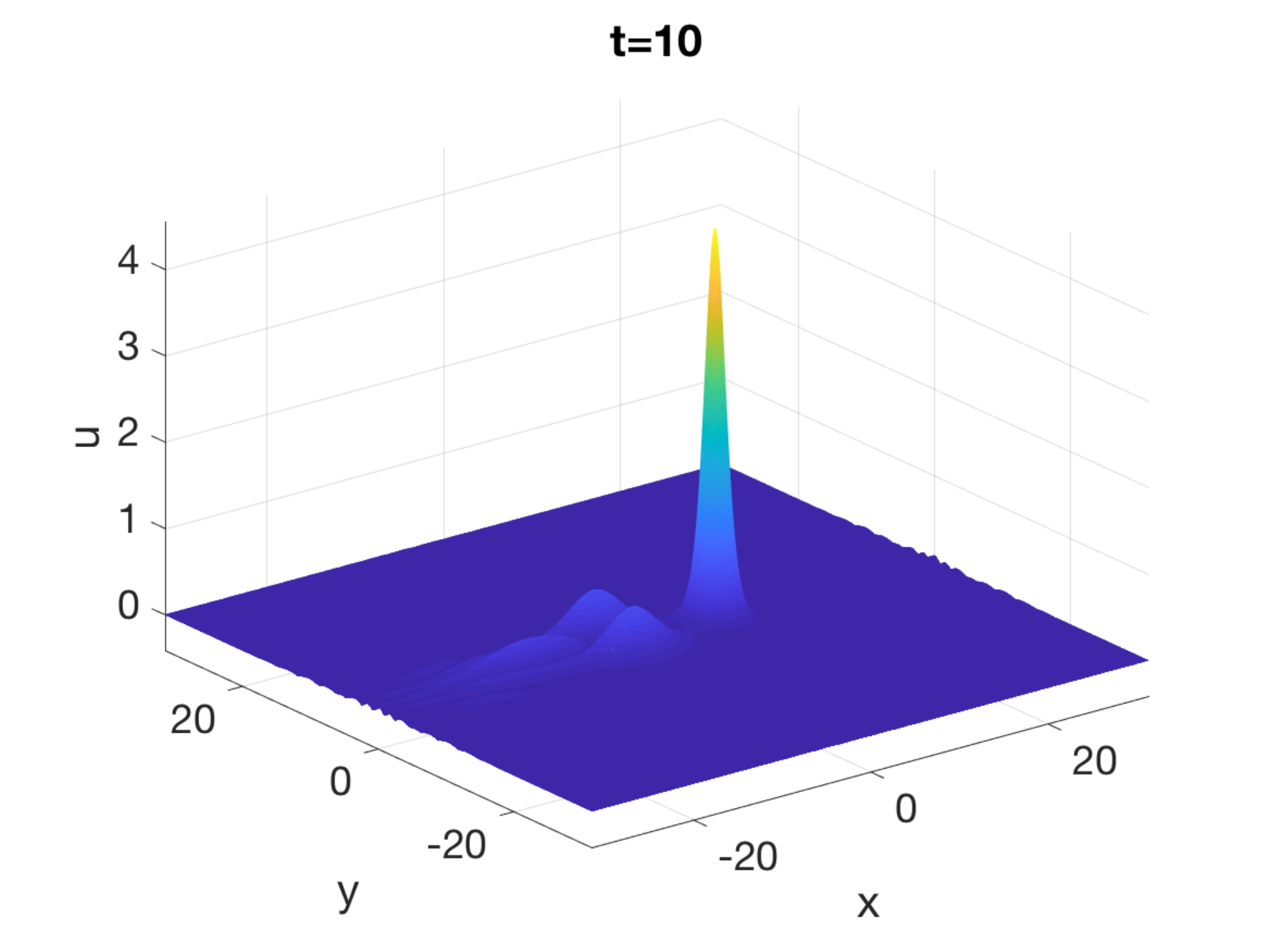}
\caption{Snapshots of the solution to \eqref{ZK} with $u(x,y,0) = 
Q(x,y-a)+Q(x,y+a)$ and $a=2$ at different times. }
\label{ZKp2x2d1}
\end{figure}

The snapshots of the solution for the same initial data with $a=2$ is shown in 
Fig.~\ref{ZKp2x2d1} at different times. 

\begin{figure}[!htb]
 \includegraphics[width=0.49\hsize]{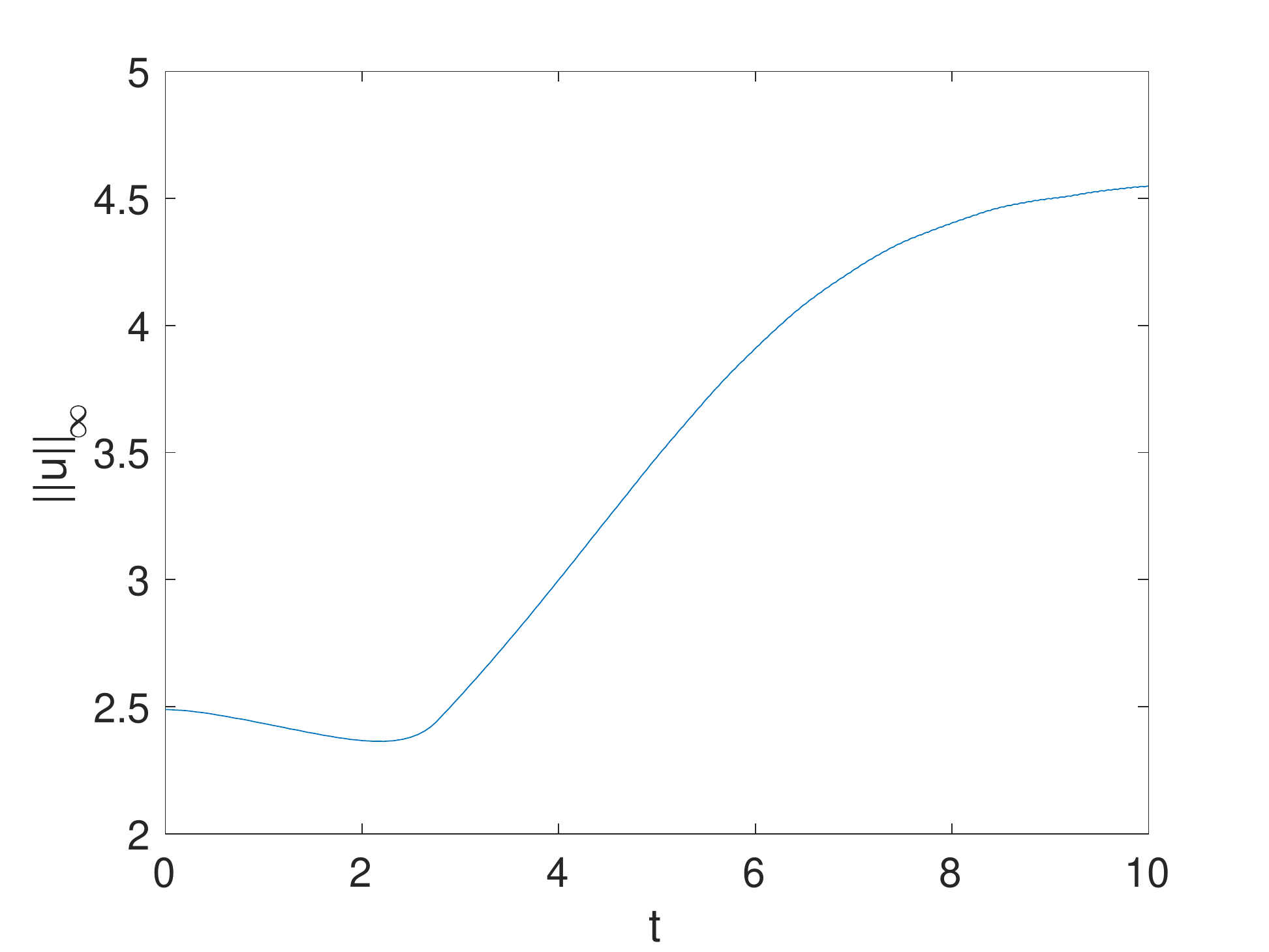}
 \includegraphics[width=0.49\hsize]{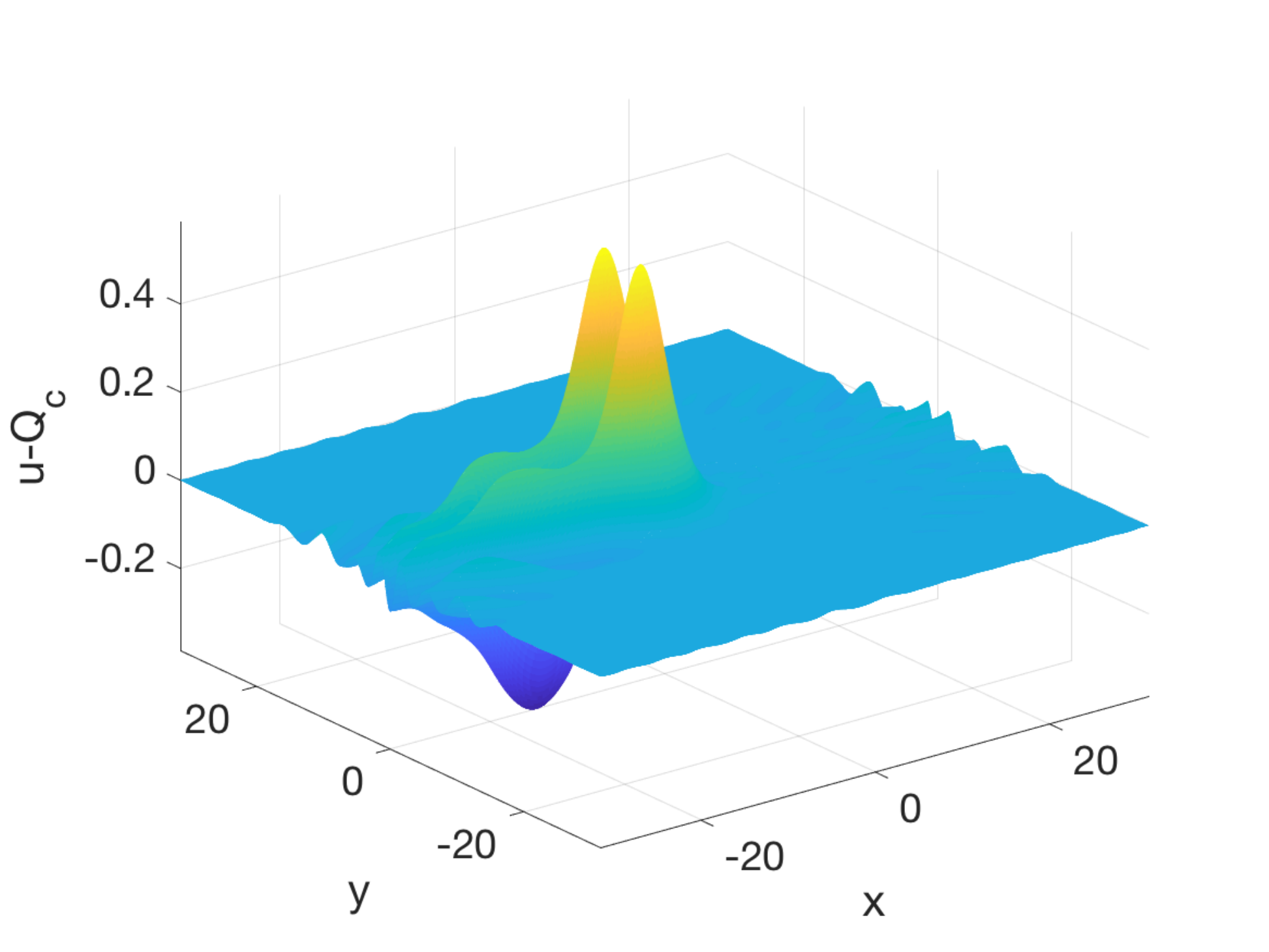}
\caption{$L^{\infty}$ norm of the solution to \eqref{ZK} with $u(x,y,0) = 
Q(x,y-a)+Q(x,y+a)$ and $a=2$ on the left; the difference between this solution at
$t=10$ and a fitted soliton solution of \eqref{Qscal} on the right. }
\label{ZKp2x2d1inf}
\end{figure}

It appears that the final state of the solution is a single soliton. This is indicated both by 
the $L^{\infty}$ norm of the solution on the left subplot of Fig.~\ref{ZKp2x2d1inf} and by the difference with a fitted soliton solution of \eqref{Qscal} on the right subplot. 

\begin{figure}[!htb]
 \includegraphics[width=0.49\hsize]{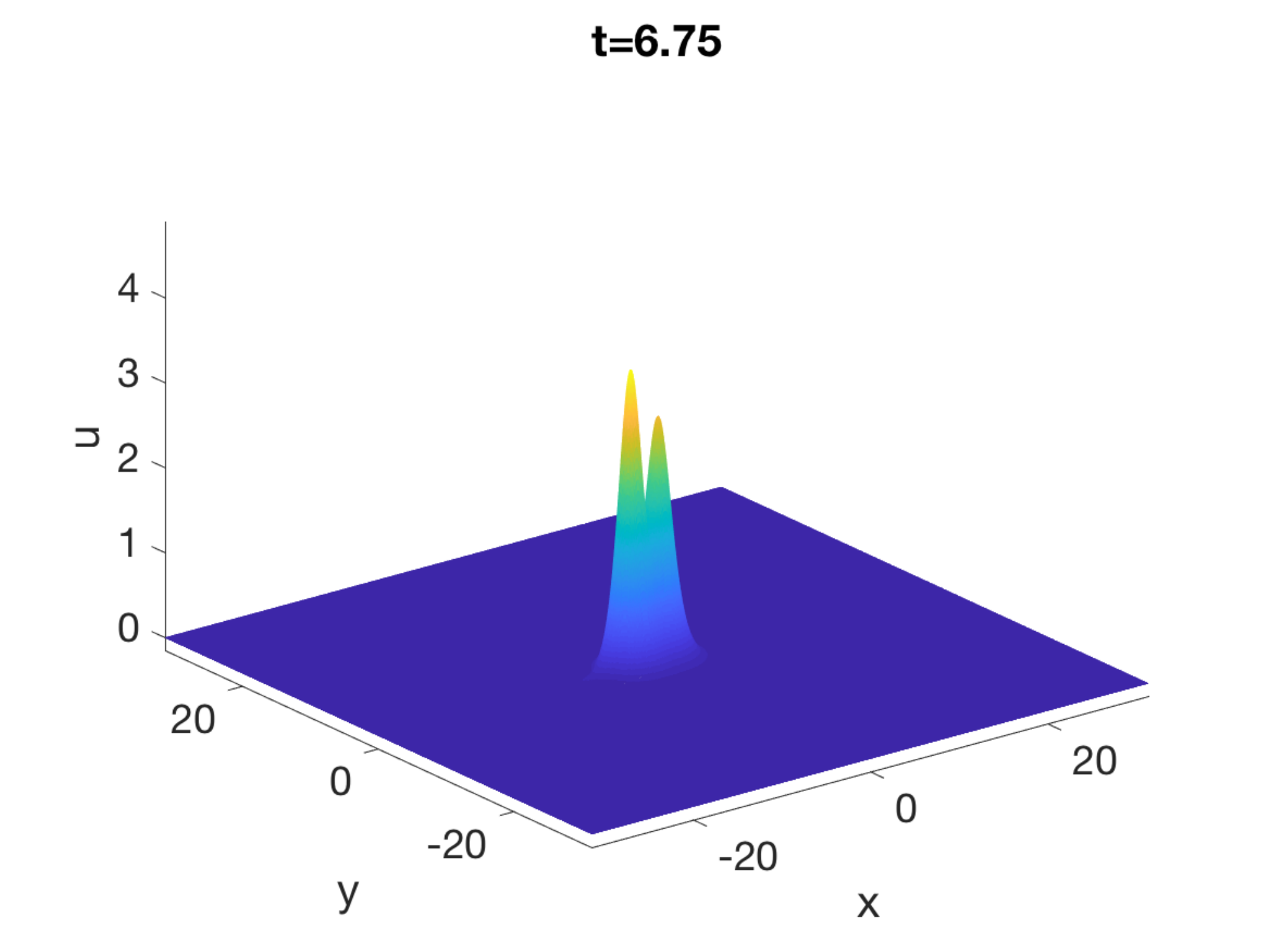}
 \includegraphics[width=0.49\hsize]{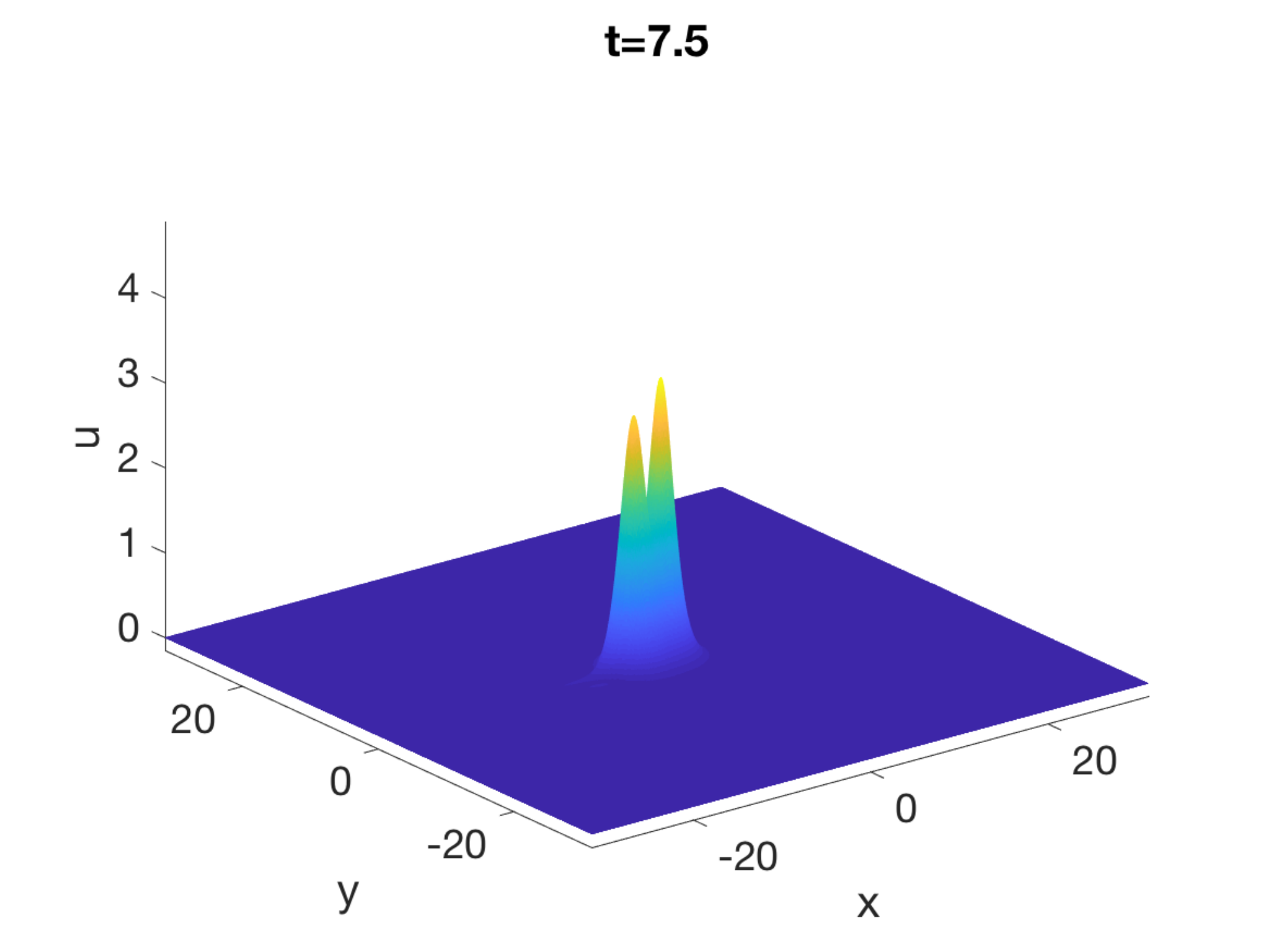}
\caption{Solution to \eqref{ZK} with $u(x,y,0) = 
Q_{2}(x+10,y+1)+Q(x,y)$ slightly before and after the collision of the 
solitons. }
\label{ZKp221}
\end{figure}

We next note that if we consider a slightly off-centered collision of solitons, that is, the solution with the initial condition $u(x,y,0)=Q_{2}(x+10,y+1)+Q(x,y)$, we 
get a very similar behavior to the collision in Fig.~\ref{ZKp22solx}, 
see snapshots in Fig.~\ref{ZKp221}. After the 
interaction, the larger soliton, which was below the 
smaller one in $y$-direction, will be above it (they essentially change roles in the 
elastic collision). 

We thus conclude that it is the separation distance, not the relative 
location in the plane, of solitons that influences their interaction in the long run. This dependence is of significant interest, but will be investigated elsewhere.

\subsection{Soliton resolution}
Since ZK solitons in the subcritical case $p=2$ are clearly stable (from what we simulated both orbitally and asymptotically), and furthermore, they even show essentially elastic collisions, one would expect that, according to the soliton resolution 
conjecture, solitons plus radiation appear in the long term evolution of localized initial data with sufficient mass. 
For the following computations in this subsection, we no longer use  co-moving frames. 

\begin{figure}[!htb]
 \includegraphics[width=0.7\hsize]{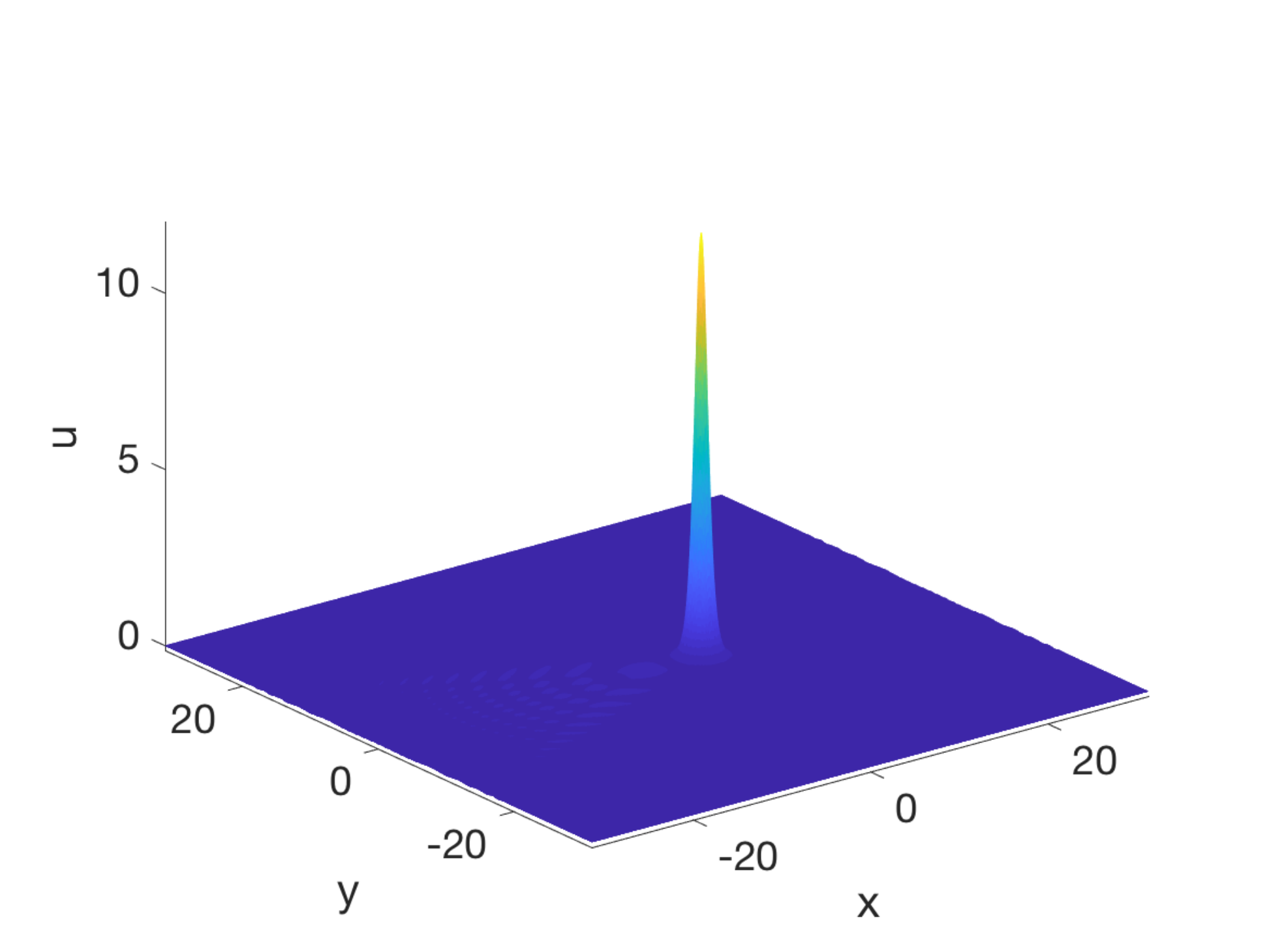}
\caption{Solution to \eqref{ZK} with $u(x,y,0)=10 \, 
e^{-(x^{2}+y^{2})}$ at $t=1$. Note the radiation emitted to the left with the total angle of $60^0$.}
\label{ZKp210Gausst1}
\end{figure}

In Fig.~\ref{ZKp210Gausst1} we show the ZK evolution (at $t=1$) of 
the Gaussian initial condition $u(x,y,0)=10 \, e^{-(x^{2}+y^{2})}$. 
It appears that a single soliton emerges from the initial bump plus some radiation. We note that the radiation is emitted to the left of the $x$-axis up to an angle of $30^0$ with the negative $x$-axis (so the total opening is $60^0$), this is in confirmation of the results in \cite{CMPS}.

\begin{figure}[!htb]
 \includegraphics[width=0.49\hsize]{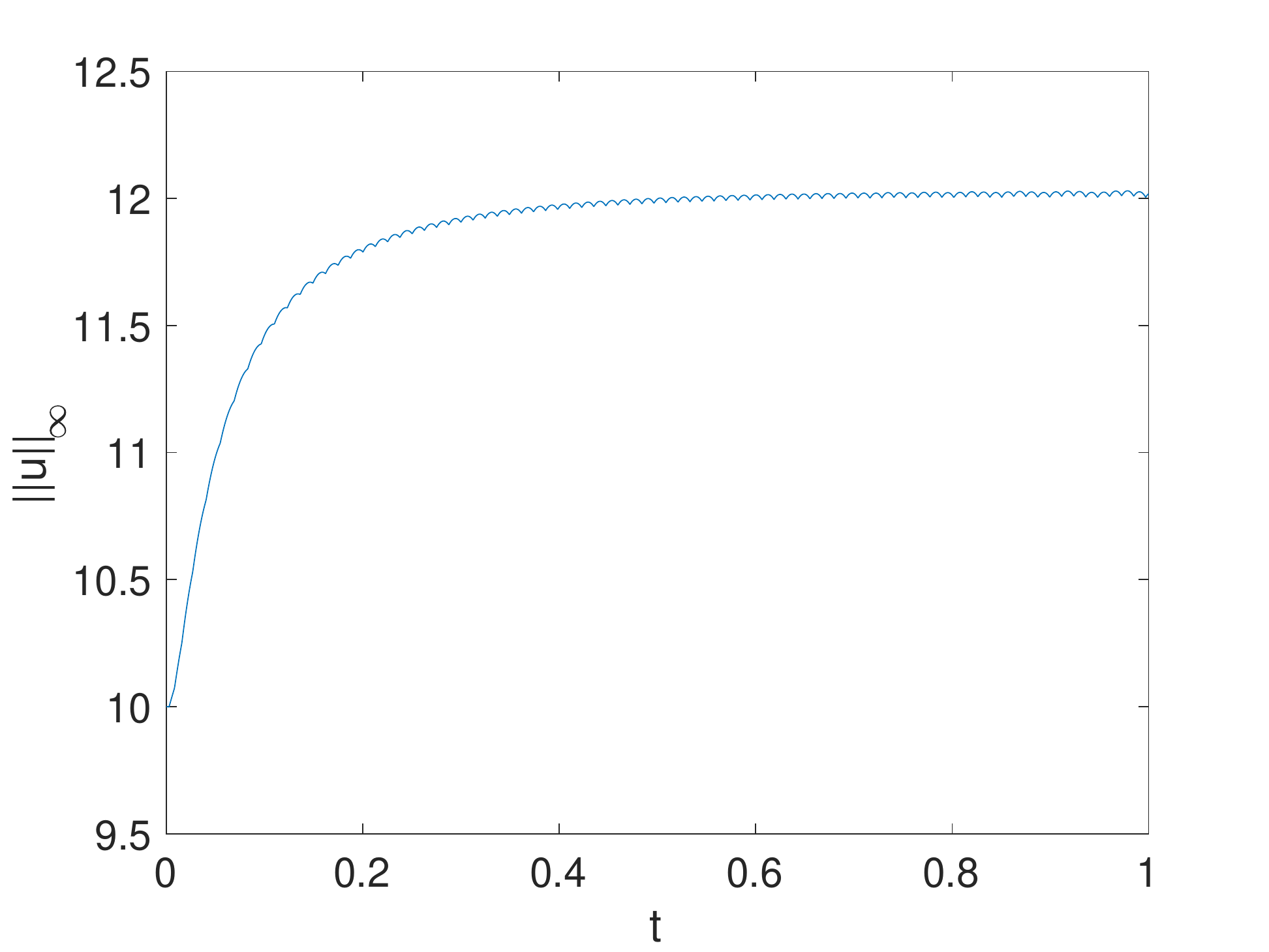}
 \includegraphics[width=0.49\hsize]{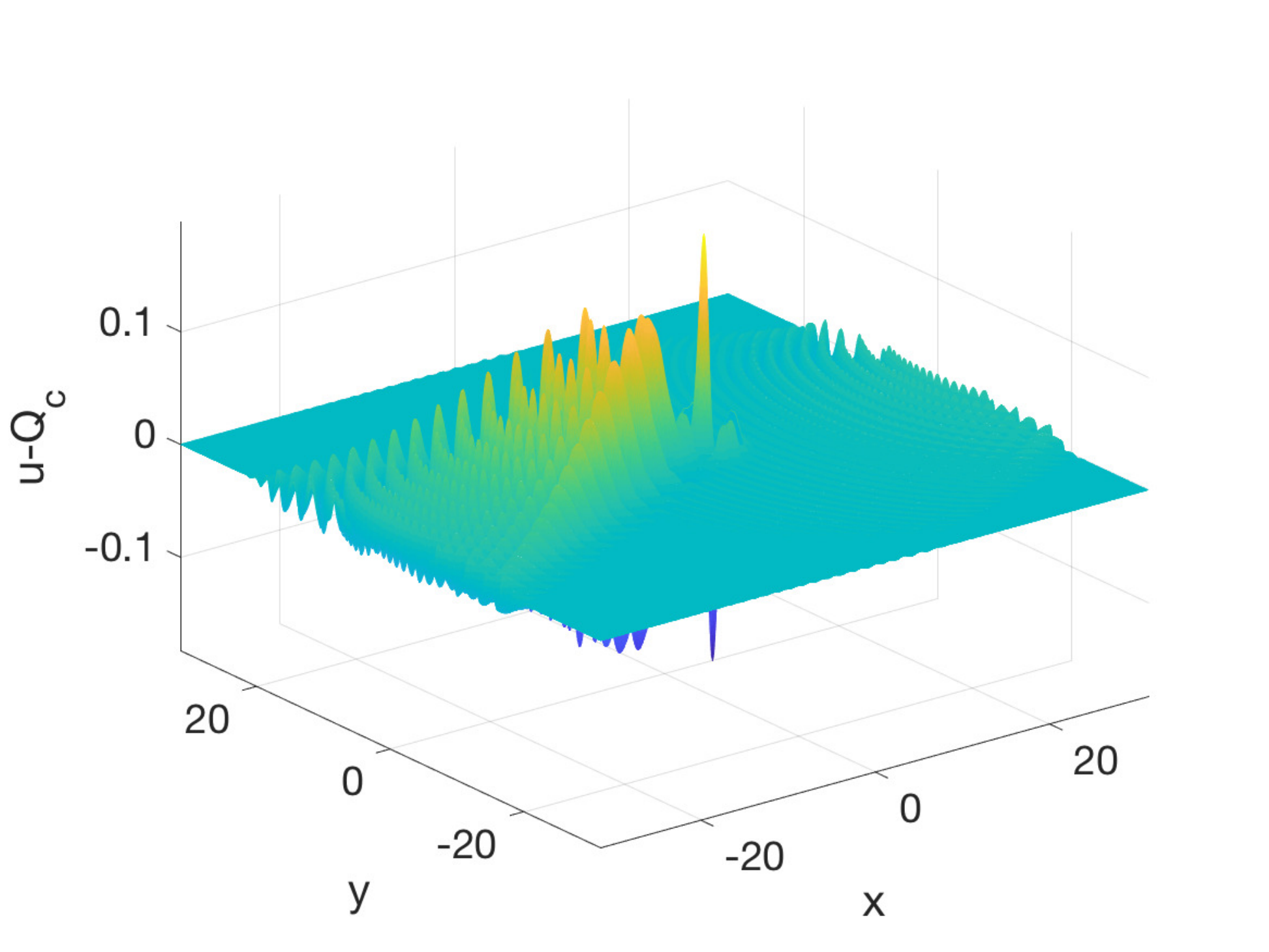}
\caption{Solution to \eqref{ZK} with $u(x,y,0)=10 \, e^{-(x^{2}+y^{2})}$: 
on the left the time dependence of the $L^{\infty}$ norm of the solution; 
on the right the difference between the solution at $t=1$ and a fitted soliton. }
\label{ZKp210Gaussmax}
\end{figure}

The $L^{\infty}$ norm of the solution shown on the left in Fig.~\ref{ZKp210Gaussmax} also indicates that a soliton appears. The difference between the numerical solution at $t=1$ and 
a fitted soliton solution of \eqref{Qscal} can be seen on the right of the same figure. It indicates that a soliton appears, but that the final state has not yet been reached, which is also clear from the presence of radiation in the figure. 

\begin{figure}[!htb]
 \includegraphics[width=0.49\hsize]{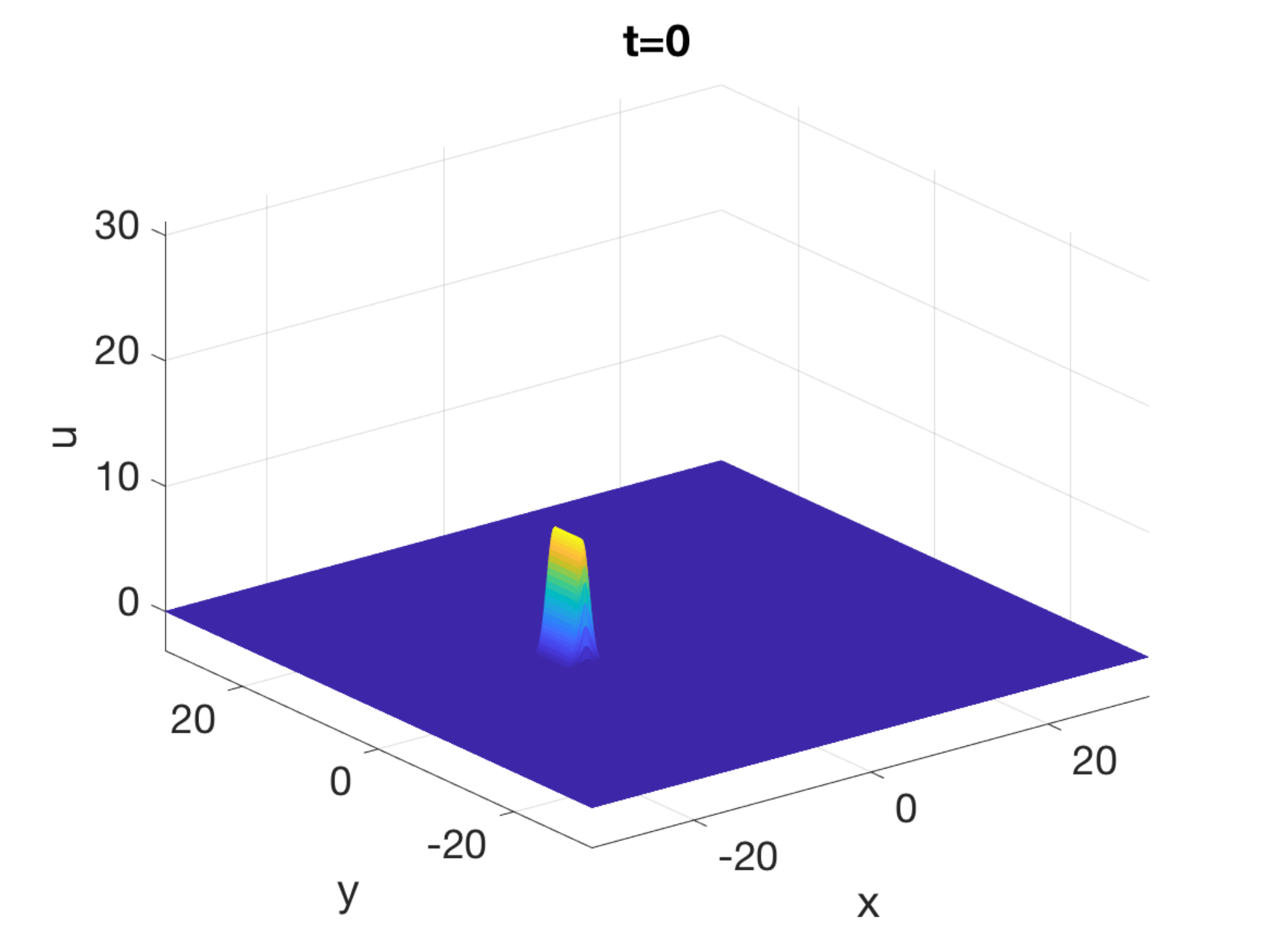}
 \includegraphics[width=0.49\hsize]{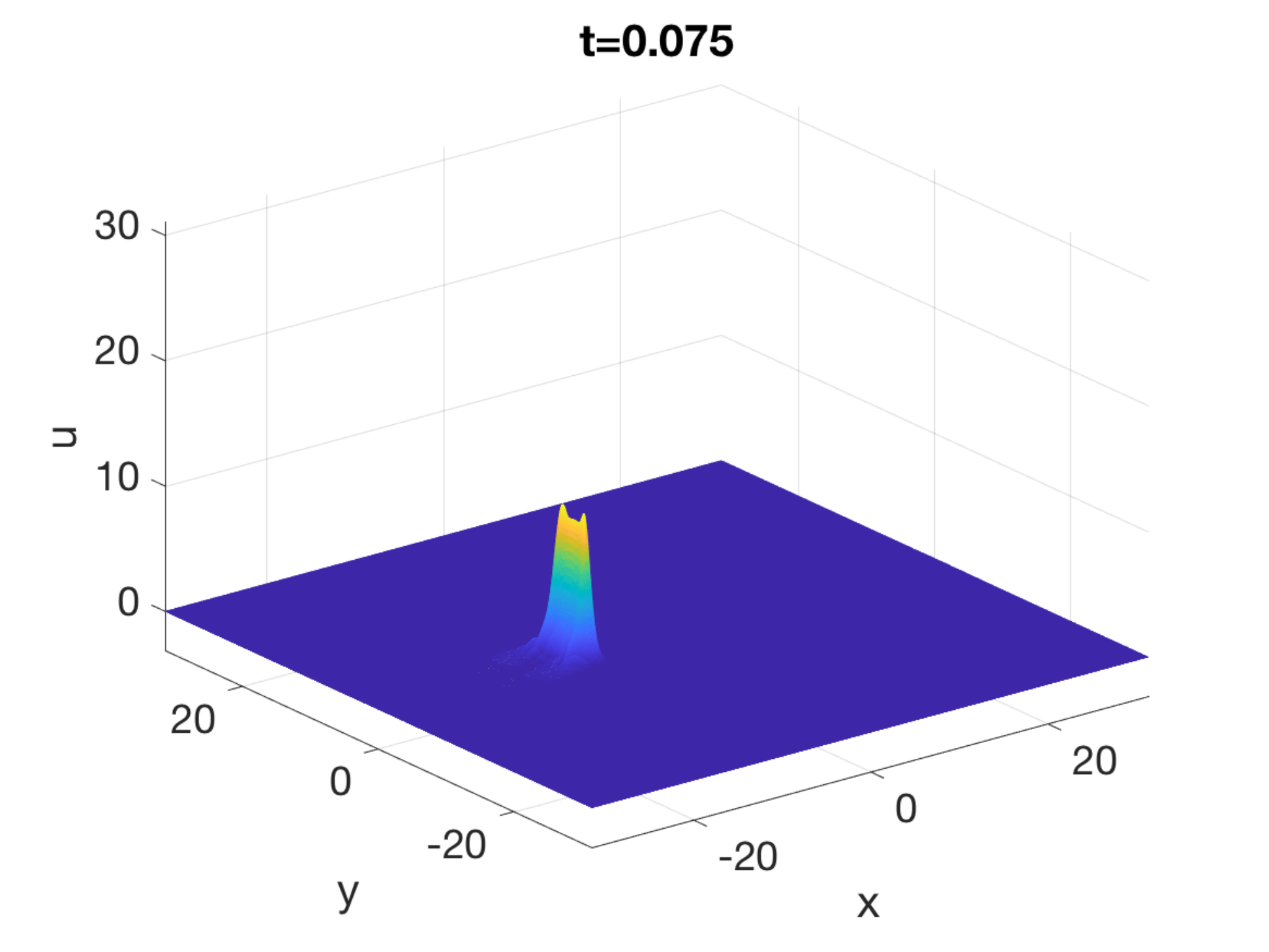}\\
 \includegraphics[width=0.49\hsize]{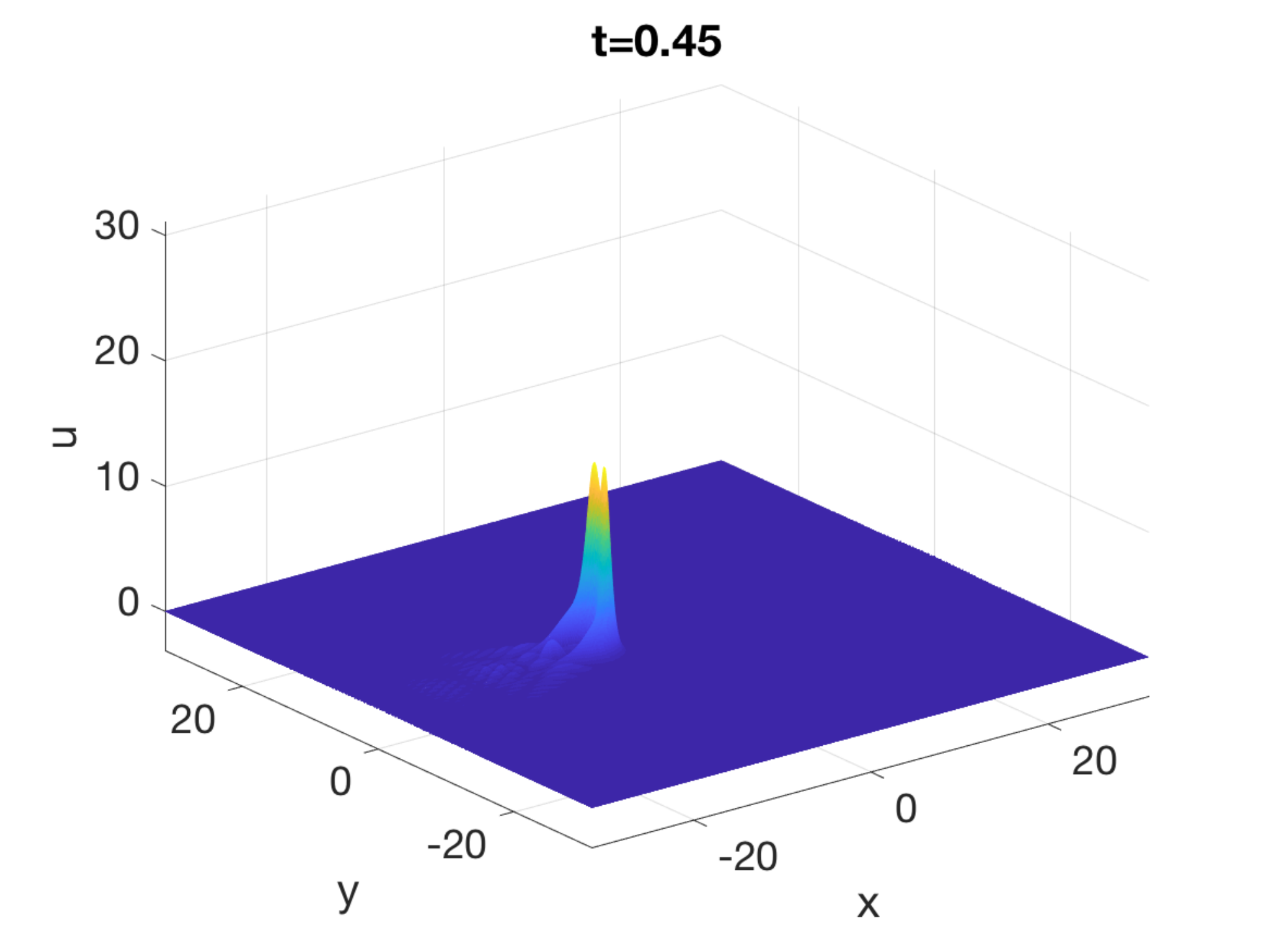}
 \includegraphics[width=0.49\hsize]{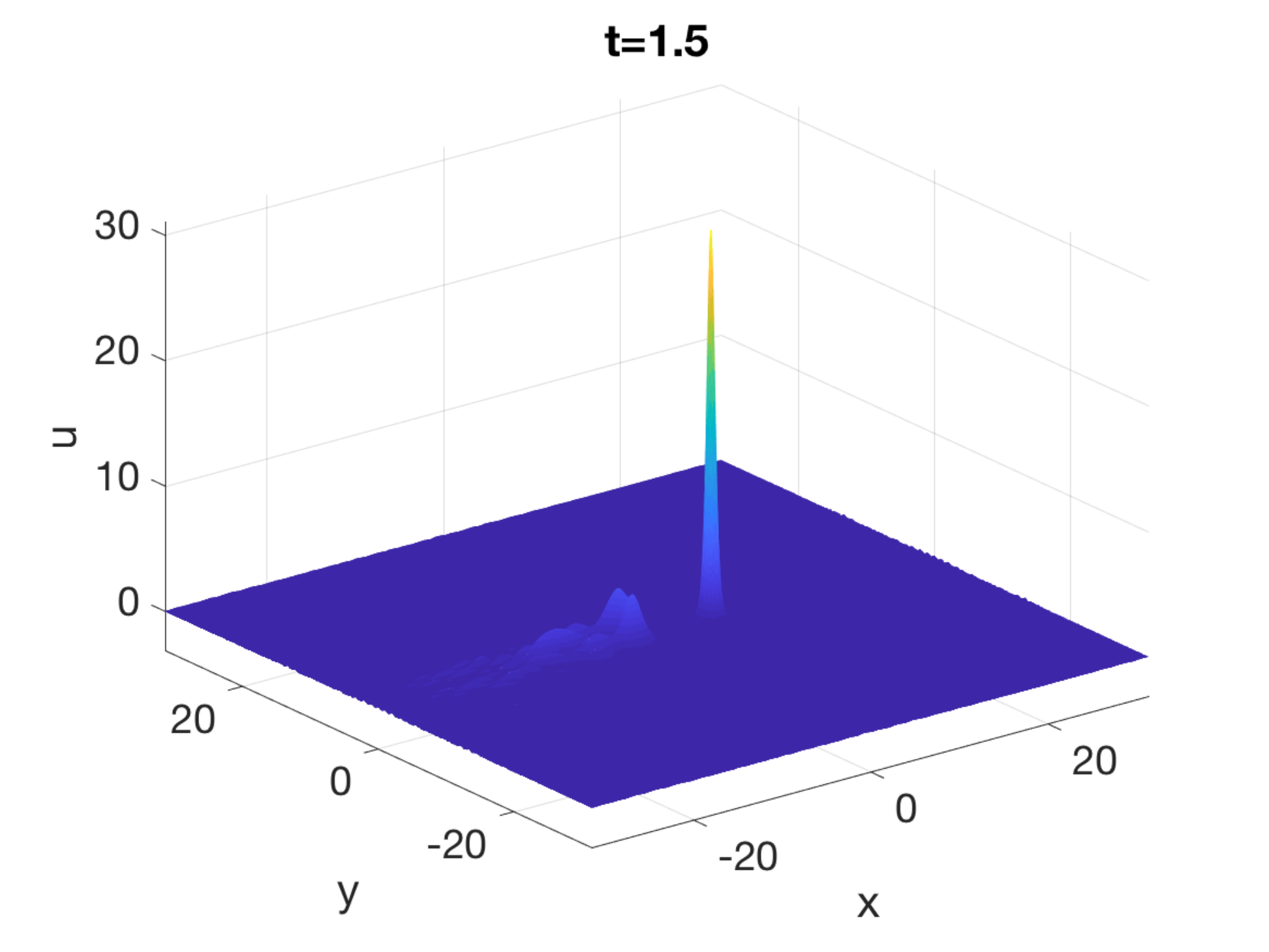}
\caption{Time snapshots of the solution to \eqref{ZK} with the wall-type initial condition \eqref{wallcond}.}
\label{ZKwall}
\end{figure}

Since we have seen in Fig.~\ref{ZKp2x2d1} that nearby solitons tend to 
merge into a single soliton, it is not surprising that the same is 
found for nearby general bumps. Therefore, it is interesting to study 
initial data with an extended maximal region, for instance, the wall-like structure
\begin{equation}
    u(x,y,0) = 
    \begin{cases}
        10 \, e^{-x^{2}} & |x|\leq 1.5 \\
         10 \, e^{-(x^{2}+(y-1.5)^{8})} & x>1.5 \\
         10 \, e^{-(x^{2}+(y+1.5)^{8})}& x<-1.5.
    \end{cases}
    \label{wallcond}
\end{equation}
The snapshots at different times of the corresponding ZK solution are 
given in Fig.~\ref{ZKwall}. The initial wall develops two peaked 
structures near the edges, which then merge into one large bump 
travelling to the right, and radiation (and possibly forming more 
smaller solitons, travelling slowly behind). 

\begin{figure}[!htb]
 \includegraphics[width=0.49\hsize]{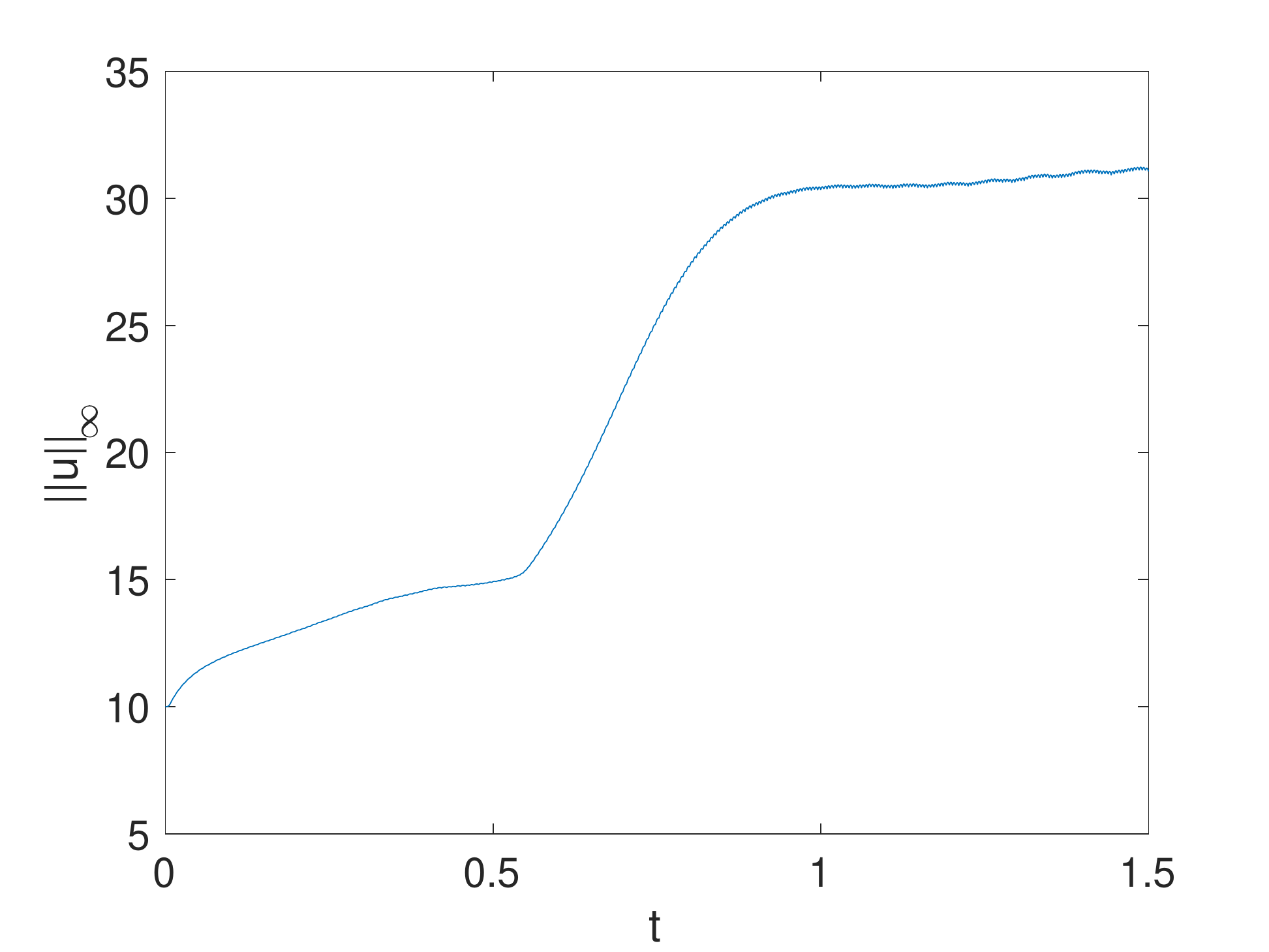}
 \includegraphics[width=0.49\hsize]{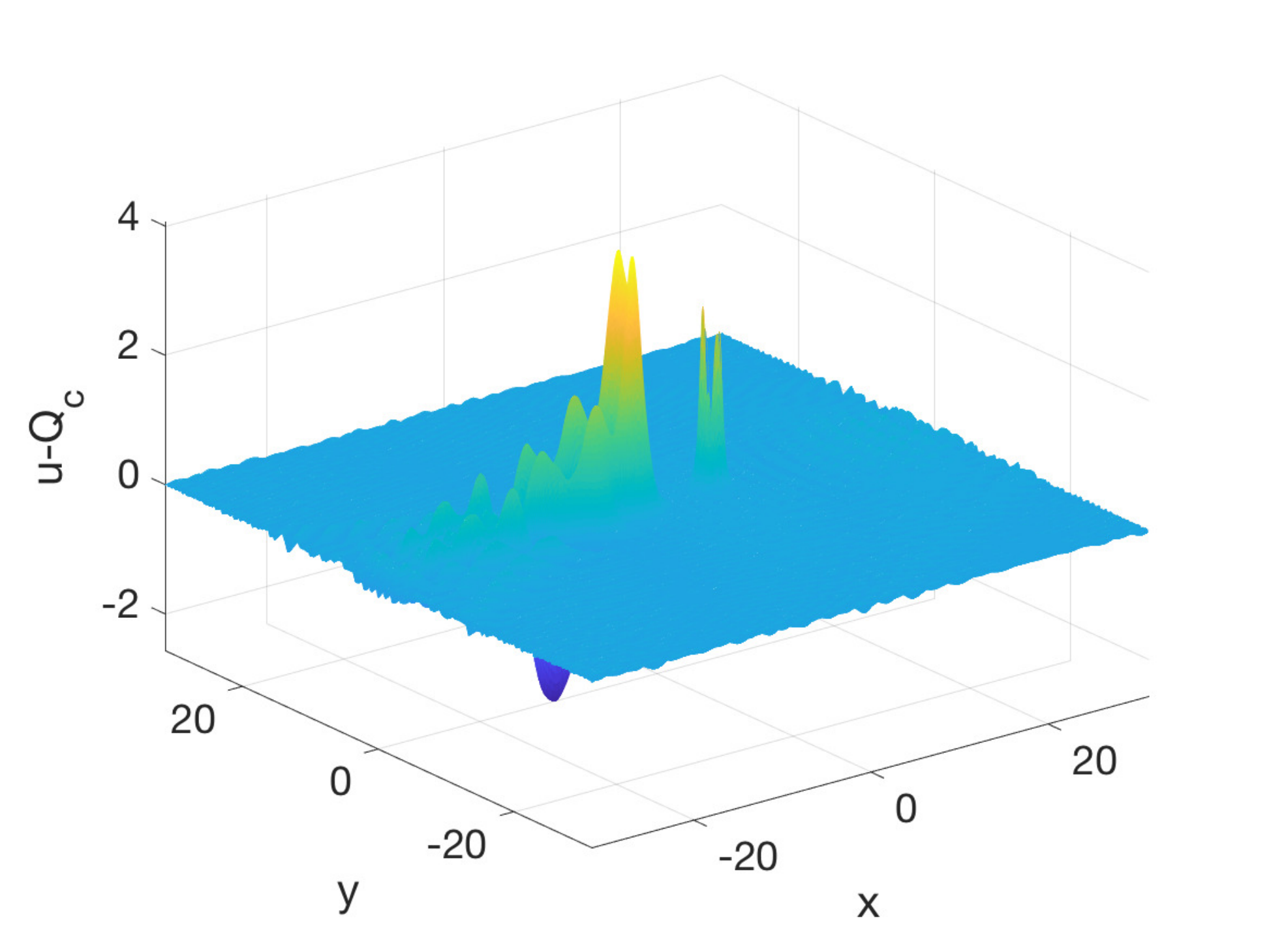}
\caption{Solution to \eqref{ZK} with $u_0$ in \eqref{wallcond}: 
on the left time dependence of the $L^{\infty}$ norm; 
on the right the difference between the solution at $t=1.5$ and a fitted soliton.}
\label{ZKwallmax}
\end{figure}

The $L^{\infty}$ norm on the left of Fig.~\ref{ZKwallmax} appears to be still slightly growing, which indicates that the final state of the main bump is not yet reached. Note, however, the difference to a fitted soliton solution of \eqref{Qscal} makes a plausible conclusion that this final state should indeed be a soliton. 

\begin{figure}[!htb]
 \includegraphics[width=0.49\hsize]{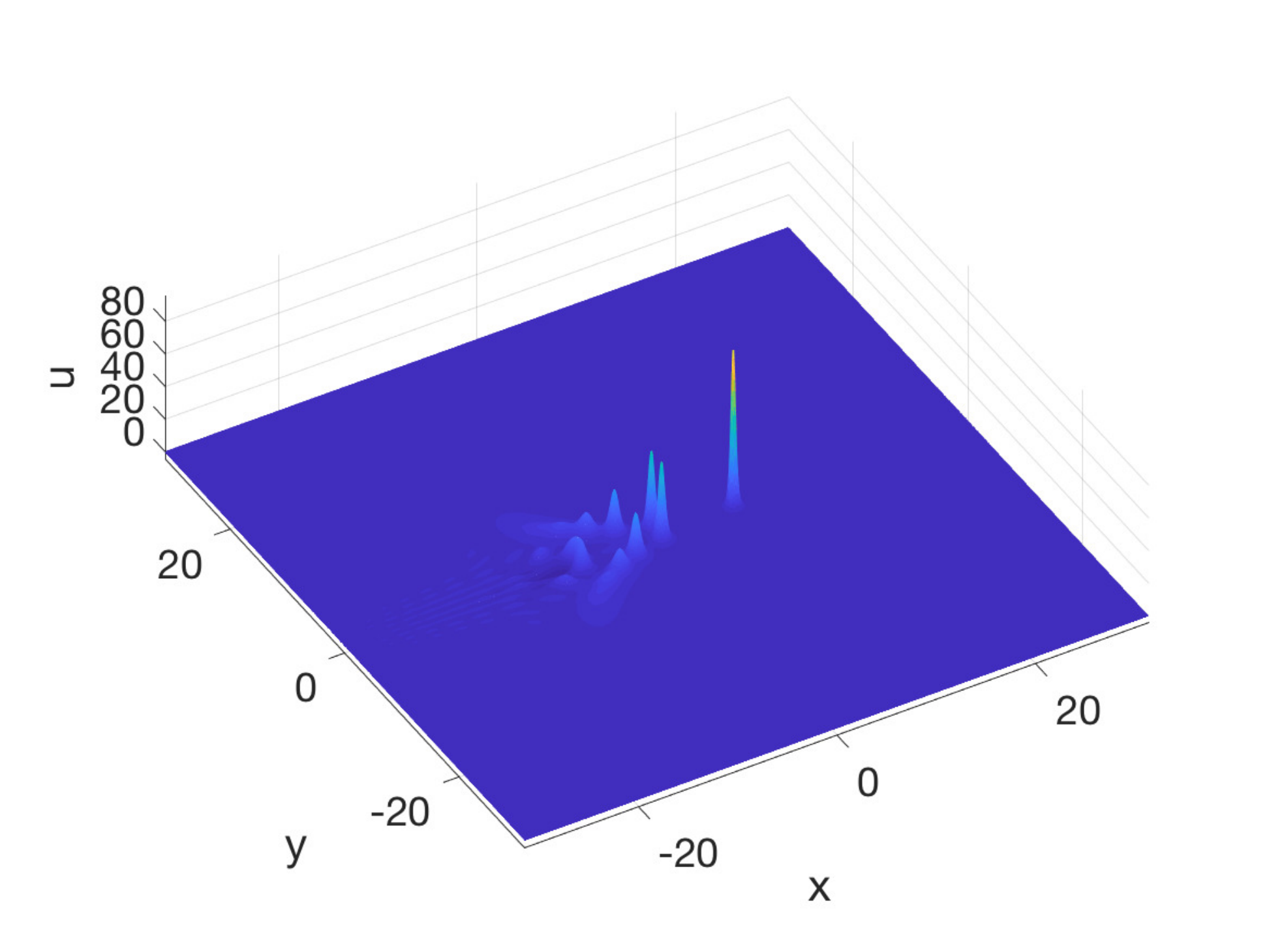}
 \includegraphics[width=0.49\hsize]{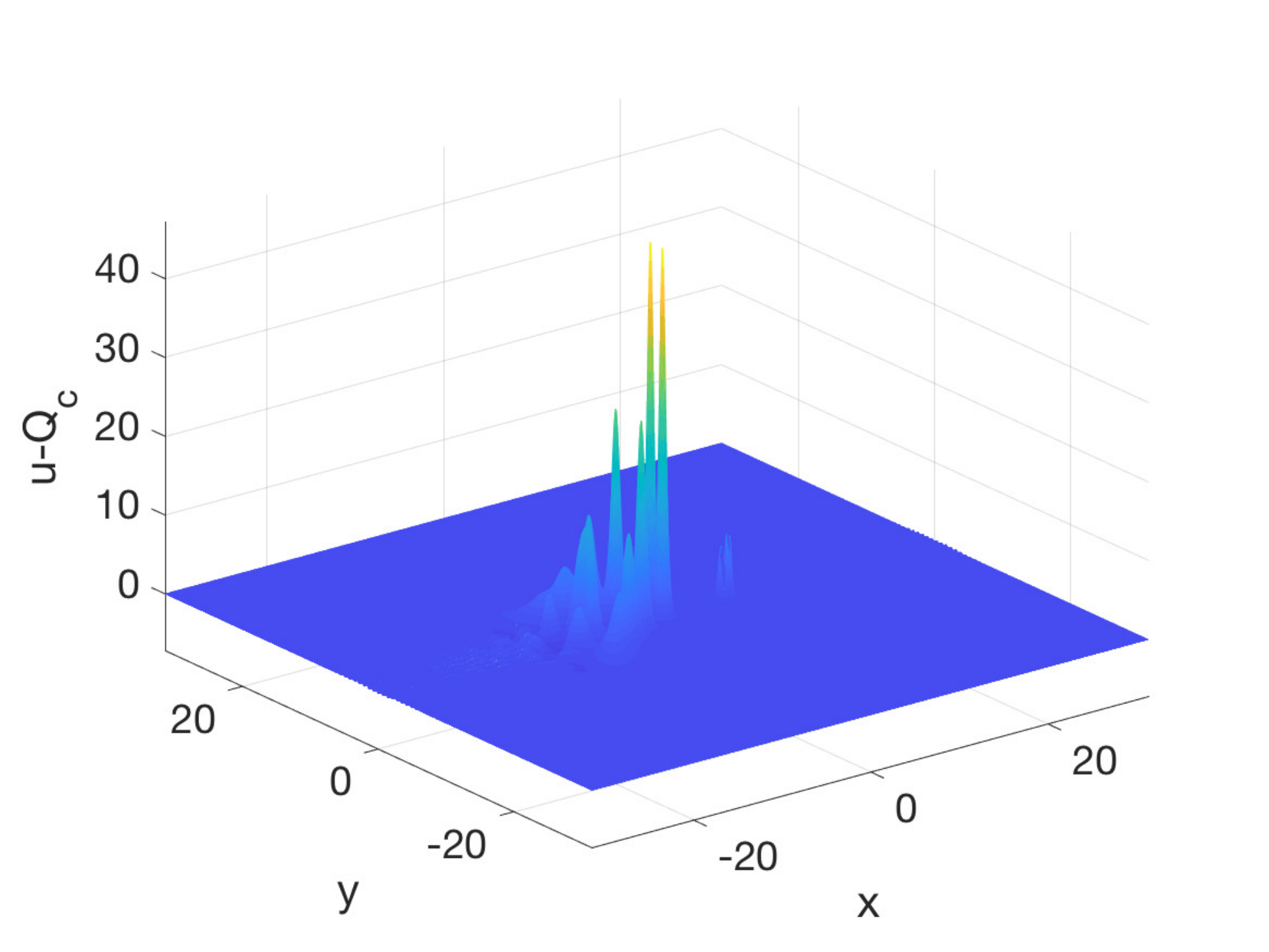}
\caption{Solution to \eqref{ZK} with $u(x,y,0)=25 \, e^{-(x^{2}-0.05y^{2})}$ 
at $t=0.5$ on the left, and on the right the difference 
between this solution and a fitted (to the first bump) soliton. }
\label{ZKp25gaussflat}
\end{figure}

From the previous simulations it is not yet clear whether more than one soliton can appear in the long term evolution of such data. To address this question, we consider 
once more initial data with a broad maximum, for example, we take 
$u(x,y,0)=25 \, e^{-(x^{2}+0.05y^{2})}$. The solution at $t=0.5$ is 
plotted on the left of Fig.~\ref{ZKp25gaussflat}. It looks as if several solitons appear in this case. On the right subplot we demonstrate the difference with a soliton fitted to the first bump, which appears to be close to a soliton. One would have to run simulations for much longer times in order to decide how 
many solitons will appear in the asymptotic solution. 


Therefore, the soliton resolution conjecture seems to hold for the 
2D subcritical ZK equation: in the long time behavior of solutions with sufficiently regular and sufficiently localized initial data only solitons and radiation appear. 

\section{The $L^{2}$-critical case}\label{S:cubic}

Since it is known that the direct integration of \eqref{ZKresc} with 
Fourier methods is challenging, we instead integrate  \eqref{ZK} and 
trace certain norms of the solution. It is expected, see \cite{FHRY}, that a blow-up is 
observed as $x\to\infty$. Therefore, we keep the term $x_{m}$ in 
\eqref{ZKresc} and solve \eqref{ZKnum} as follows: we 
choose $x_{m}$ in such a way that the maximum of the solution is at 
$x=x_{0}$, $y=0$ for all times. The quantity $x_{0}$ is chosen so that the radiation, propagating in negative $x$-direction, will hit the computational boundary 
(because of the imposed periodicity) only at a time shortly before the blow-up time, thus, 
its influence on the blow-up is negligible. 

We study perturbations of the soliton as in the previous section and 
Gaussian initial data. The results of this section 
confirm the Conjecture \ref{C:2}, which in some sense resembles the 1D critical gKdV equation. 
In particular, the gKdV examples showed that the blow-up mechanism for the studied norms 
($L^2$ or $L^\infty$) can be well captured, the more challenging task is to understand the 
velocity of the blow-up profile.  

\subsection{Perturbations of the soliton}

\begin{figure}[!htb]
\includegraphics[width=0.49\hsize]{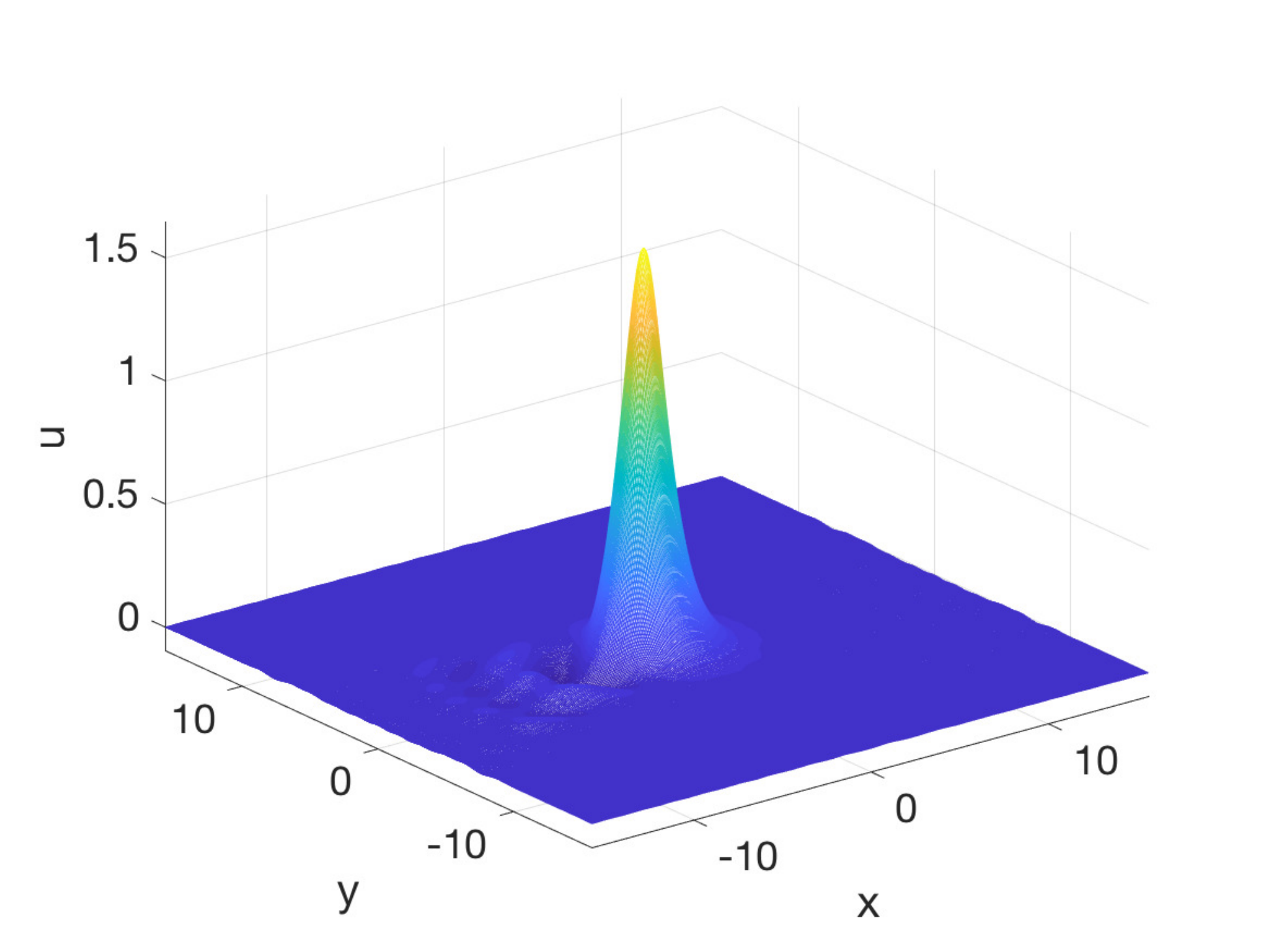}
\includegraphics[width=0.49\hsize]{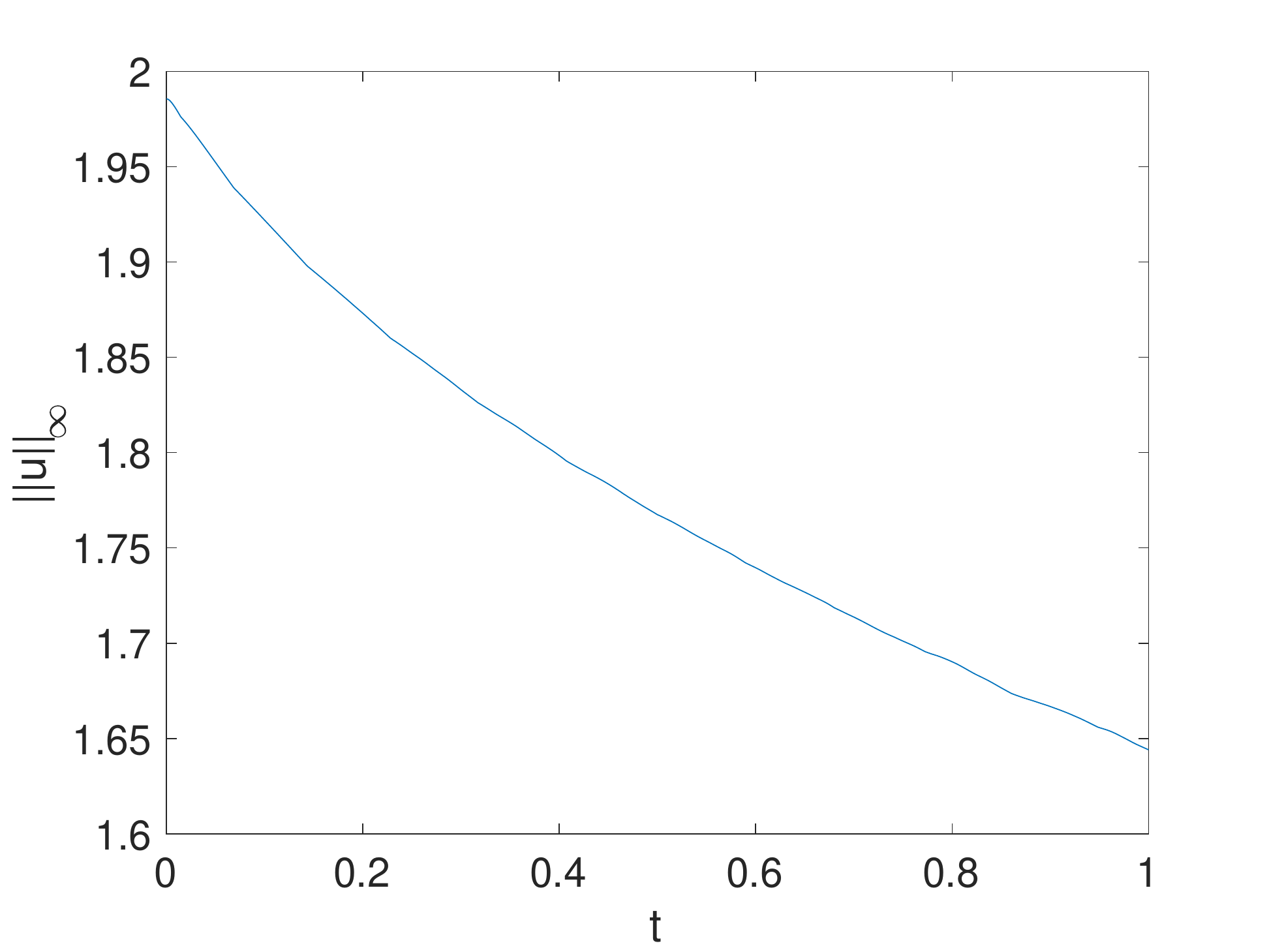}
\caption{Solution to \eqref{ZK} with $u(x,y,0)=0.9Q(x,y)$ at $t=1$ on the left 
and its $L^{\infty}$ norm depending on time on the right. }
\label{ZKp309sol}
\end{figure}

We start investigating evolution of ZK flow with initial data of 
perturbed solitons, of the form  $u(x,y,0)=\lambda \, Q(x,y)$, where 
$Q$ is numerically constructed soliton solution to \eqref{Q} with 
$c=1$. We start with $u(x,y,0)=0.9\, Q(x,y)$. A snapshot of the ZK 
evolution of this $u_0$ (at $t=1$) is plotted on the left subplot of 
Fig.~\ref{ZKp309sol}. The right subplot shows the $L^\infty$ norm 
depending on time, which appears to be monotonically decreasing, 
therefore, this solution disperses to infinity (of course, one could 
debate if the $L^\infty$ norm ever stabilizes at a certain value, as 
for example in Fig. \ref{ZKp2sol09t15}, in the present situation we 
see that the $L^\infty$ norm has a definite negative slope and that 
there is no increase after some time as in the stable cases. One 
could run this example for longer times, but since we approximate a 
situations in $\mathbb{R}^{2}$ by a situation on $\mathbb{T}^{2}$, 
the $L^{\infty}$ norm can never tend to zero, but will saturate at 
the level of the noise; this question should also be investigated analytically.)  In our simulations, because of the imposed periodicity, radiation (propagating towards negative values of $x$) reenters the computational domain on the right after some time (and hence, we have to stop our simulations at a certain time). We also note that dispersion propagates leftward in some wedge around the negative $x$-axis ($30^0$ as shown in \cite{CMPS}). Therefore, we conclude that the 
soliton is unstable against dispersion, as expected, for perturbations with a 
smaller mass than that of the soliton.

\begin{figure}[!htb]
 \includegraphics[width=0.49\hsize]{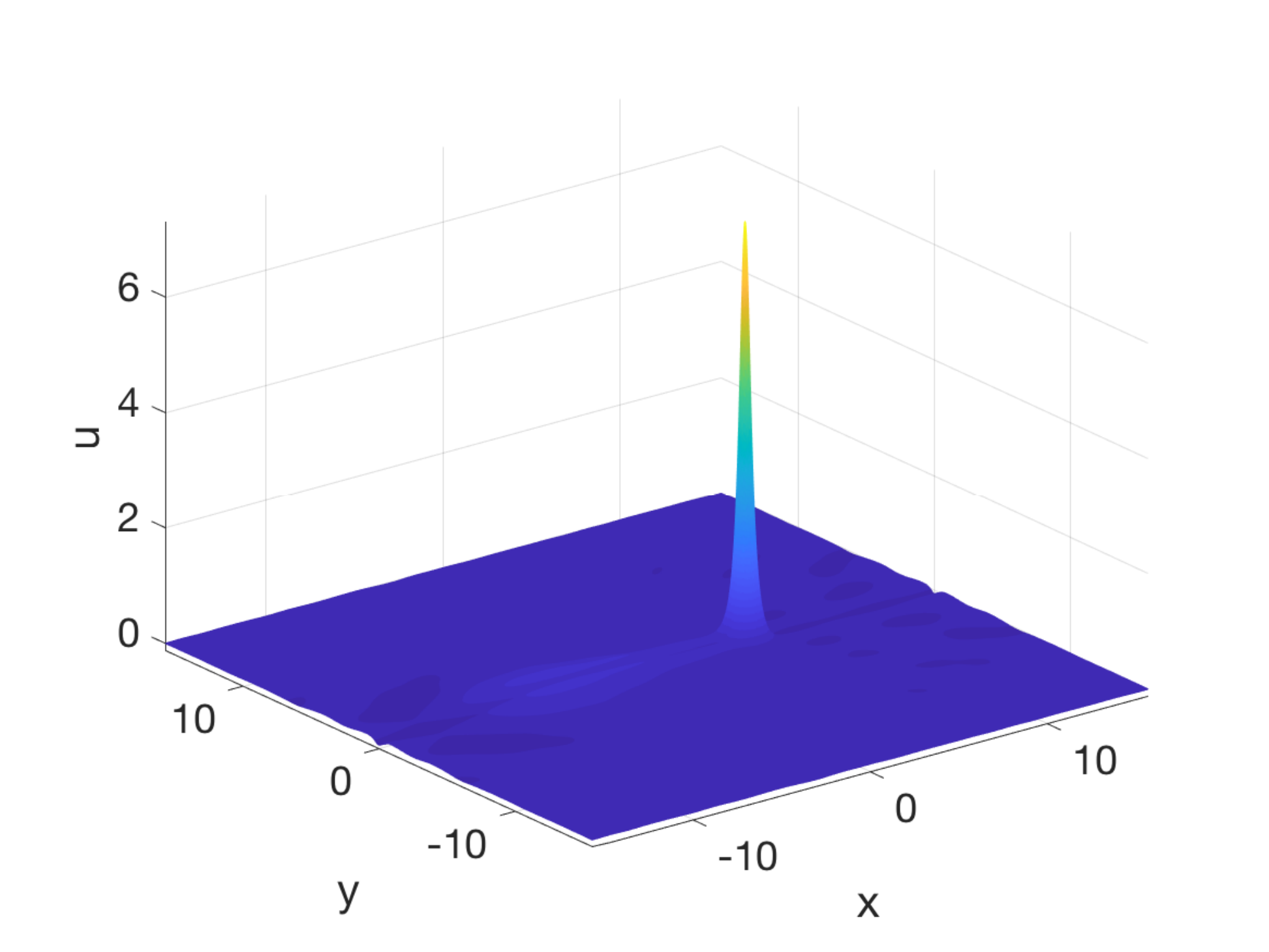}
 \includegraphics[width=0.49\hsize]{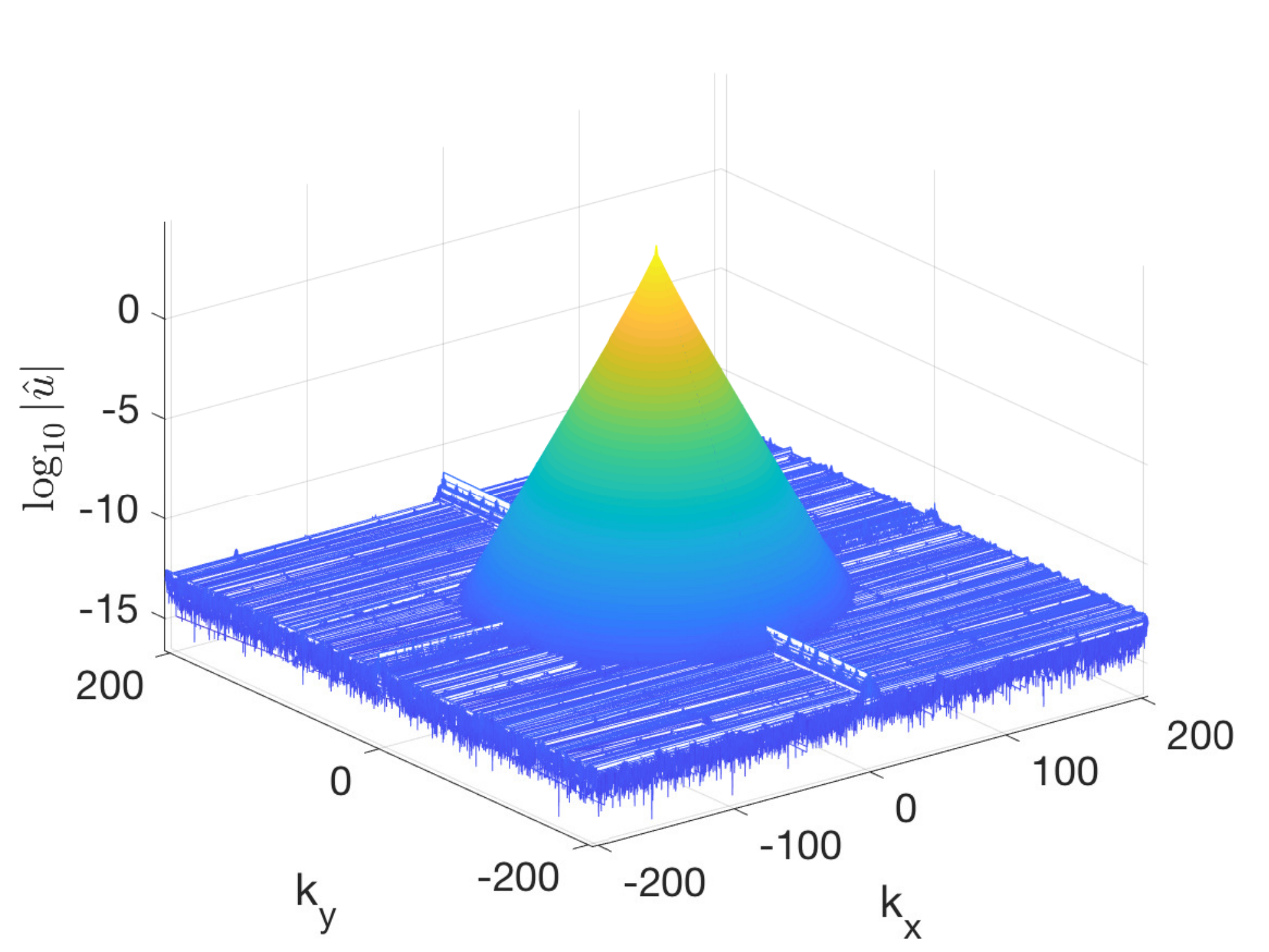}
\caption{Solution to \eqref{ZK} with $u(x,y,0)=1.1Q(x,y)$ at $t=2.6$ on the left, and the Fourier coefficients on the right. }
\label{Q11}
\end{figure}

Perturbations of solitons with larger mass, for instance, with the initial condition 
$ u(x,y,0) = 1.1 \, Q(x,y) $, lead to a blow-up solution with various diverging
norms. To approach this blow-up whilst maintaining at 
least plotting accuracy, we run the code first with \( 
N_{x}=N_{y}=2^{10} \) Fourier modes and \( N_{t}=2000 \) time steps 
for $t<2.6$. 
The snapshot of this ZK evolution at $t=3$ is shown in Fig.~\ref{Q11} on the left. 
Noting that the height of the bump is already at least 3 times larger than the initial height, implies that the soliton is unstable: a strong peak has formed (and moving with an increasing speed) as well as some bulk of radiation propagating in the negative wedge of the $x$-direction (recall that we are in a frame co-moving with the maximum, which is kept fixed at $(x_0,0)$). The Fourier coefficients of the solution on the right of Fig.~\ref{Q11} indicate that it is resolved to the order of the rounding error. 

\begin{figure}[!htb]
 \includegraphics[width=0.49\hsize]{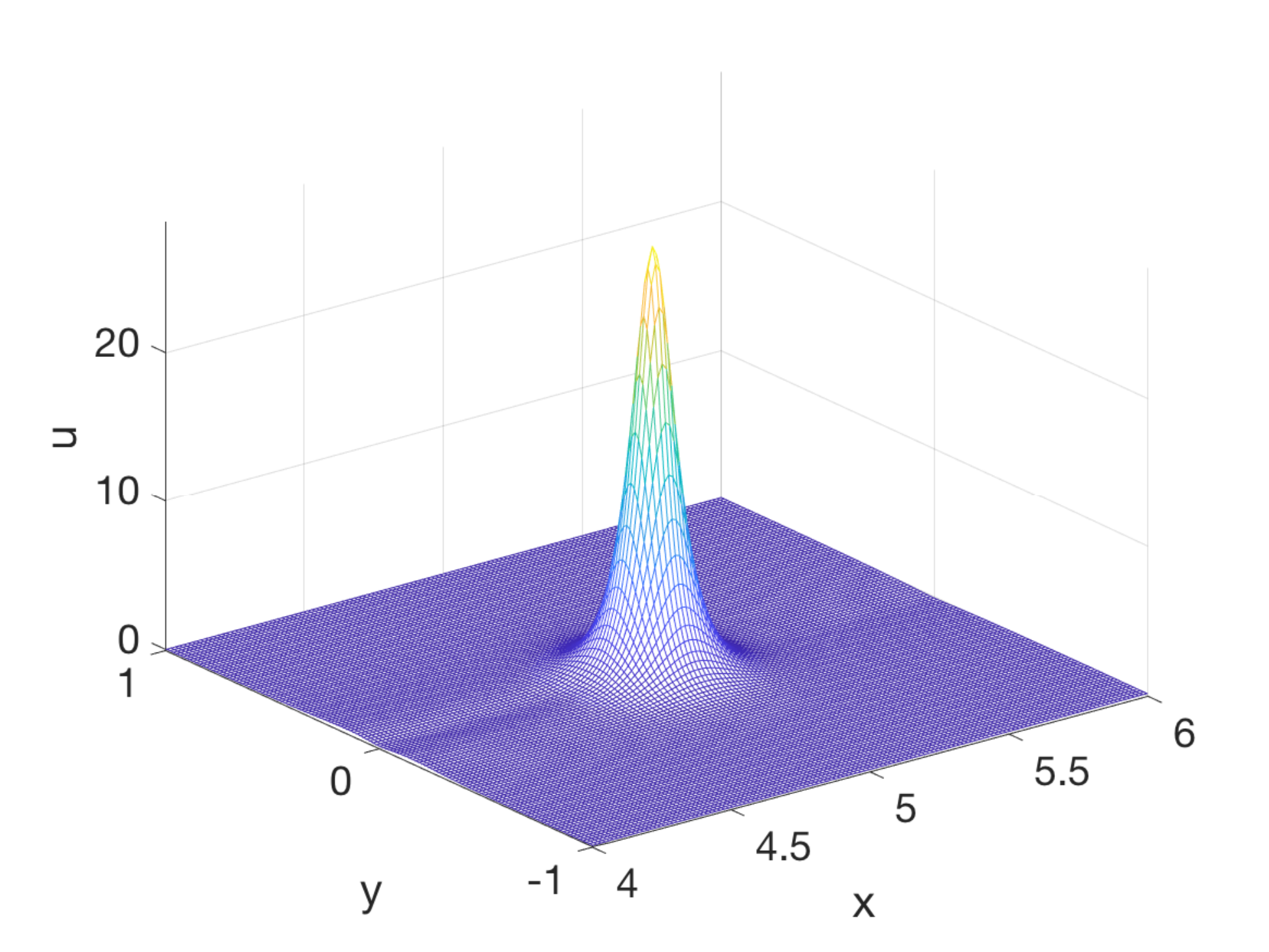}
 \includegraphics[width=0.49\hsize]{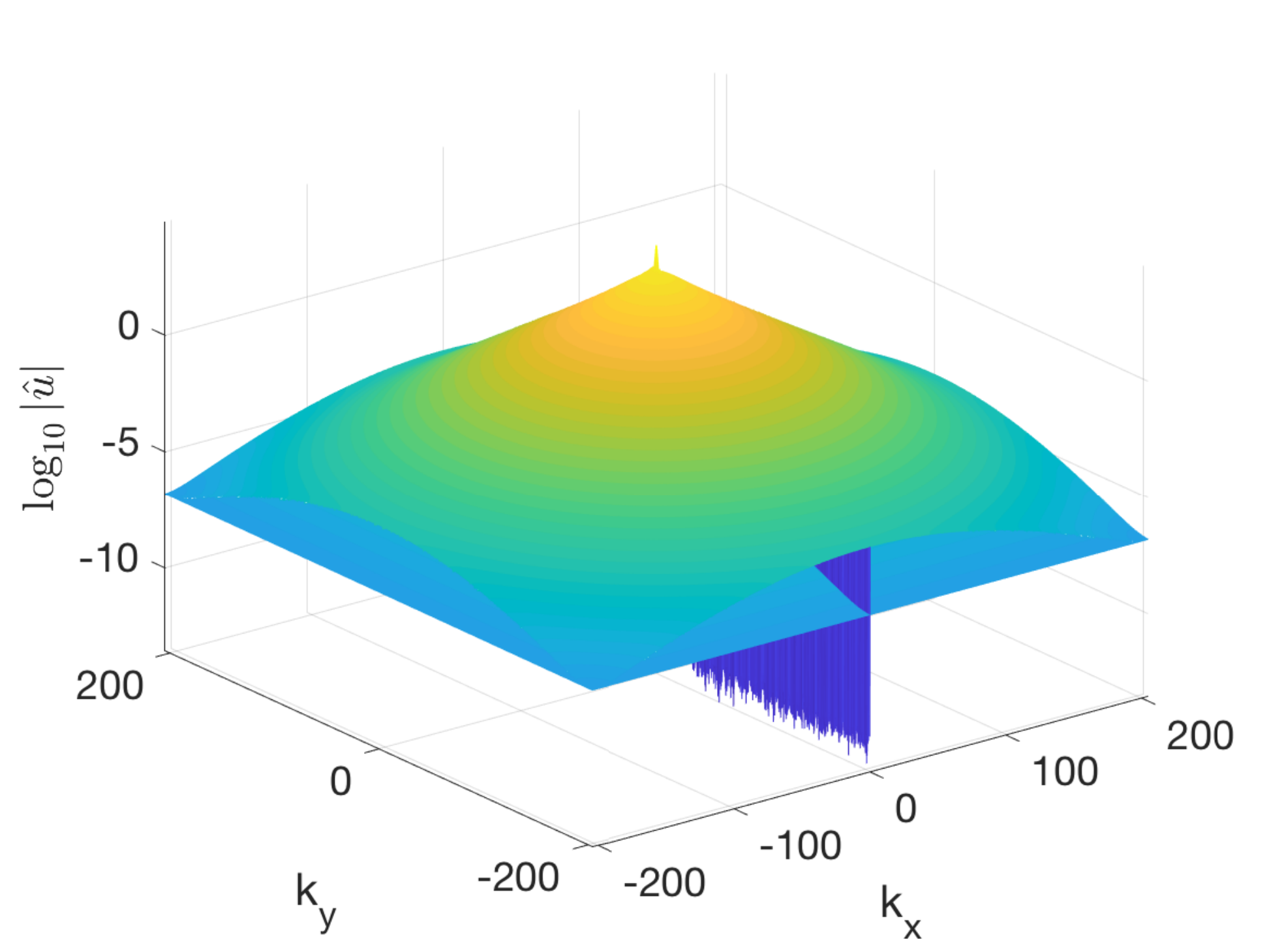}
\caption{Solution to \eqref{ZK} with $u(x,y,0)=1.1Q(x,y)$ at 
$t=3.1445$
on the left, and the 
corresponding Fourier coefficients on the right. }
\label{ZKp311sol}
\end{figure}

The solution shown in Fig~\ref{Q11} is then used as the initial condition for 
an ensuing computation with \( N_{x}=N_{y}=2^{11} \) and \( 
N_{t}=10^{4} \) time steps for \( t\leq 0.55 \). The code breaks at 
\( t\sim 0.541 \). In Fig.~\ref{ZKp311sol}, we show the solution at 
the last recorded time $t=0.5445$ on the left. The Fourier 
coefficients on the right of the same figure indicate that there is 
still spatial resolution beyond plotting accuracy at that time. This 
means that (similar to the case of blow-up in the Novikov-Veselov equation 
\cite{KK}) the resolution is first lost in time (compare this to the case 
of blow-up solutions in DS II system \cite{KS}, where the limiting factor is spatial resolution). 
The loss of resolution in time leads eventually to a breaking of the code.

\begin{figure}[!htb]
 \includegraphics[width=0.49\hsize]{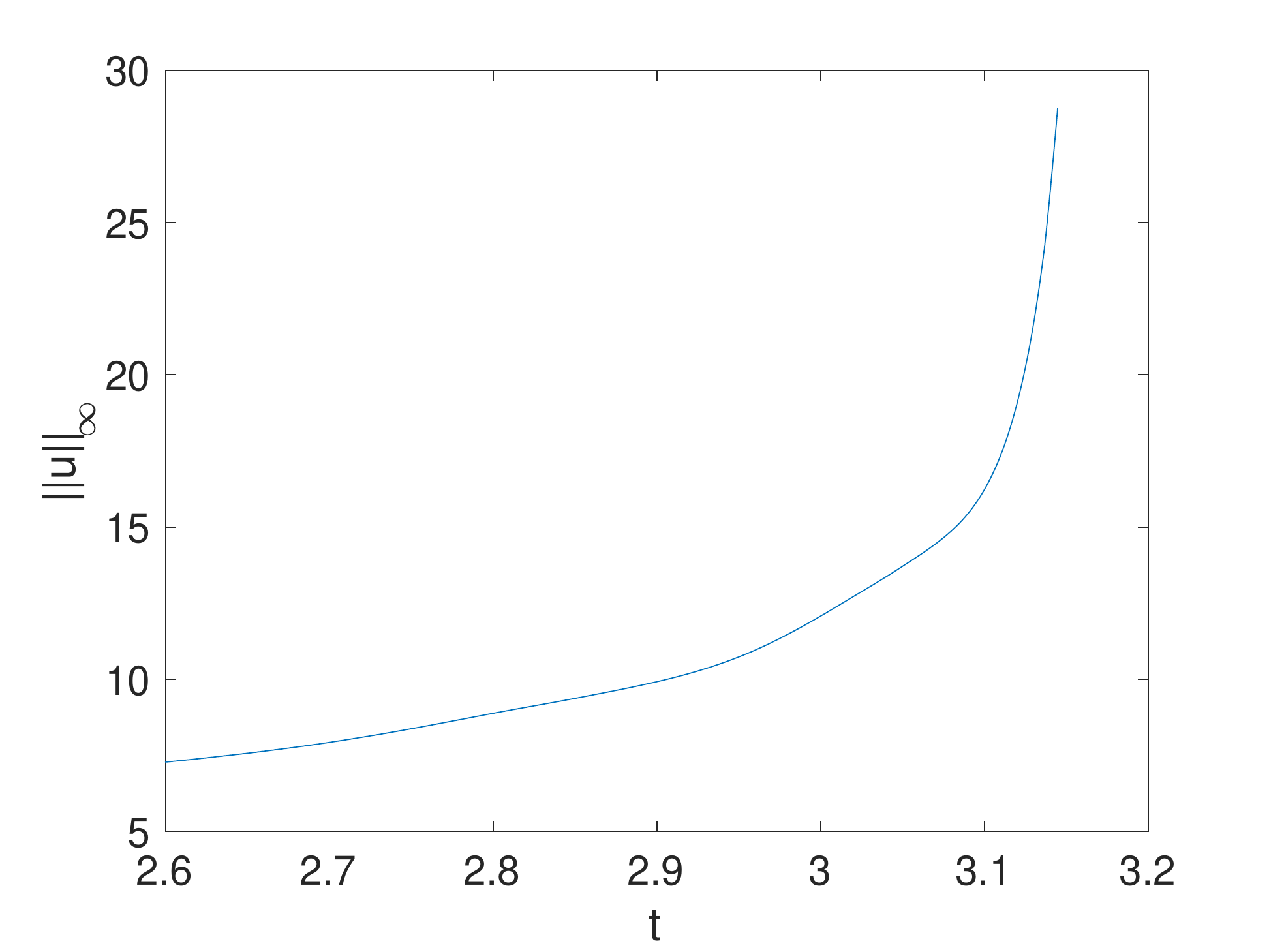}
 \includegraphics[width=0.49\hsize]{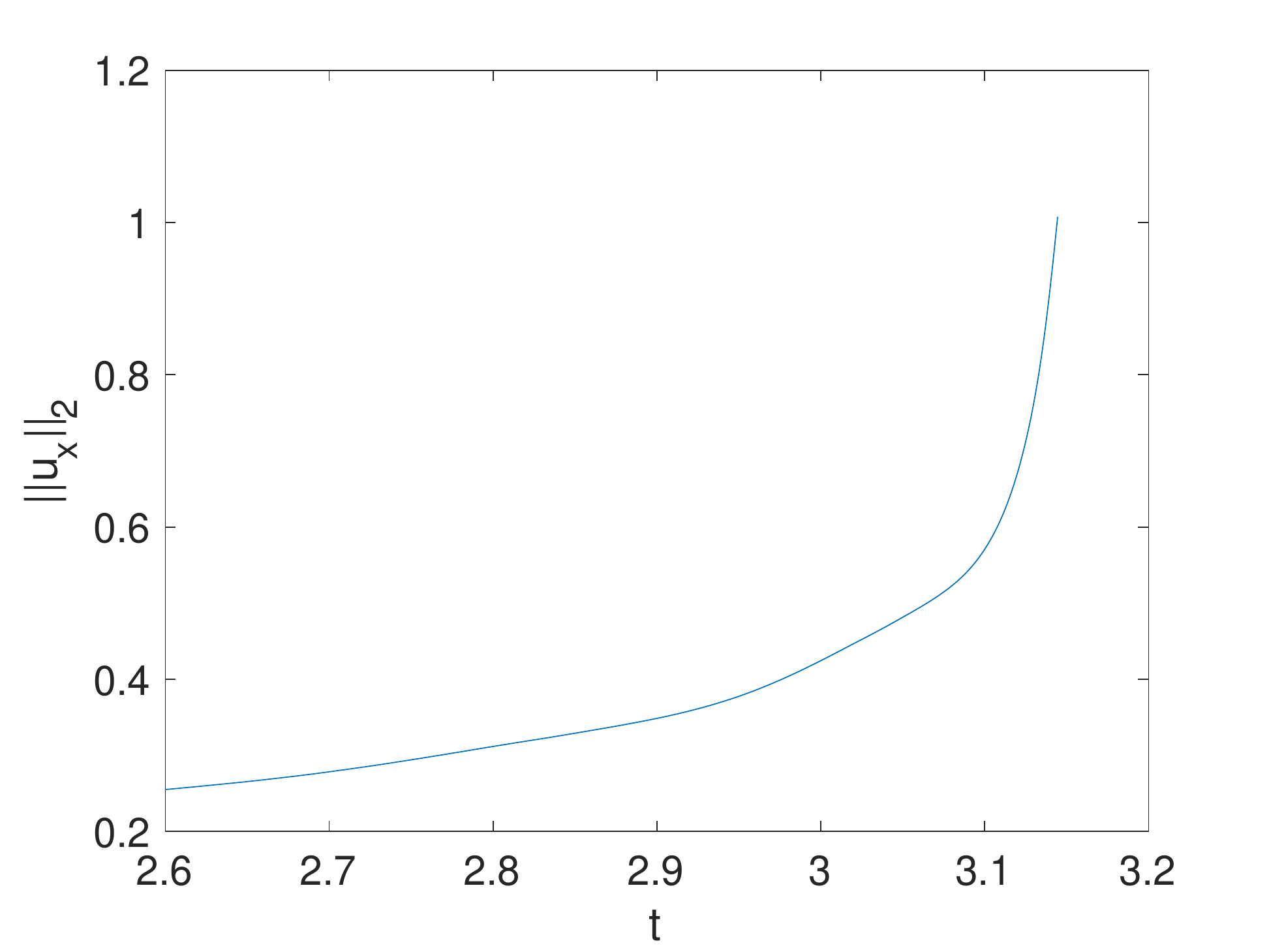}
\caption{Solution to \eqref{ZK} with $u(x,y,0)=1.1Q(x,y)$: 
on the left the $L^{\infty}$ norm of the solution; 
on the right the $L^{2}$ norm of $u_{x}$. }
\label{Q11norm}
\end{figure}

\begin{figure}[!htb]
 \includegraphics[width=0.49\hsize]{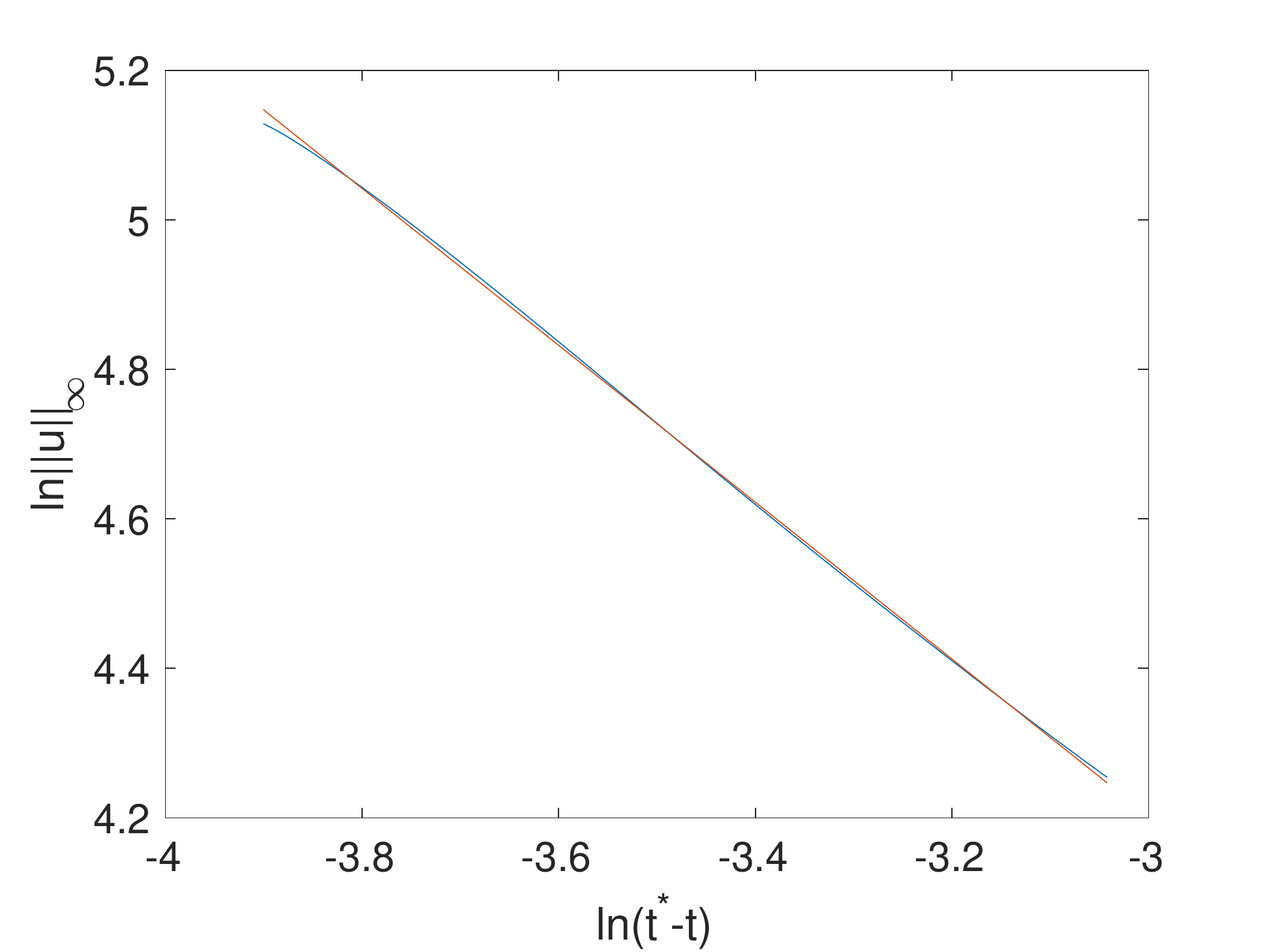}
 \includegraphics[width=0.49\hsize]{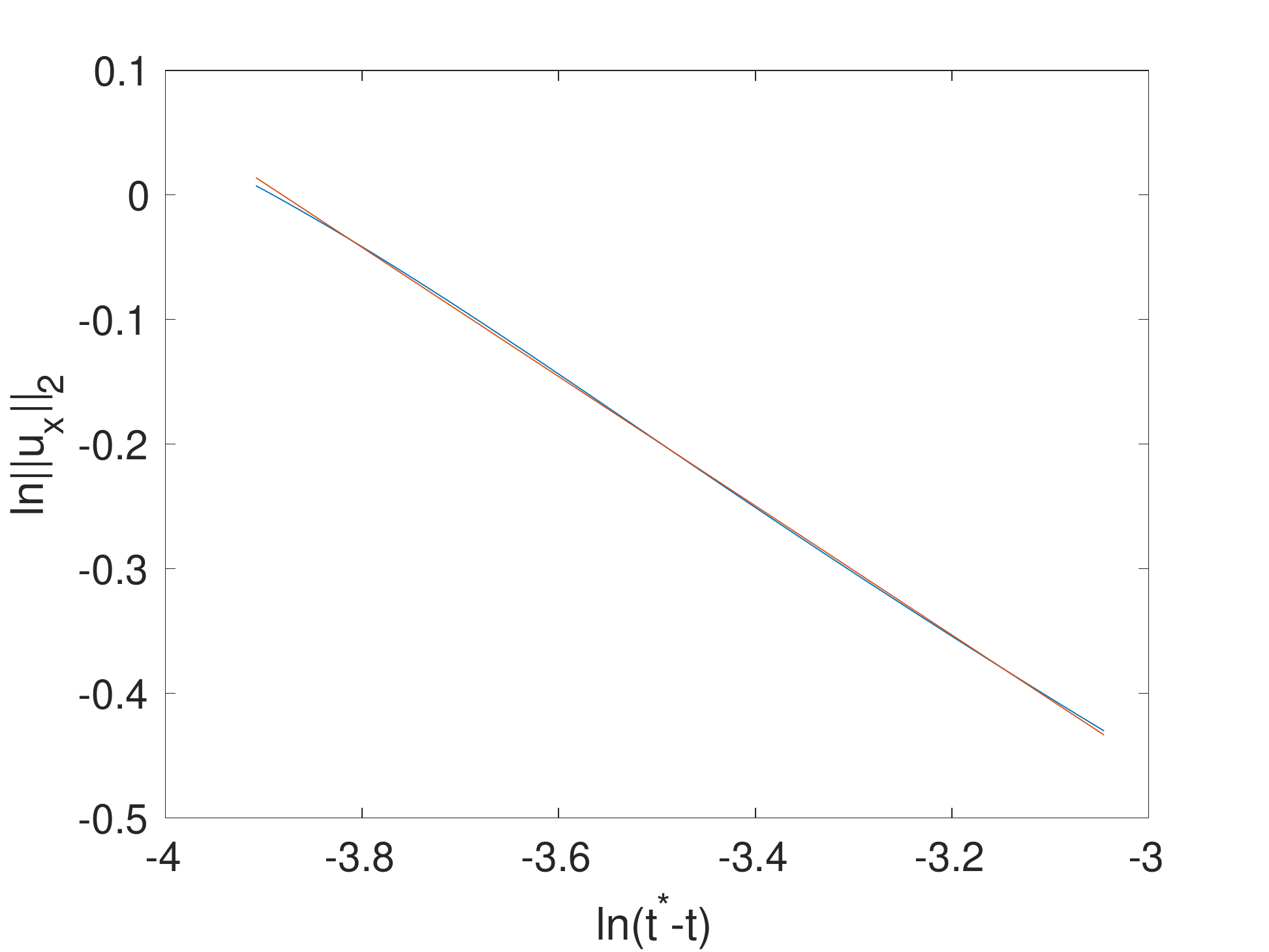}
\caption{Fitting of various norms of the solution to \eqref{ZK} with 
$u(x,y,0)=1.1Q(x,y)$ to $\ln g(t)\sim a\ln(t^{*}-t)+b$: 
on the left the $L^{\infty}$ norm of the solution fitted to 
$y= -0.4824 x +1.4266$; 
on the right the $L^{2}$ norm of $u_{x}$ fitted to 
$y= -0.5185 x -2.0124$; in red the fitted 
line. }
\label{Q11normfit}
\end{figure}

The divergence of the $L^{\infty}$ norm of the solution $u(t)$ and of the $L^{2}$ norm of 
$u_{x}(t)$ (shown in Fig.~\ref{Q11norm}) confirms the blow-up behavior. 
Furthermore, the growth of the norms in Fig.~\ref{Q11norm} provides 
information on the blow-up mechanism. Assuming that the blow-up core (the first bump) has a self-similar structure (a dynamic rescaling of the $Q$ profile), 
we fit various norms $g(t)$ close to the blow-up time to the following law
\begin{equation}\label{E:g}
\ln g(t)\sim a\ln(t^*-t)+b.
\end{equation}
The fitting is done for the last 500 recorded time steps (obtained results do not change 
significantly if slightly more or less points are used for the 
fitting) with the algorithm \cite{fminsearch} implemented in Matlab 
as \emph{fminsearch}. For the $L^{\infty}$ norm we find $a=-0.4824$,     
$b=1.4266$ and $t^{*}=0.5625$. For the $L^{2}$ norm of $u_{x}$ we get $a=-0.5185$,     
$b=-2.0124$ and $t^{*}=0.5646$. The quality of both fittings is shown in Fig.~\ref{Q11normfit}. 

The growth of the quantity or $v_x$ as in \eqref{ZKnum}, 
the speed of the frame co-moving with the maximum (see \eqref{ZKresinfty}),  is shown
on the left of Fig.~\ref{Q11norm}, which suggests that the 
blow-up takes place at infinity.  Fitting to $\ln v(t)\sim 
a\ln(t^{*}-t)+b$, we find $a=-1.0491$,     
$b=1.0553$ and $t^{*}=0.5647$. Note the agreement of the 
blow-up times, which shows the consistency of the used approach, though 
as mentioned in Remark \ref{buremark}, it is rather difficult 
to identify the blow-up rate of the velocity. 
\begin{figure}[!htb]
\includegraphics[width=0.49\hsize]{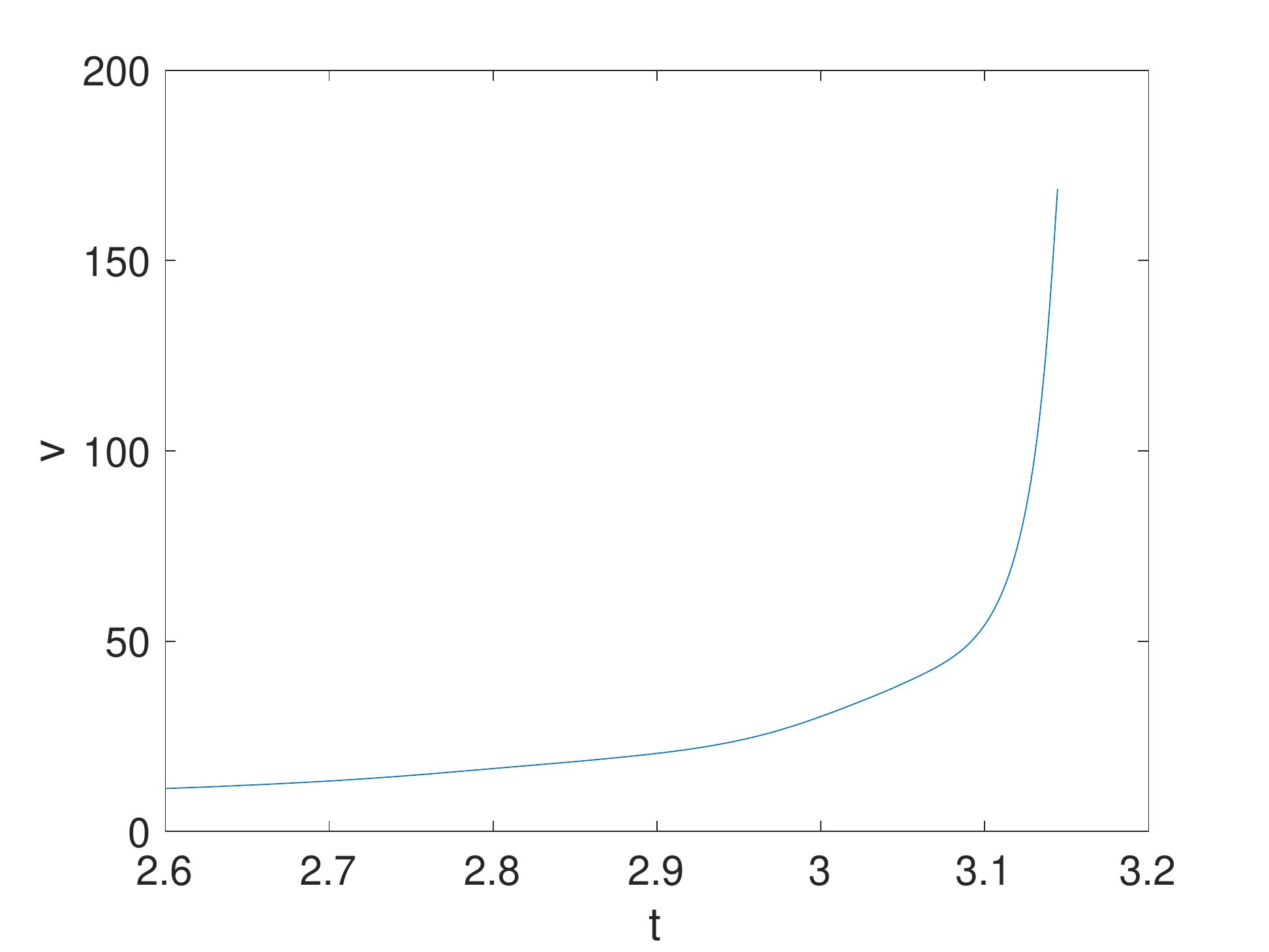}
\includegraphics[width=0.49\hsize]{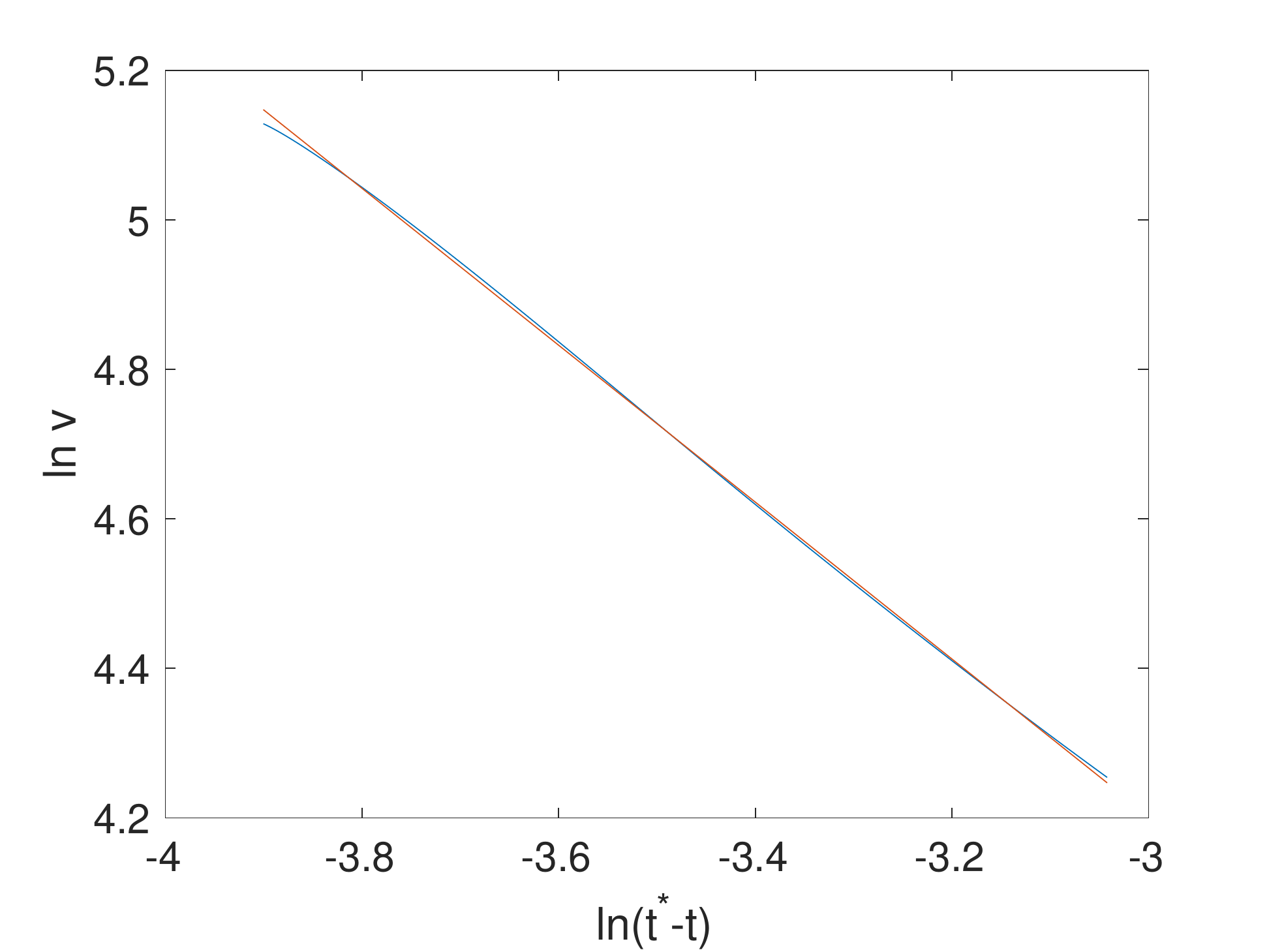}
\caption{Solution to \eqref{ZK} with $u(x,y,0)=1.1Q(x,y)$:  
on the left time dependence of $v_{x}$; 
on the right a fitting to $\ln v(t)\sim a\ln(t^{*}-t)+b$. }
\label{Q11v}
\end{figure}

\begin{figure}[!htb]
\includegraphics[width=0.7\hsize]{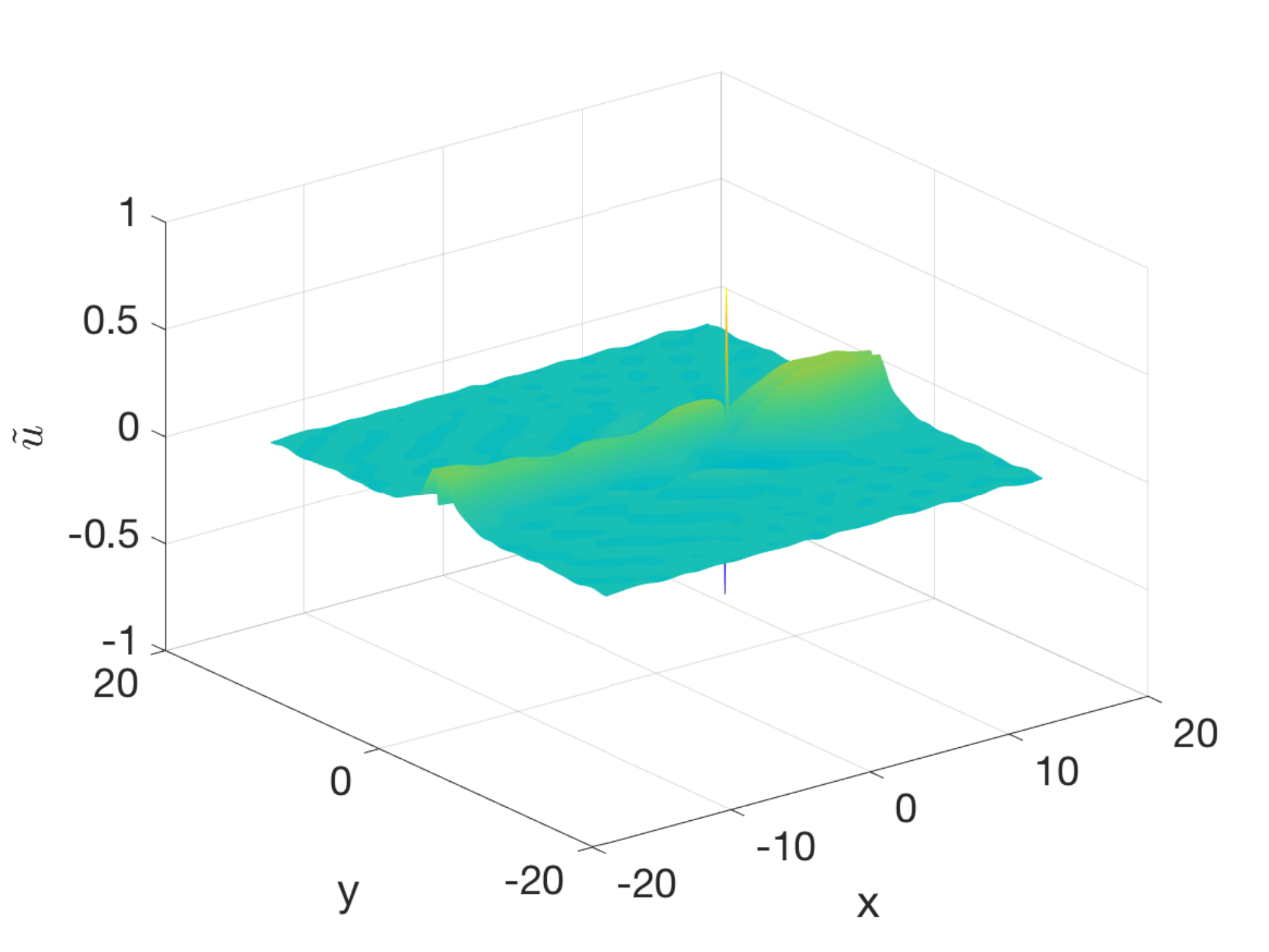}
\caption{The quantity $\tilde{u}$ in Conj.~\ref{C:2} giving the 
blow-up profile for the solution to \eqref{ZK} with  
$u(x,y,0)=1.1Q(x,y)$ at $t=3.1445$. }
\label{critprofile}
\end{figure}

To numerically determine the blow-up profile, we compute the quantity 
$\tilde{u}$, representing the limiting object in Conjecture ~\ref{C:2}. 
To do this,  we numerically determine the maximum and its location and determine via 
interpolation the dynamically rescaled soliton $Q$ according to \eqref{selfs}. 
The result is shown in Fig.~\ref{critprofile}. 
Remarkably, the difference between the expected  and computed blow-up 
profile is of the same order as the radiation profile. 
This indicates that the 
numerical estimate accurately represents
the actual blow-up.

\subsection{Gaussian initial data}
\begin{figure}[!htb]
 \includegraphics[width=0.49\hsize]{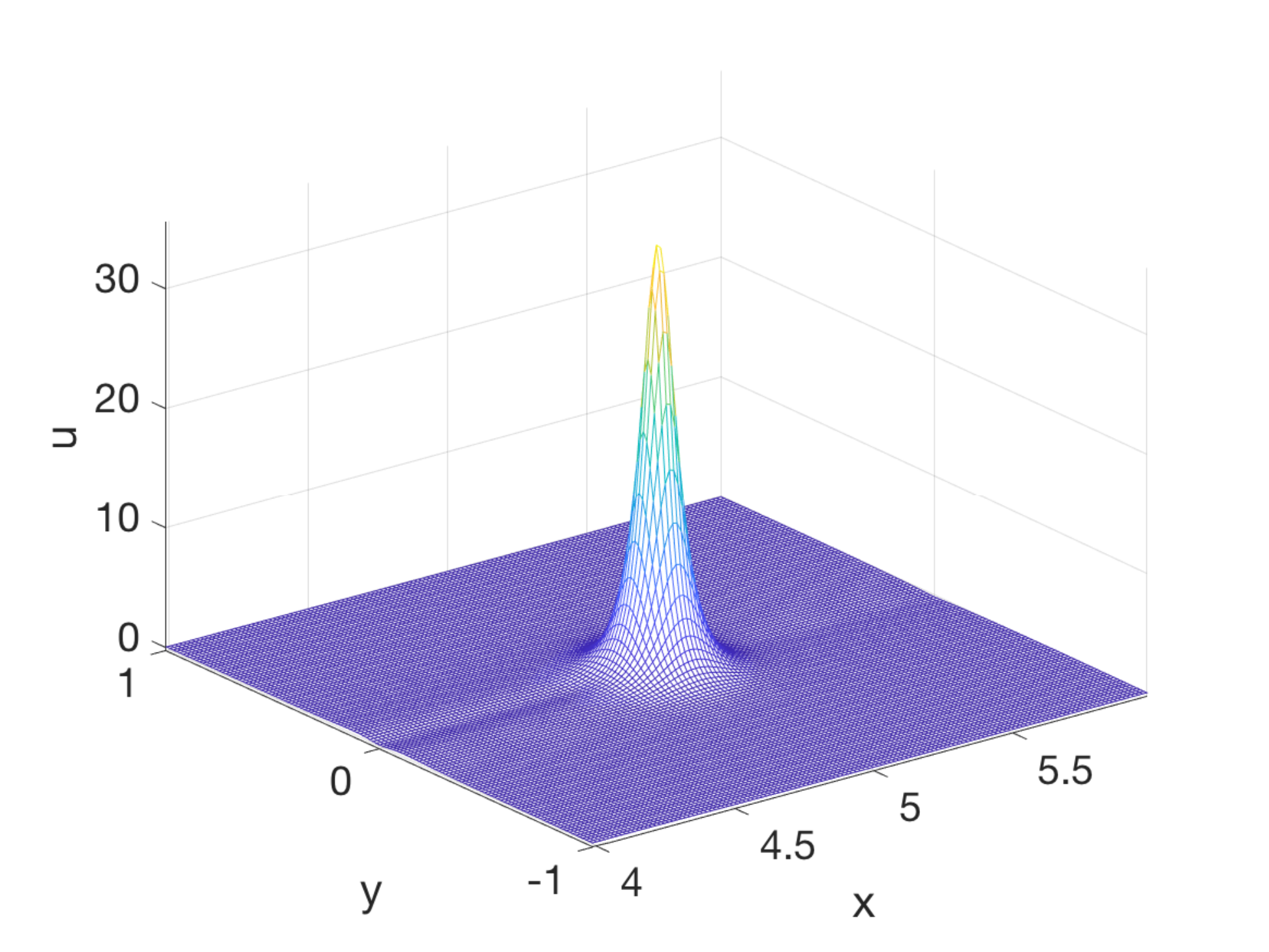}
 \includegraphics[width=0.49\hsize]{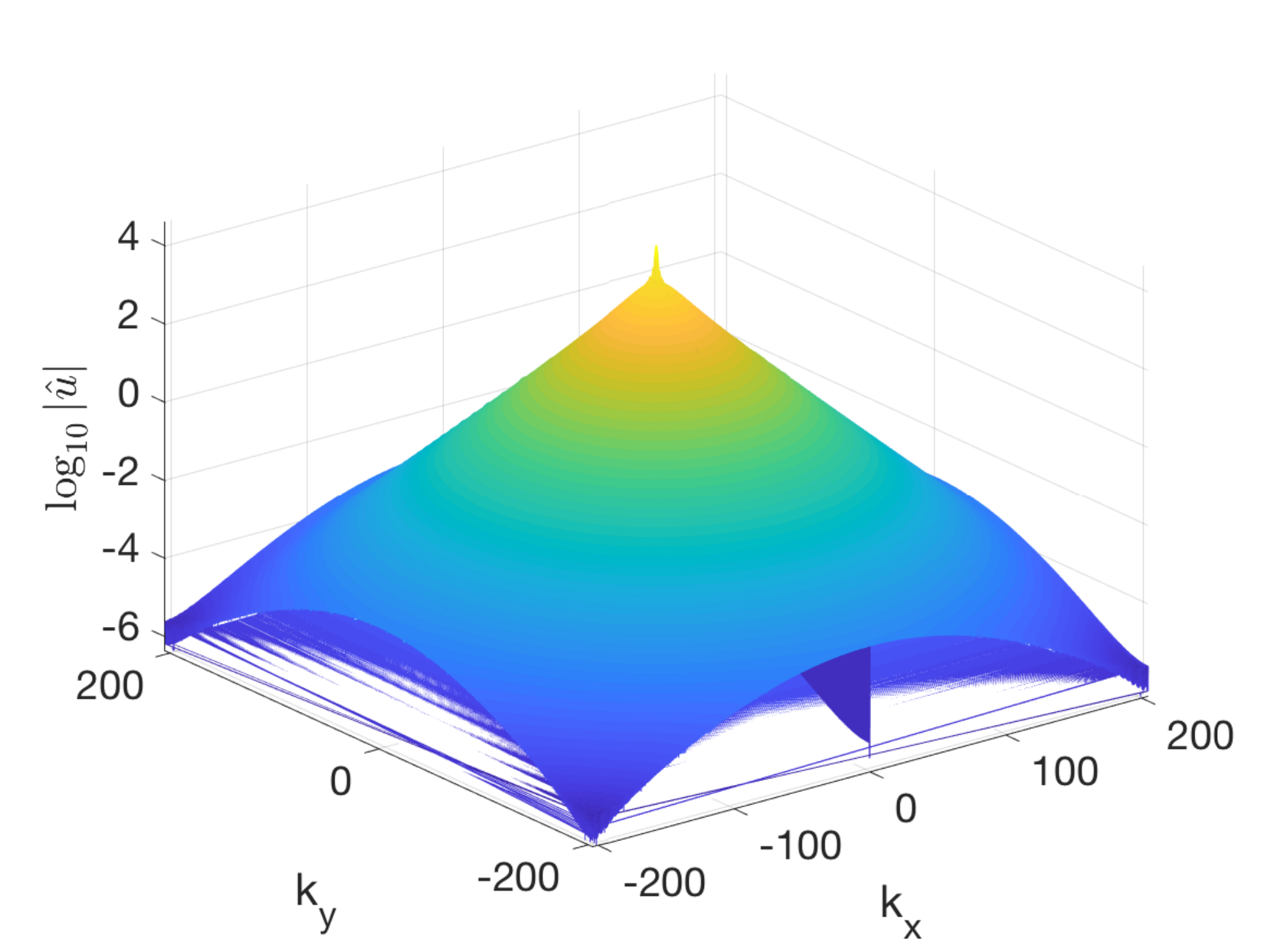}
\caption{Solution to \eqref{ZK} with $u(x,y,0)=3 \, e^{-(x^{2}+y^{2})}$ at $t=0.9175$ on the left, and the corresponding Fourier coefficients on the right. }
\label{ZKp33gauss}
\end{figure}

Formation of blow-up in finite time can be observed not only for 
perturbations of the soliton, but for more general initial data. We 
illustrate this on an example of Gaussian initial data, 
$u(x,y,0)=\lambda \, e^{-(x^{2}+y^{2})}$. For smaller $\lambda$, for 
instance $\lambda=2$, the evolution is dispersed. For 
$\lambda=3$, the evolution blows up in finite time. To study this, we use again 
$N_{x}=N_{y}=2^{10}$ Fourier modes on 
$5[-\pi,\pi]\times 5[-\pi,\pi]$ and $N_{t}=2000$ time steps for 
$t\leq0.5$. The resulting solution is then taken as the initial condition for 
an ensuing computation with $N_{x}=N_{y}=2^{11}$ Fourier modes and 
$N_{t}=10^{4}$ time steps for $t\leq 0.42$. The code breaks at 
$t=0.4175$. The solution at this specific time is shown in 
Fig.~\ref{ZKp33gauss} on the left. The Fourier coefficients for this 
time are given on the right of the same figure and show that the 
solution is still well resolved in the Fourier domain. Thus again, one 
runs out of resolution in time with the KdV-type equations. 

As far as the norms, the $L^{\infty}$ norm of $u$, the $L^{2}$ norm of 
$u_{x}$ and the velocity $v_{x}$ appear to blow-up. Fitting these 
norms as before to the law \eqref{E:g} 
for the last 500 recorded time steps gives $a=-0.53$, $b=1.44$ and $t^{*}=0.9436$ for 
the $L^{\infty}$ norm of $u$, $a=-0.64$, $b=-0.73$ and $t^{*}=0.944$ for 
the $L^{2}$ norm of $u_{x}$, and $a=-1.4$, $b=0.49$ and $t^{*}=0.9446$ for 
the $v_{x}$. The fitting errors are slightly larger than in the case of the 
perturbed soliton studied above (on the order of a few percent for 
the norms, though around 20\% for the velocity, which as we mentioned before  
is difficult to trace). However, there is good agreement of the fitted 
blow-up times in all cases. 

\begin{figure}[!htb]
\includegraphics[width=0.7\hsize]{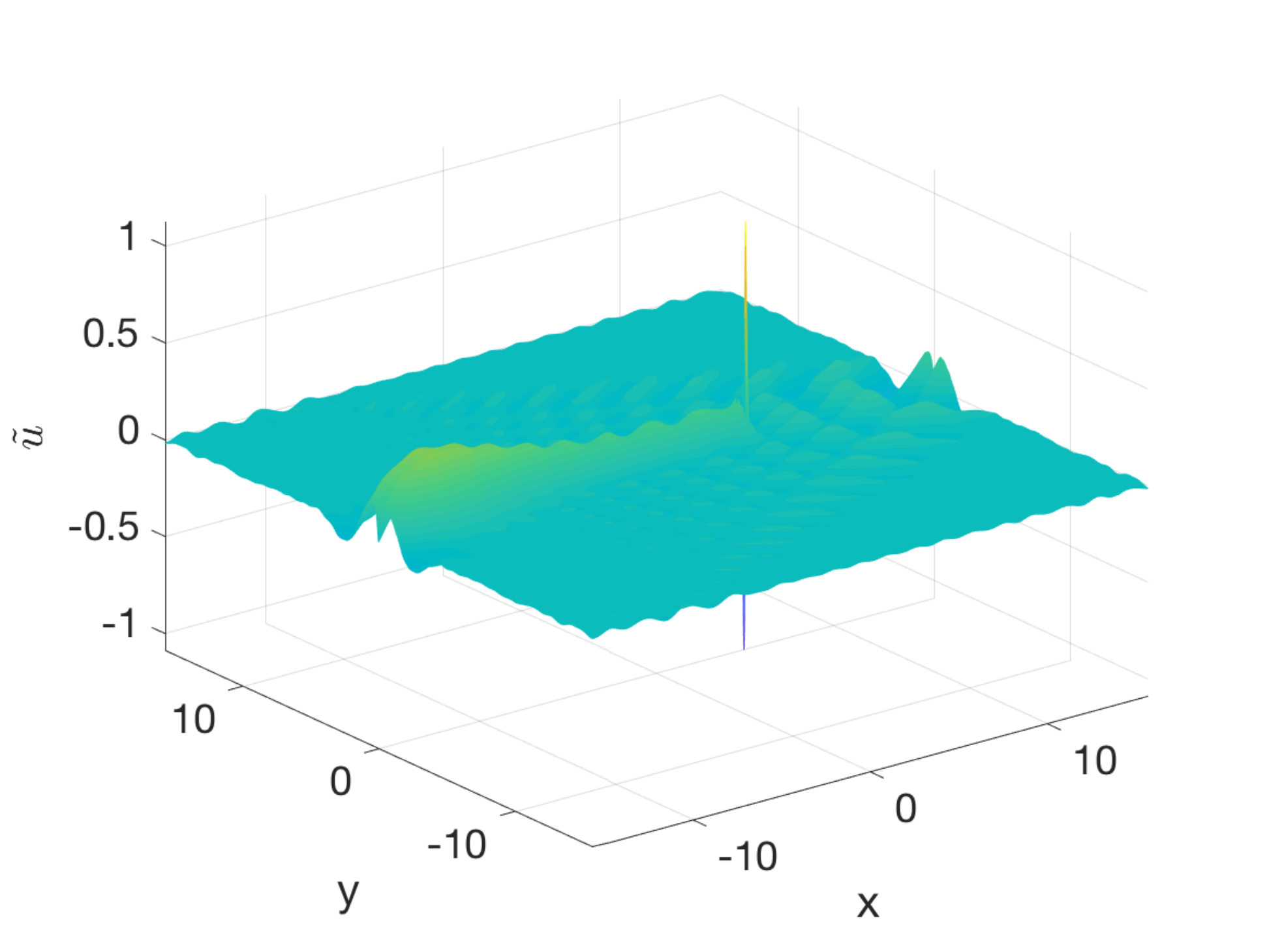}
\caption{The quantity $\tilde{u}$ in Conj.~\ref{C:2} giving the 
blow-up profile for the initial condition $u(x,y,0)=3 \, e^{-(x^{2}+y^{2})}$ at $t=0.9175$. }
\label{ZKp33gaussprofile}
\end{figure}

The blow-up profile appears to be again a dynamically rescaled 
soliton $Q$ as in Conjecture \ref{C:2}, see 
Fig.~\ref{ZKp33gaussprofile} for the residual. The residual (in the center) is 
slightly above the radiation background shows that the final phase of 
the blow-up is close, though not yet fully reached. 

\section{The $L^{2}$-supercritical case}\label{S:last}
In this section we study the $L^2$-supercritical case. 
As far as the rates of blow-up are concerned, it is 
numerically less challenging than in the critical case, since the blow-up 
happens on smaller time scales and centered at finite values of $x$. Though 
blow-up is always numerically challenging, at least 
some of the complications from the $L^2$-critical case are absent. 
As in the previous section we study perturbations of the soliton and 
Gaussian initial data. 
In this section we give positive confirmation to Conjecture \ref{C:3}. 


\subsection{Perturbations of the soliton}
As in the previous sections, we first study the stability of the 
soliton by considering initial data of the form $u(x,y,0)=\lambda 
Q(x,y)$, $\lambda>0$ in a co-moving frame with $c=1$. 

\begin{figure}[!htb]
 \includegraphics[width=0.49\hsize]{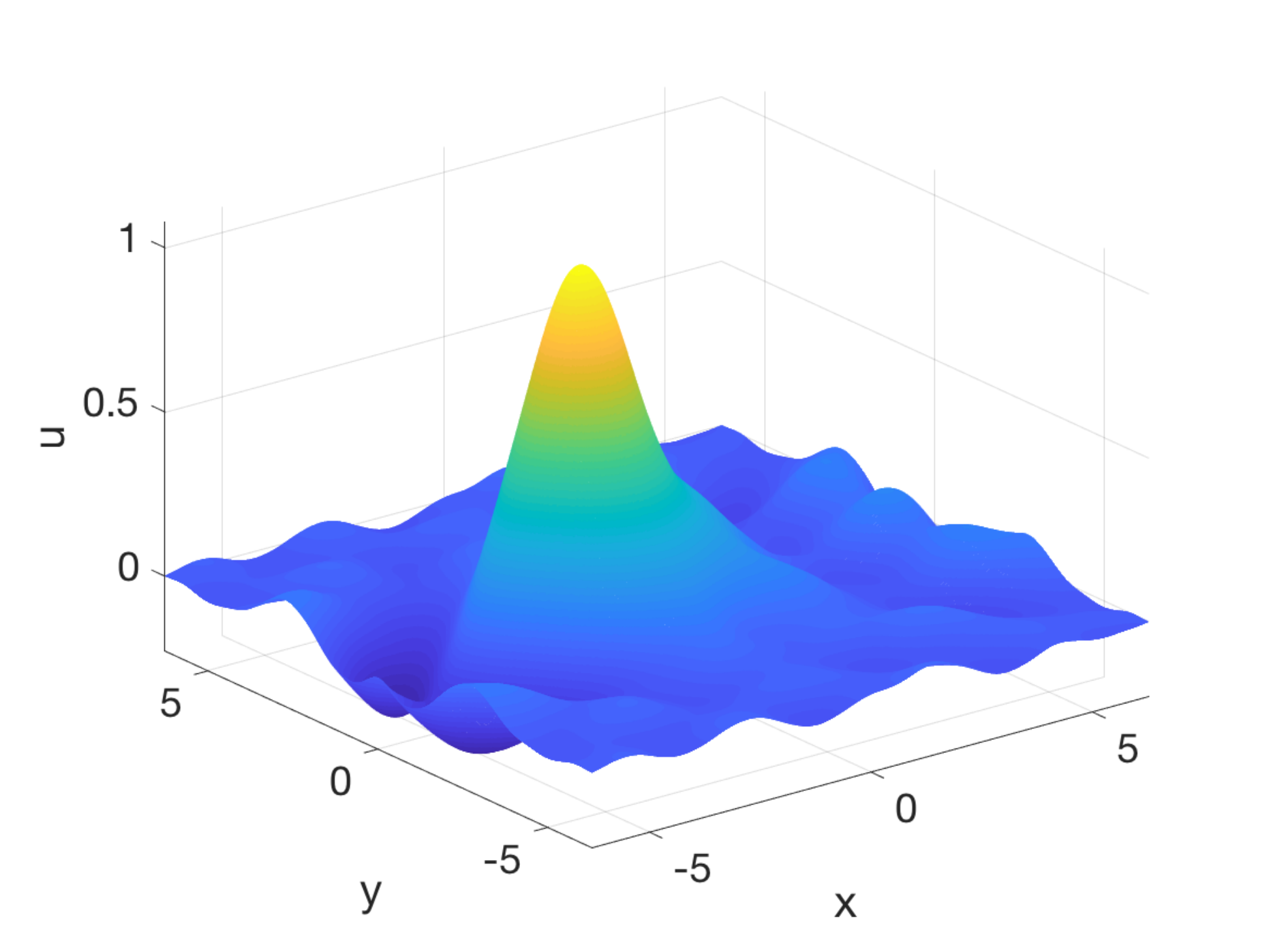}
 \includegraphics[width=0.49\hsize]{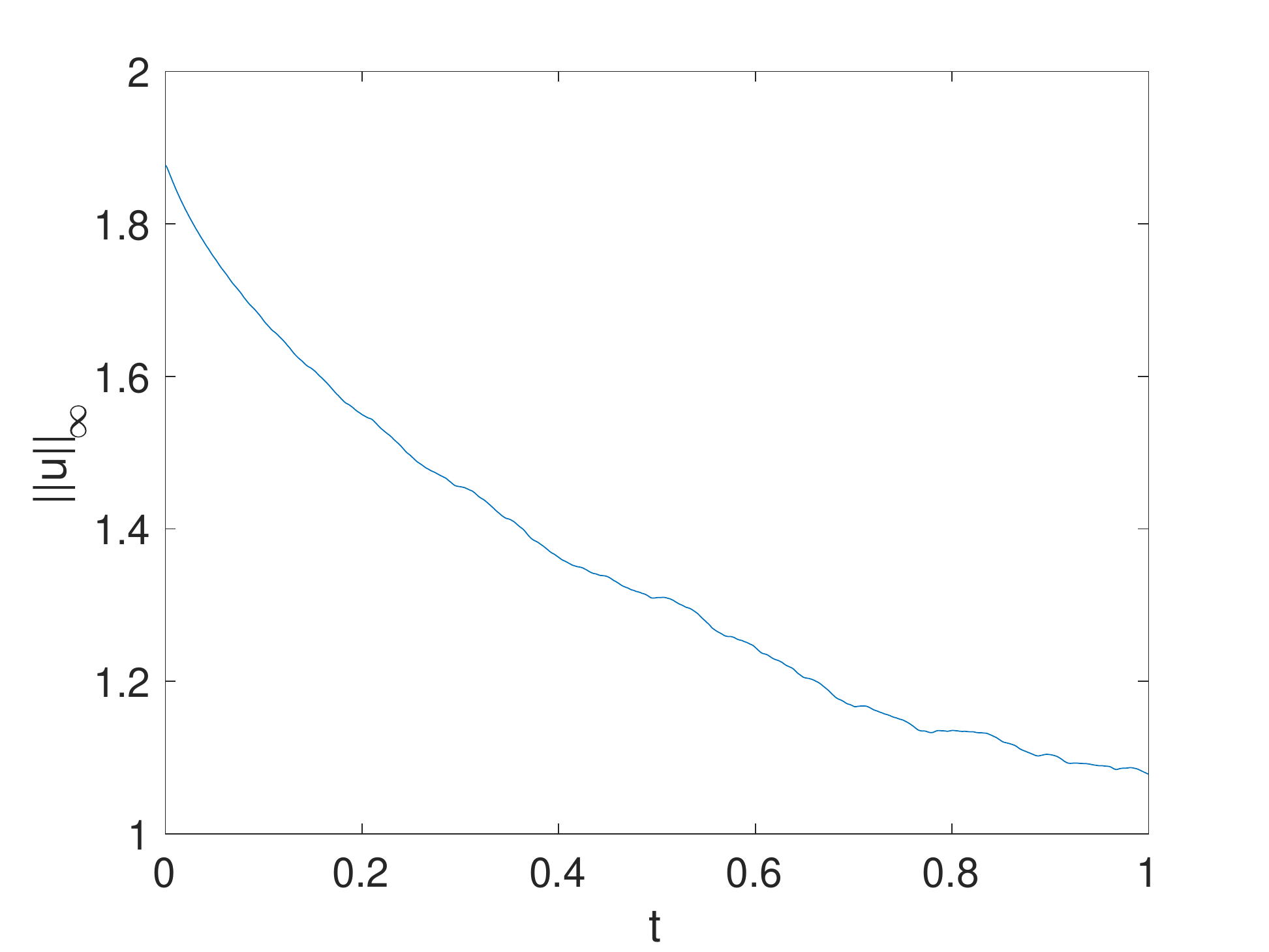}
\caption{Solution to \eqref{ZK}, $p=4$, with $u(x,y,0)=0.9 Q(x,y)$ 
at $t=1$ on the left, and time dependence of the corresponding $L^{\infty}$ norm
on the right. }
\label{ZKp4sol09}
\end{figure}

We start with the case $\lambda=0.9$ for $t=0$. 
The solution at $t=1$ is shown in Fig.~\ref{ZKp4sol09} on the left. The soliton is 
clearly unstable and disperses as time increases. Since we work on 
$\mathbb{T}^{2}$ here, the radiation cannot escape to infinity and 
forms a noisy background, into which the soliton will finally 
disappear. The time dependence of the $L^{\infty}$ norm of the solution is plotted on the right of Fig.~\ref{ZKp4sol09}. The norm is monotonically decreasing (the ripples being due to radiation reappearing on the other side of the computational domain and 
interacting then with the remaining peak). 

\begin{figure}[!htb]
 \includegraphics[width=0.49\hsize]{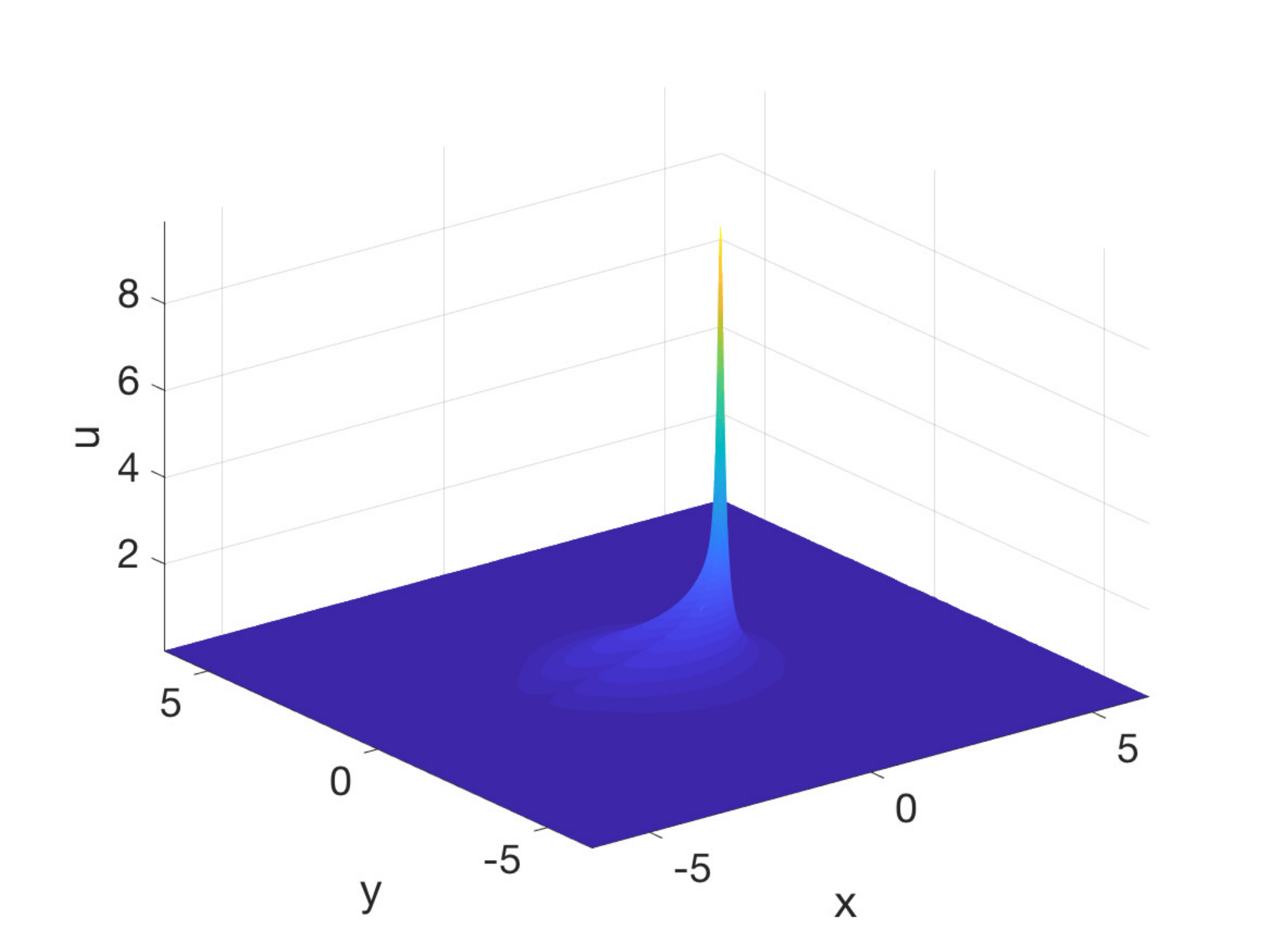}
 \includegraphics[width=0.49\hsize]{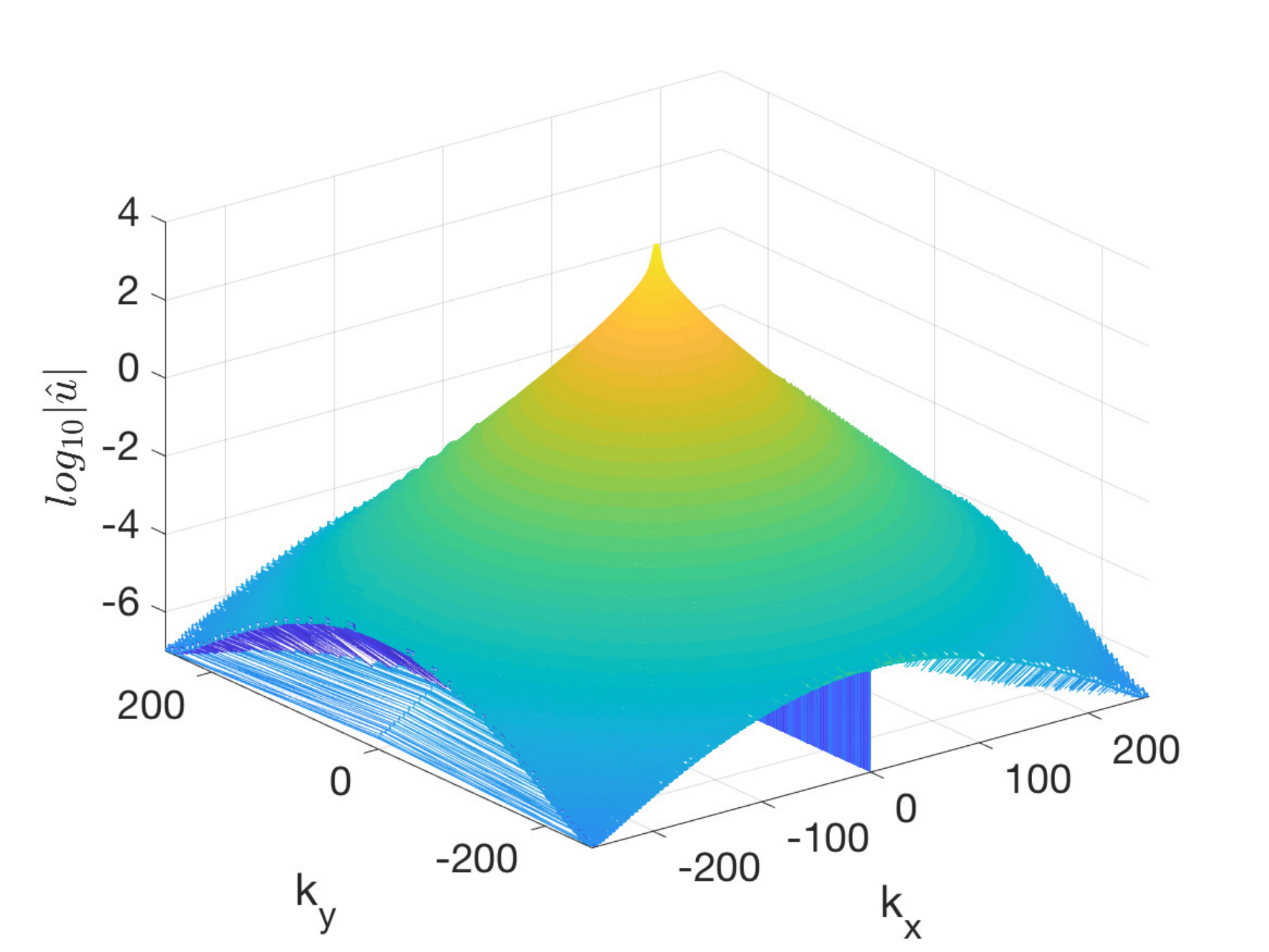}
\caption{Solution to \eqref{ZK}, $p=4$, with $u(x,y,0)=1.1 Q(x,y)$ at $t=0.1674$ on the left, and the corresponding Fourier coefficients on the right. }
\label{ZKp4sol11t101674}
\end{figure}

When $\lambda=1.1$, we use a smaller computational 
domain than in the $L^{2}$-critical case, since the blow-up will 
happen at finite values of $x$ and $y$. In practical terms this means 
that the blow-up profile will stay close to the bulk of the 
radiation. Thus, we can work on the domain $2[-\pi,\pi]\times 
2[-\pi,\pi]$. Consequently we operate with much higher spatial 
resolution than in the previous section, where we used a considerably 
larger domain. We use  $N_{x}=N_{y}=2^{10}$ Fourier modes and 
$N_{t}=20000$ time steps for $t\leq 0.17$. The code breaks at 
$t=0.1674$. The solution at the final recorded time can be seen on 
the left of Fig.~\ref{ZKp4sol11t101674}. The Fourier coefficients of 
the solution at the final time are shown on the right of the same 
figure. Note that the solution is still very well resolved spatially, 
though the code breaks a few time steps after the last recorded one. 
Thus, again the resolution in time is the limiting factor in blow-up 
computations for ZK. 

\begin{figure}[!htb]
 \includegraphics[width=0.49\hsize]{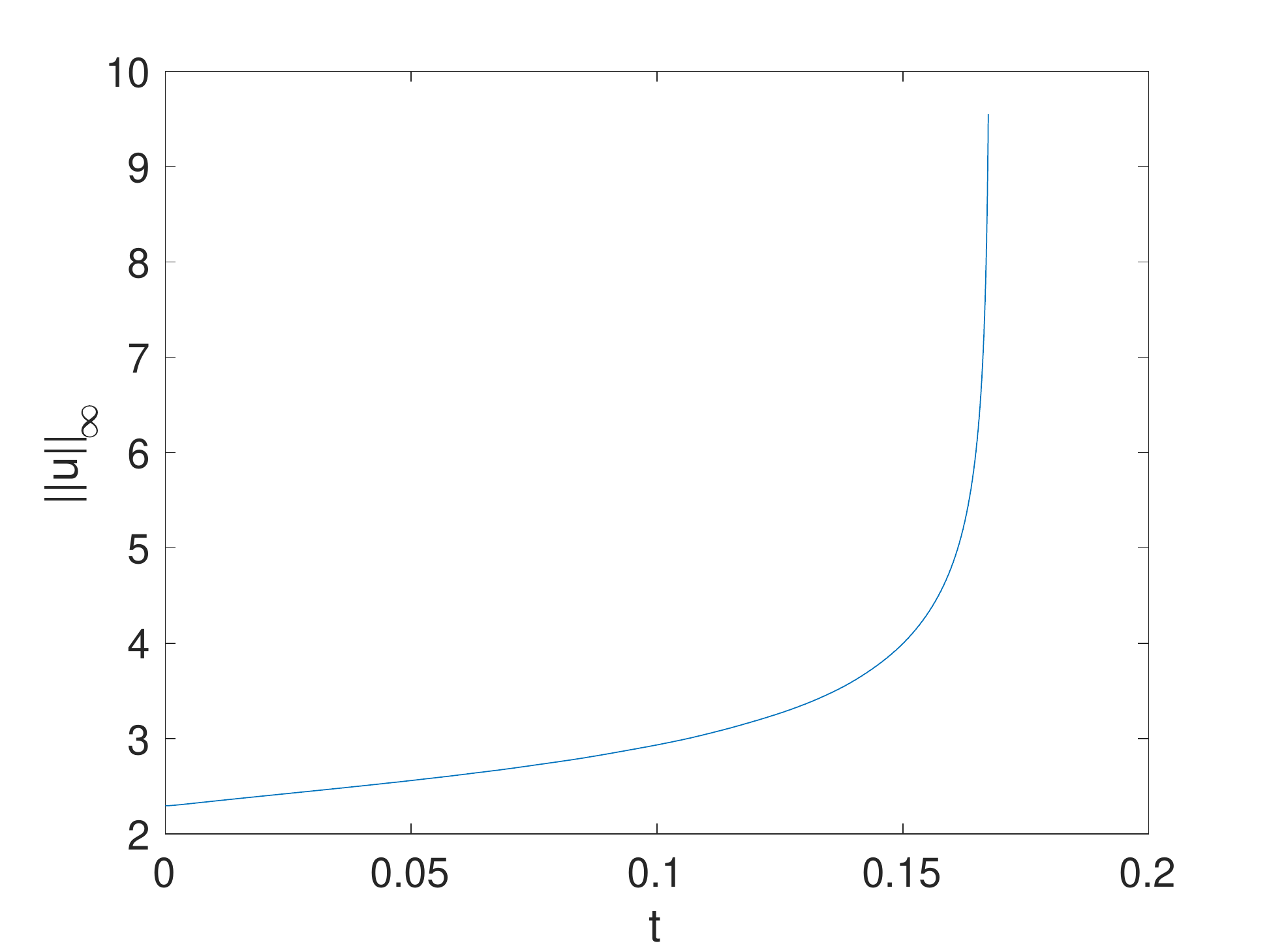}
 \includegraphics[width=0.49\hsize]{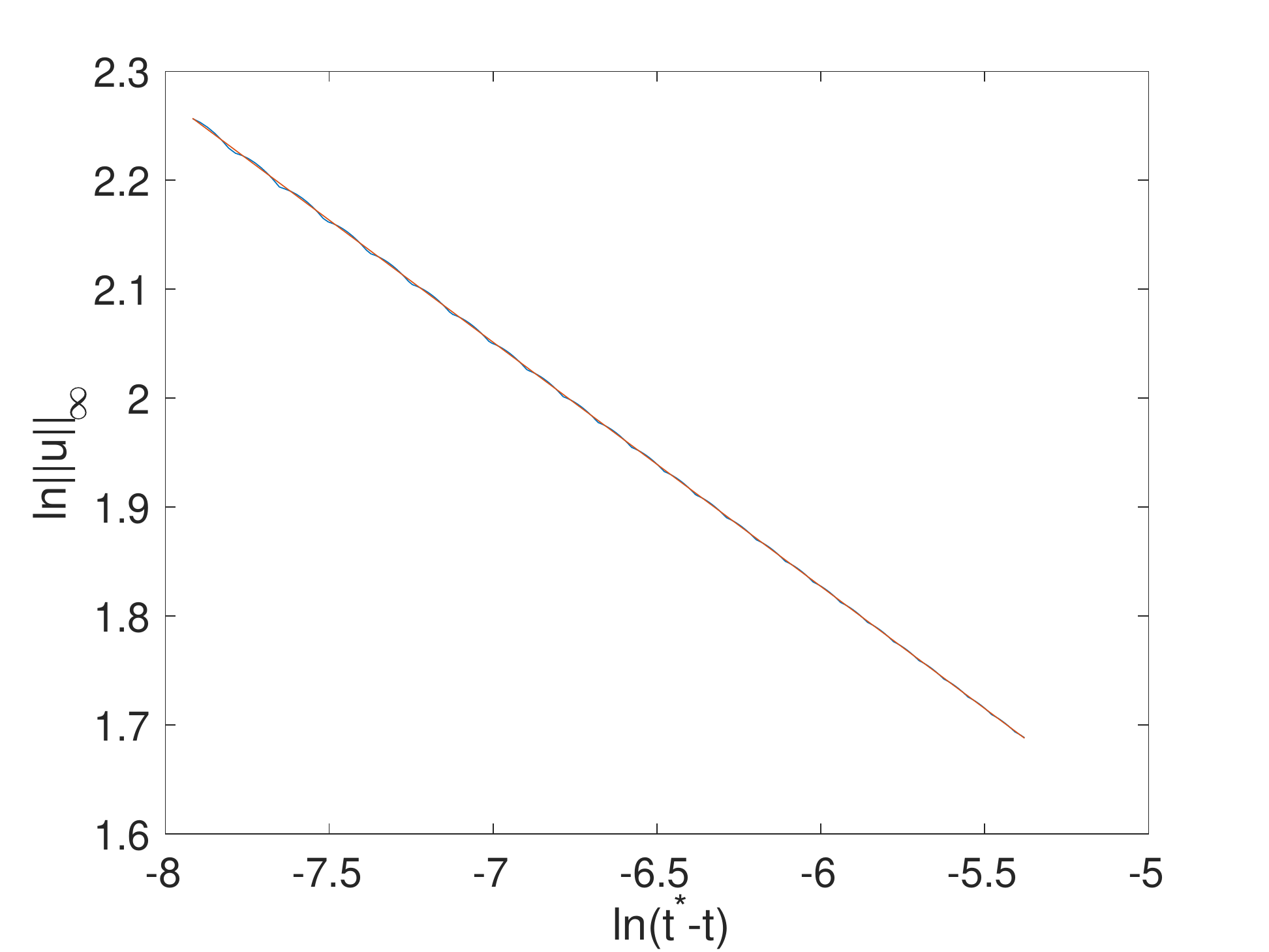}
\caption{The $L^{\infty}$ norm  of the solution to \eqref{ZK}, $p=4$, with  
$u(x,y,0)=1.1Q(x,y)$ on the left, and a fit  to  $\ln g(t)\sim 
a\ln(t^{*}-t)+b$ on the right  (in red the fitted 
line). }
\label{ZKp4solmax}
\end{figure}

The $L^{\infty}$ norm of the solution is shown in 
Fig.~\ref{ZKp4solmax}. It can be seen that the blow-up happens on 
much smaller time scales than in the critical case. A fit of the norm 
for the last 500 recorded time steps to \( \ln g(t)\sim 
a\ln(t^{*}-t)+b \)  yields $a=-0.22$    $b=0.48$ and     
$t^{*}=0.1678$. We observe that the result is virtually the same if we only fit 
the last 100 time steps. This is in accordance with the 
expectation that the power exponent for the blow-up rate $u_x(t)$ in terms of $t^*-t$ is $a=\frac29 = \frac12 \cdot \frac23 \approx 0.22$, see Conjecture \ref{C:3}, \eqref{ux2s}. 

The figures for the growth of the $L^{2}$ norm of $u_{x}$ are very 
similar to the ones for the $L^{\infty}$ norm of $u$, therefore, we do 
not present them here. The fitting of this norm for the 
last 500 time steps gives $a=-0.22$, $b=-0.38$ and $t^{*}=0.1678$ 
with a fitting error of the order of $10^{-4}$. The 
excellent agreement of the results by both norms confirms that 
the results for the blow-up in the $L^{2}$-supercritical case are more 
accurate and more stable than in the critical case. 

\subsection{Gaussian initial data}
Finally, we consider Gaussian initial data of the form $u(x,y,0)=\lambda \,
e^{-(x^{2}+y^{2})}$ with $\lambda>0$ in a stationary frame. First we observe that for 
smaller $\lambda$ (for instance $\lambda=1$) the initial bump is dispersed similar to the $L^2$-critical case of the perturbed soliton with the initial mass smaller than the soliton mass. 

\begin{figure}[!htb]
 \includegraphics[width=0.49\hsize]{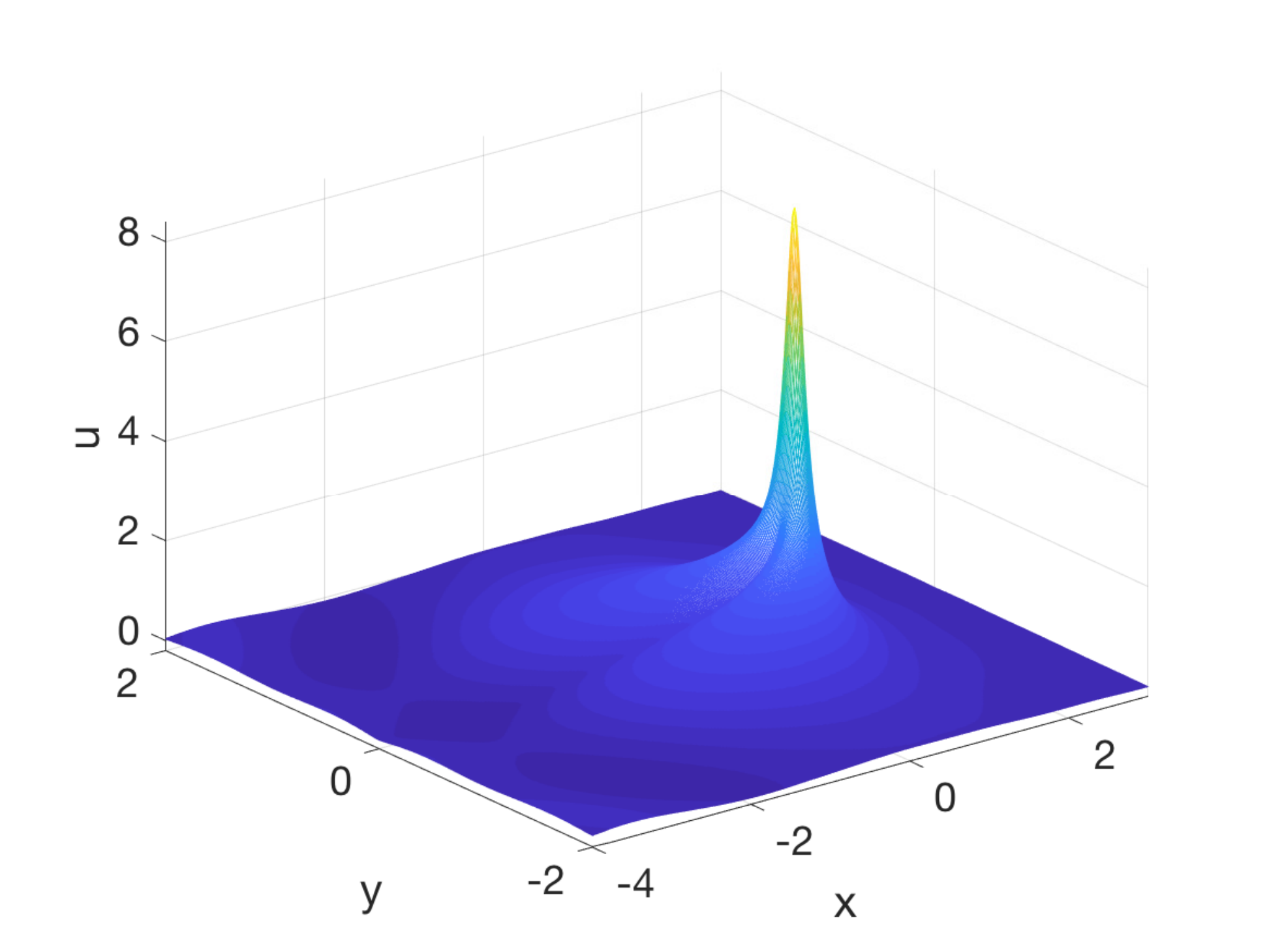}
 \includegraphics[width=0.49\hsize]{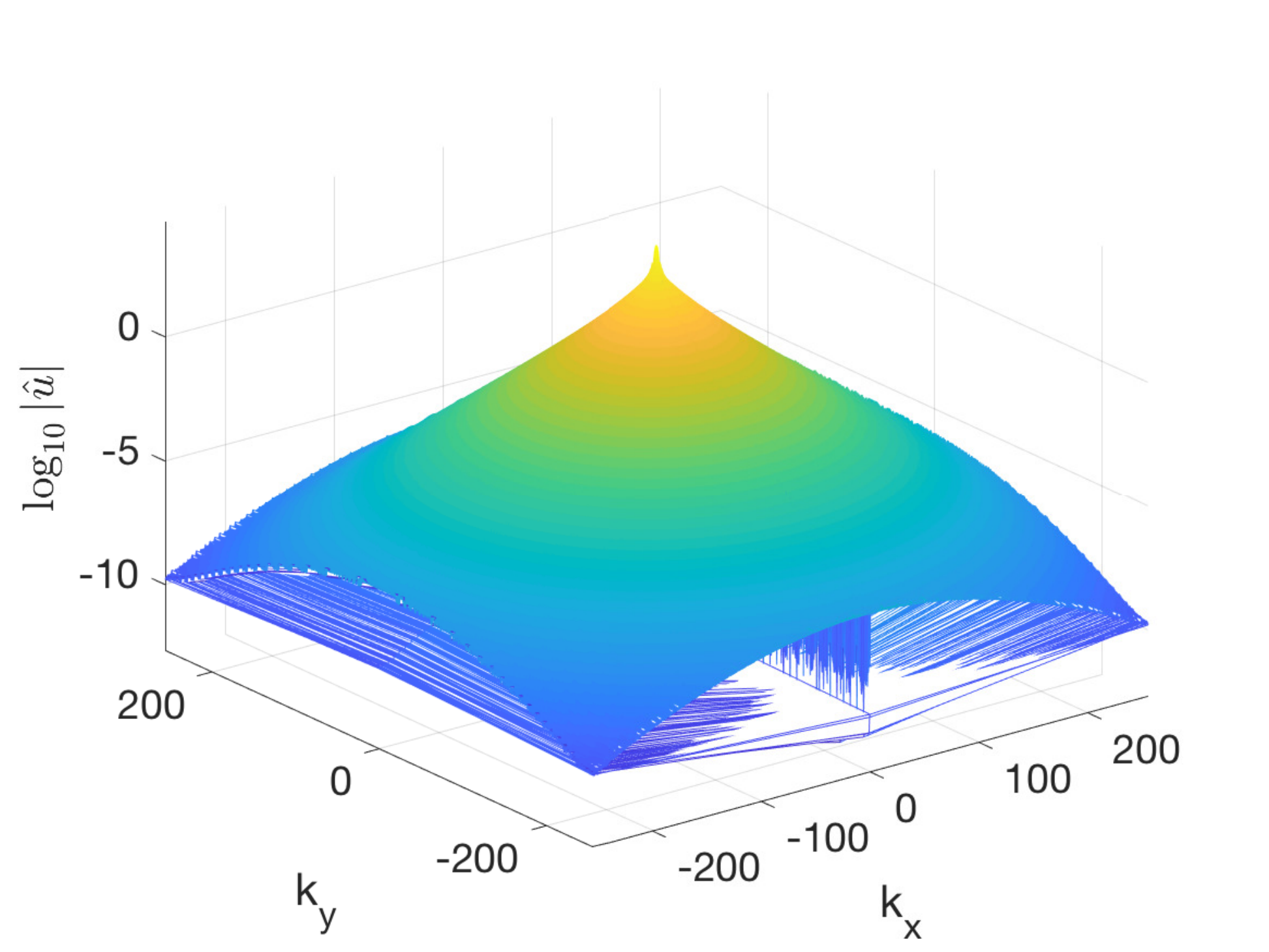}
\caption{Solution to \eqref{ZK}, $p=4$, with $u(x,y,0)=1.1 Q(x,y)$ 
at $t=0.1399$ on the left, and the corresponding Fourier 
coefficients on the right. }
\label{ZKp42gauss}
\end{figure}

For larger $\lambda$, the evolution blows up in finite time. 
For example, for $\lambda=2$, we use the same parameters as for the 
blow-up computation in the previous subsection. The code breaks at 
$t=0.1399$. The ZK solution at that time can be seen in 
Fig.~\ref{ZKp42gauss}.

Fitting of the norms is of the same quality as in 
Fig.~\ref{ZKp4solmax}. We find the fitting to the law \eqref{E:g} for the $L^{\infty}$ norm of $u$ gives the values  $a=-0.229$, $b=0.45$ and $t^{*}=0.1405$ with a fitting error of the order of $10^{-2}$. 
For the $L^{2}$ norm of $u_{x}$, we get  $a=-0.216$, $b=-0.39$ and     
$t^{*}=0.1405$ with a fitting error of the order of $10^{-4}$. 

The good agreement of both norms as well as the agreement with the 
perturbed soliton in the previous subsection gives strong evidence 
for Conjecture \ref{C:3}. Note also the similarity of the 
blow-up profiles in Fig. ~\ref{ZKp4sol11t101674} and Fig. ~\ref{ZKp42gauss} 
suggesting a universal blow-up profile $P(X,Y)$ which, however, we do not investigate it in this paper.

\end{document}